
\def\arXiv{oui}
  
 \def\oui{oui} 
  
\ifx\arXiv\oui
\else
 \pdfpagewidth=210truemm
 \pdfpageheight=297truemm 
\fi
  
  %
  %

  %
  %

  \catcode`@=12 

 \def\defrefnote#1{\definexref{#1}{{\the\footnotenumber}}{refnotes}}

  %
  %


\ifx\couleurs\oui
\input graphicx
 \pdfpagewidth=210truemm
 \pdfpageheight=297truemm 
 \voffset=-5mm
\fi

\input eplain.tex
\expandafter\def\expandafter\newdimen\expandafter{\newdimen}

\ifx\couleurs\oui
\beginpackages
\usepackage{color}
\endpackages 
 \pdfpagewidth=210truemm
 \pdfpageheight=297truemm 
\long\def\rge#1{{\color{red}#1}}

\definecolor{bleu-iecn}{cmyk}{.98,.13,.1,.55}

\else
\long\def\rge#1{#1}

\fi

\makeatletter
\def\numberedfootnote{%
ÊÊ\global\advance\footnotenumber by 1
ÊÊ\@eplainfootnote{{\number\footnotenumber}}%
}%
\def\makecolumns#1/#2 {\par \begingroup
ÊÊ \@columndepth = #1
ÊÊ \advance\@columndepth by -1
ÊÊ \divide \@columndepth by #2
ÊÊ \advance\@columndepth by 1
ÊÊ \@linestogoincolumn = \@columndepth
ÊÊ \@linestogo = #1
ÊÊ \currentcolumn = 1
ÊÊ \def\@endcolumnactions{%
ÊÊÊÊÊÊ\ifnum \@linestogo<2
ÊÊÊÊÊÊÊÊ \the\crtok \egroup \endgroup \par 
ÊÊÊÊÊÊ\else
ÊÊÊÊÊÊÊÊ \global\advance\@linestogo by -1
ÊÊÊÊÊÊÊÊ \ifnum\@linestogoincolumn<2
ÊÊÊÊÊÊÊÊÊÊÊÊ\global\advance\currentcolumn by 1
ÊÊÊÊÊÊÊÊÊÊÊÊ\global\@linestogoincolumn = \@columndepth
ÊÊÊÊÊÊÊÊÊÊÊÊ\the\crtok
ÊÊÊÊÊÊÊÊ \else
ÊÊÊÊÊÊÊÊÊÊÊÊ&\global\advance\@linestogoincolumn by -1
ÊÊÊÊÊÊÊÊ \fi
ÊÊÊÊÊÊ\fi
ÊÊ }%
ÊÊ \makeactive\^^M
ÊÊ \letreturn \@endcolumnactions
ÊÊ \@columnwidth = \hsize
ÊÊÊÊ \advance\@columnwidth by -\parindent
ÊÊÊÊ \divide\@columnwidth by #2
ÊÊ \penalty\abovecolumnspenalty
ÊÊ \noindent 
ÊÊ \valign\bgroup
ÊÊÊÊ &\hbox to \@columnwidth{\strut \hsize = \@columnwidth ##\hfil}\cr
}%
\makeatother

\lefteqnumbers
   \def\testd{oui}
   \def\choixlat{\ifx\numadroite\testd\righteqnumbers
            \else  \lefteqnumbers\fi}
    \choixlat

\catcode`@=\letter
\def\@eplainfootnote#1{\let\@sf\empty 
  \ifhmode\edef\@sf{\spacefactor\the\spacefactor}\/\fi
  \global\advance\hlfootlabelnumber by 1
  \hlstart@impl{foot}{\hlfootlabel}%
  \hldest@impl{footback}{\hlfootbacklabel}%
  \hbox{$^{(#1)}$}%
  \hlend@impl{foot}%
  \@sf\vfootnote{#1.}%
}%
\catcode`@=\other

  \interfootnoteskip=0pt
  \let\note=\numberedfootnote
  \everyfootnote={\eightpoint\leftskip=5truemm\rightskip5truemm}
  
  \hsize150truemm\vsize 240truemm\hoffset=5truemm

  \def\dimart{\hsize126truemm\vsize186truemm\hoffset16truemm\voffset=24truemm}
  \pretolerance=500\tolerance=1000\brokenpenalty=5000
  \parindent3mm
  
  \countdef\temps=170
  \temps=\time
  \countdef\nminutes=171{\nminutes = \time}
  \countdef\nheure=172
  \def\heure{\begingroup                     
     \temps = \time \divide\temps by 60
     \nheure = \temps                        
     \nminutes = \time
     \multiply\temps by 60
     \advance\nminutes by -\temps            
     \ifnum\nminutes<10 \toks1 = {0}%
     \else\toks1 = {}%
     \fi
     \number\nheure h\the\toks1 \number\nminutes  
  \endgroup}%

  \newcount\chstart
  \chstart=\pageno
 \headline={\ifnum\pageno=\chstart {\hfill} \else {\hss \tenrm --\ \folio\ --\hss}\fi}
  \footline={\hfill}
  \normalbaselines
  \frenchspacing
    \def\dater{\vglue-10mm\rightline{(\the\day/\the\month/\the\year)}}
  \def\dateheure{\vglue-10mm\rightline{(\the\day/\the\month/\the\year,\ \heure)}}

  \newif\ifpagetitre \pagetitretrue
\newtoks\hautpagetitre \hautpagetitre={\hfill}
\newtoks\baspagetitre \baspagetitre={\hfill}
\newtoks\auteurcourant \auteurcourant={\hfill}
\newtoks\titrecourant \titrecourant={\hfill}
\newtoks\hautpagegauche
\newtoks\hautpagedroite
\newtoks\hautpagemilieu
\hautpagemilieu={\tenrm\hfil -- \folio\ -- \hfil}
\hautpagegauche={\ifx\midfolio\oui\the\hautpagemilieu\else\tenrm\folio\hfill\the\auteurcourant\hfill\fi}
\hautpagedroite={\ifx\midfolio\oui\the\hautpagemilieu\else\hfill\the\titrecourant\hfill\tenrm\folio\fi}
\newtoks\baspagegauche \baspagegauche={\hfil}
\newtoks\baspagedroite \baspagedroite={\hfil}
\headline={\ifpagetitre\the\hautpagetitre
\else\ifodd\pageno\the\hautpagedroite\else\the\hautpagegauche\fi\fi }
\footline={\ifpagetitre\the\baspagetitre
\else\ifodd\pageno\the\baspagedroite
\else\the\baspagegauche\fi\fi \global\pagetitrefalse}

\def\pageblanche{\vfill\eject\pagetitretrue
\null\vfill\eject
\pagetitretrue
}
\def\chgtpage{\ifodd\pageno \else
\pageblanche \fi \pagetitretrue\titreun=0\footnotenumber=0}

\def\chgtpageincrtitreun{\ifodd\pageno \else
\pageblanche \fi \pagetitretrue\footnotenumber=0}

\def\majnombres{\ifodd\pageno \else
\pageblanche \fi \pagetitretrue\hautpoly\titreun=0\footnotenumber=0}

\def\hautspages#1#2{\auteurcourant={\ninepcap#1}\titrecourant={\nineit#2}}

\ifnum\chstart=\pageno \pagetitretrue\fi
  

  \def\PAR{\par}
  
  \def\leftnote#1{\vadjust{\setbox1=\vtop{\hsize 20mm\parindent=0pt\eightpoint
  \baselineskip=9pt\rightskip=4mm plus 4mm\vskip-4.7mm#1}\hbox{\kern-2cm\smash{\box1}}}}

  
  \def\raggedcenter{\leftskip=20pt plus 10em  
       \rightskip=\leftskip 
        \parfillskip=0pt 
         \spaceskip=.3333em \xspaceskip=.5em 
          \pretolerance=9999 \tolerance=9999
           \hyphenpenalty=9999 \exhyphenpenalty=9999 }
           
  \def\titrecentre#1{{\parindent0mm\raggedcenter
       \spaceskip=.6em plus .2em minus .2em\xspaceskip=.6em plus .2em minus .2em
        \tit#1\par}}
        


  \def\oui{oui}
  
\def\fontetitreun{\ifx\paradouze\oui\douzepts\gpdouze\twelvebf\textfont1=\twelveib\else
\quatorzepts\gpquatorze\fourteenbf\fi}

\def\fontetitreunl{\douzepts\textfont1=\twelveib\scriptfont1=\tenib\fourteenti}
 
 \def\fontetitredeux{\textfont1=\eleveni\ifx\paradouze\oui\onzepts\scriptfont1=\ninei\elevenit\else
                        \douzepts\twelveit\fi}
 
   \def\fontetitredeuxb{\ifx\paradouze\oui\onzepts\eleventi\gponze\textfont1=\elevenib\scriptfont1=\nineib
                         \else\douzepts\twelveti\scriptfont1=\twelveib\scriptfont1=\tenib\gpdouze\fi}
                         
\def\fontetitredeuxl{\onzepts\textfont1=\elevenbf\scriptfont1=\ninebf\twelvebf}
  
\def\fontetitretrois{\textfont0=\elevenrm\scriptfont0=\eightrm\textfont1=\eleveni
                      \scriptfont1=\eighti\scriptscriptfont1=\sixi\elevenit}
                      
\def\fontetitrequatre{\textfont0=\elevenrm\scriptfont0=\eightrm\textfont1=\eleveni
                      \scriptfont1=\eighti\scriptscriptfont1=\sixi\elevenrm}
  
  \newcount\titreun\titreun=0
  \newcount\titredeux\titredeux=0
  \newcount\titretrois\titretrois=0
  \newcount\titrequatre\titrequatre=0
  \newcount\enonce\enonce=0
  
  \def\incr#1{\global\advance#1 by 1 {\the #1}}
  \def\avance#1{\global\advance#1 by 1}
  \def\init#1{\global#1=0}
  
  \long\def\Indentation#1#2{\setbox10=\hbox{\fontetitreun#1}
  	                    \ifdim\wd10 < 4mm
                         \setbox10=\hbox to 4mm{\box10\hfill}
                       \else\ifdim\wd10 < 6mm
                         \setbox10=\hbox to 6mm{\box10\hfill}
  	                    \else\ifdim\wd10 < 8mm
                         \setbox10=\hbox to 8mm{\box10\hfill}
                       \else\ifdim\wd10 < 12mm
                         \setbox10=\hbox to 12mm{\box10\hfill}
                       \fi\fi\fi\fi
                       \dimen10=\hsize
                       \advance \dimen10 by -\wd10
                       \noindent \box10 %
                       \ignorespaces
                       \hbox{\vtop{\hsize=\dimen10\raggedright\noindent\fontetitreun#2}}}

  \long\def\paraun#1{\removelastskip\par\medskip\goodbreak\vskip0pt plus.01\vsize\penalty-100
                \vskip0pt plus-.01\vsize
  	              \init{\titredeux}\ifnum\optionparag=1{\init\eqnumber\init\enonce}\else{}\fi
                  \goodbreak{\fontetitreun
  	                \Indentation{\incr{\titreun}.\ }{\fontetitreun #1\par}}\nobreak\medskip}

 %
 %
 \long\def\paraunc#1{\removelastskip\par\bigskip\goodbreak\vskip0pt plus.01\vsize\penalty-100
                \vskip0pt plus-.01\vsize
  	              \init{\titredeux}
                 \ifnum\optionparag=1{\init{\eqnumber}\init\enonce}\else{}\fi
                  \goodbreak
  	                {\parindent0mm\raggedcenter\fontetitreun\incr{\titreun}.\ 
                     \fontetitreun #1\par}\nobreak\medskip}
                     
\newtoks\titreunl
\titreunl={\ifnum\titreun=1{I}\fi%
\ifnum\titreun=2{II}\fi%
\ifnum\titreun=3{III}\fi%
\ifnum\titreun=4{IV}\fi%
\ifnum\titreun=5{V}\fi%
\ifnum\titreun=6{VI}\fi%
\ifnum\titreun=7{VII}\fi%
\ifnum\titreun=8{VIII}\fi%
\ifnum\titreun=9{IX}\fi%
\ifnum\titreun=10{X}\fi%
\ifnum\titreun=11{XI}\fi%
\ifnum\titreun=12{XII}\fi%
\ifnum\titreun=13{XIII}\fi%
}
\long\def\paraunl#1{\removelastskip\par\bigskip\bigskip\goodbreak\vskip0pt plus.01\vsize\penalty-100
                \vskip0pt plus-.01\vsize
  	              \init{\titredeux}\ifnum\optionparag=1{\init\eqnumber\init\enonce}\else{}\fi
                  \goodbreak{\fontetitreunl
  	                \Indentation{\global\advance\titreun by 1{\the\titreunl}.\ }{\fontetitreunl #1\par}}\nobreak\smallskip}

  
  \long\def\paradeux#1{\init{\titretrois}\vskip0pt plus.01\vsize\penalty-10
                \vskip0pt plus-.01\vsize\ifx \elie\oui\medskip\ifnum\titredeux>0\medskip\fi\fi
                 \Indentation{\fontetitredeux\the\titreun${\cdot}$\incr{\titredeux}.}
                              {\fontetitredeux\textfont1=\eleveni#1}\nobreak\par }
  
  \long\def\paradeuxb#1{\init{\titretrois}\vskip0pt plus.001\vsize\penalty-10
                \vskip0pt plus-.01\vsize{\ifx \elie\oui\medskip\ifnum\titredeux>0\medskip\fi\fi
                  \Indentation
  {\fontetitredeuxb\the\titreun${\cdot}$\incr{\titredeux}.}{ \fontetitredeuxb#1}}\nobreak
\smallskip}

\newtoks\titredeuxl
\titredeuxl={\ifnum\titredeux=1{A}\fi%
\ifnum\titredeux=2{B}\fi%
\ifnum\titredeux=3{C}\fi%
\ifnum\titredeux=4{D}\fi%
\ifnum\titredeux=5{E}\fi%
\ifnum\titredeux=6{F}\fi%
\ifnum\titredeux=7{G}\fi%
\ifnum\titredeux=8{H}\fi%
\ifnum\titredeux=9{I}\fi%
\ifnum\titredeux=10{J}\fi%
\ifnum\titredeux=11{K}\fi%
\ifnum\titredeux=12{L}\fi%
\ifnum\titredeux=13{M}\fi%
}
 \long\def\paradeuxl#1{\init{\titretrois}\vskip0pt plus.001\vsize\penalty-10
                \vskip0pt plus-.01
                \vsize \bigskip%
                  \Indentation
     {\fontetitredeuxl\global\advance\titredeux by 1
  \quad \the\titreunl${\cdot}$\the\titredeuxl.}{ \fontetitredeuxl#1}
  \removelastskip\nobreak\smallskip}
  

  \long\def\paratrois#1{\init{\titrequatre}\ifdim\lastskip<\smallskipamount
                \removelastskip\smallskip\fi
                 \vskip0pt plus.01\vsize\penalty-10
                  \vskip0pt
plus-.01\vsize{\ifx \elie\oui\ifnum\titretrois>0\medskip\fi\fi
\Indentation{\fontetitretrois\the\titreun${\cdot}$\the\titredeux${\cdot}$\incr{\titretrois}.\ }
  {\hskip0mm\baselineskip=14pt\fontetitretrois#1}\nobreak\smallskip}}
  
  
  \long\def\paratroisl#1{\init{\titrequatre}\ifdim\lastskip<\smallskipamount
                \removelastskip\fi
                 \vskip0pt plus.01\vsize\penalty-10
                  \vskip0pt
plus-.01\vsize\ifx \elie\oui\bigskip
\fi
\Indentation{\fontetitretrois\quad \quad \the\titreunl{${\cdot}$}\the\titredeuxl${\cdot}$\incr{\titretrois}.\ }
  {\hskip0mm\fontetitretrois#1}\nobreak\smallskip}


  \long\def\paraquatre#1{\ifdim\lastskip<\smallskipamount
                \removelastskip\smallskip\fi
                 \vskip0pt plus.01\vsize\penalty-10
                  \vskip0pt
                  plus-.01\vsize\par
 
\Indentation{\fontetitrequatre \the\titreun{${\cdot}$}\the\titredeux${\cdot}$\the\titretrois${\cdot}$\incr{\titrequatre}.\ }
{\hskip0mm\fontetitrequatre#1}\nobreak\smallskip}


\newtoks\titrequatrel
\titrequatrel={\ifnum\titrequatre=1{a}\fi%
\ifnum\titrequatre=2{b}\fi%
\ifnum\titrequatre=3{c}\fi%
\ifnum\titrequatre=4{d}\fi%
\ifnum\titrequatre=5{e}\fi%
\ifnum\titrequatre=6{f}\fi%
\ifnum\titrequatre=7{g}\fi%
\ifnum\titrequatre=8{h}\fi%
\ifnum\titrequatre=9{i}\fi%
\ifnum\titrequatre=10{j}\fi%
\ifnum\titrequatre=11{k}\fi%
\ifnum\titrequatre=12{l}\fi%
\ifnum\titrequatre=13{m}\fi%
}
\long\def\paraquatrel#1{\ifdim\lastskip<\smallskipamount
                \removelastskip\smallskip\fi
                 \vskip0pt plus.01\vsize\penalty-10
                  \vskip0pt
                  plus-.01\vsize{\bigskip
\Indentation{\global\advance\titrequatre by 1
\fontetitrequatre\quad \quad \quad \the\titreunl${\cdot}$\the\titredeuxl${\cdot}$\the\titretrois${\cdot}$\the\titrequatrel.\ }
{\hskip0mm\fontetitrequatre#1}\nobreak\smallskip}}

\ifx\optionkeys\oui
\def\drefun#1{\definexref{¤#1}{{\the\titreun}}{}} 
\def\drefdeux#1{\definexref{¤#1}{{\the\titreun}.{\the\titredeux}}{}}
\def\dreftrois#1{\definexref{¤#1}{{\the\titreun}.{\the\titredeux}.{\the\titretrois}}{}}
\else
\def\drefun#1{\definexref{prg#1}{{\the\titreun}}{}} 
\def\drefdeux#1{\definexref{prg#1}{{\the\titreun}.{\the\titredeux}}{}}
\def\dreftrois#1{\definexref{prg#1}{{\the\titreun}.{\the\titredeux}.{\the\titretrois}}{}}
\fi

%


  \long\def\propdeux#1#2#3#4{%
       \avance{\enonce}
       \leavevmode\edef\temp{#2}%
         \ifx\temp\empty 
          \else
           \definexref{#2}{#1~{\the\titreun.\the\enonce}}{enonces}
            \definexref{s#2}{{\the\titreun.\the\enonce}}{enonces}
             \fi
\smallskip
      \noindent{\bf#1\ {\bf\the\titreun.\the\enonce{#3}.}\enspace}{\sl#4\par}%
      \ifdim\lastskip<\medskipamount \removelastskip\penalty55\par \fi
   }

  \long\def\propun#1#2#3#4{%
      \avance{\enonce}
       \leavevmode\edef\temp{#2}%
        \ifx\temp\empty 
          \else
           \definexref{#2}{#1~{\the\enonce}}{enonces}
            \definexref{{s#2}}{{\the\enonce}}{enonces}
             \fi
   \par 
     \noindent{\bf#1\ {\bf\the\enonce{#3}.}\enspace}{\sl#4\par}%
     \ifdim\lastskip<\medskipamount \removelastskip\penalty55\medskip\fi
  }
  
  \long\def\prop#1#2#3#4{\ifnum\optionparag=1
                          \propdeux{#1}{#2}{\textfont1=\elevenib#3}{#4} \else\propun{#1}{#2}{\textfont1=\elevenib#3}{#4}\fi}

  \long\def\propt#1#2#3{\ifx\tpf\oui \prop{Th\'eo\-r\`eme}{#1}{#2}{#3}\par
                       \else\prop{Theorem}{#1}{#2}{#3}\par\fi}
  \long\def\Propt#1#2{\propt{#1}{}{#2}}
  \long\def\propl#1#2#3{\ifx\tpf\oui\prop{Lem\-me}{#1}{#2}{#3}\par
                         \else\prop{Lemma}{#1}{#2}{#3}\par\fi}
  \long\def\Propl#1#2{\propl{#1}{}{#2}}
  \long\def\propc#1#2#3{\ifx\tpf\oui\prop{Corol\-laire}{#1}{#2}{#3}\par
                         \else\prop{Corollary}{#1}{#2}{#3}\par\fi}
  \long\def\Propc#1#2{\propc{#1}{}{#2}}

  \long\def\propd#1#2#3{\ifx\tpf\oui\prop{D\'efi\-nition}{#1}{#2}{#3}\par
                       \else\prop{Definition}{#1}{#2}{#3}\par\fi} 
  
  \long\def\proptd#1#2#3{\ifx\tpf\oui\prop{Th\'eor\`eme et d\'efi\-nition}{#1}{#2}{#3}\par
                       \else\prop{Theorem and definition}{#1}{#2}{#3}\par\fi}


  
  \newcount\optionparag\optionparag=1
  
  \long\def\section#1#2{\ifnum\optionparag=1 \paraun{#2} 
                        \else\goodbreak{\fontetitreun
  	                \Indentation{#1.\ }{#2}}\nobreak\smallskip\fi}

  \def\eqconstruct#1{\ifnum\optionparag=1{\the\titreun\hbox{$\cdot$}#1}\else{#1}\fi}

  
  
  \def\numref{oui}  
  
  \newcount\mesref\mesref=0 
  \def\defbib#1{\ifx\numref\oui\global\advance\mesref by 1 \definexref{#1}{{\the
                 \mesref}}{}\else\definexref{#1}{#1}{}\fi}
  \def\bibtem#1{\defbib{#1}\item{\citer{#1}}}
  \def\citer#1{[\ref{#1}]}

  
  \font\seventeenmsa=msam10 at 17pt    
  \font\fourteenmsa=msam10 at 14pt
  \font\twelvemsa=msam10 at 12pt
  \font\tenmsa=msam10                 
  \font\ninemsa=msam10 at 9pt 
  \font\eightmsa=msam10 at 8pt 
  \font\sevenmsa=msam7 
  \font\sixmsa=msam10 at 6pt
  \font\fivemsa=msam5
  \newfam\msafam\textfont\msafam=\tenmsa\scriptfont\msafam=\sevenmsa\scriptscriptfont\msafam=\fivemsa
  
  \font\seventeenbb=msbm10 at 17pt     
  \font\fourteenbb=msbm10 at 14pt
  \font\twelvebb=msbm10 at 12pt
  \font\tenbb=msbm10                   
  \font\ninebb=msbm10 at 9pt 
  \font\eightbb=msbm10 at 8pt 
  \font\sevenbb=msbm7 
  \font\sixbb=msbm10 at 6pt
  \font\fivebb=msbm5 
  \newfam\bbfam\textfont\bbfam=\tenbb\scriptfont\bbfam=\sevenbb\scriptscriptfont\bbfam=\fivebb
  \def\bb{\fam\bbfam\tenbb}%

  \font\seventeenscaln=eusm10 at 17pt   
  \font\twelvescaln=eusm10 at 12pt
  \font\tenscaln=eusm10                
  \font\ninescaln=eusm10 scaled 900
  \font\eightscaln=eusm10 scaled 800
  \font\sevenscaln=eusm10 scaled 700
  \font\sixscaln=eusm10 scaled 600
   
  \newfam\scalnfam\textfont\scalnfam=\tenscaln\scriptfont\scalnfam=\sevenscaln\scriptscriptfont\scalnfam=\sixscaln
  \def\scaln{\fam\scalnfam\tenscaln}%
  \def\scal{\scaln}
  
  \font\tenscalb=eusb10                

  \font\sevenscalb=eusb10 scaled 700

  \newfam\scalbfam\textfont\scalbfam=\tenscalb\scriptfont\scalbfam=\sevenscalb
  %
  
  %
  %
  \font\fourteenrm=cmr12 scaled 1200
  \font\elevenrm=cmr10 at 11pt
  \font\twelverm=cmr12
  \font\ninerm=cmr9
  \font\eightrm=cmr8      
  \font\sevenrm=cmr7
  \font\sixrm=cmr6

  \font\seventeenpcap=cmcsc10 at 17pt
  \font\tenpcap=cmcsc10                        
  \font\ninepcap=cmcsc9
  \font\eightpcap=cmcsc8
  \font\sevenpcap=cmcsc10 scaled 700
  
  \newfam\pcapfam\textfont\pcapfam=\tenpcap\scriptfont\pcapfam=\sevenpcap
  \def\pcap{\fam\pcapfam\tenpcap}
  
  \font\seventeenrm=cmbx12 scaled 1400

  \font\fourteenbf=cmbx10 scaled 1400
  
  \font\twelvebf=cmbx12
  \font\elevenbf=cmbx10 at 11pt
  \font\ninebf=cmbx9  
  \font\eightbf=cmbx8
  \font\sixbf=cmbx6
  
  \font\tengot=eufm10                           
   
  \font\eightgot=eufm10 at 8truept 
  \font\sevengot=eufm7 
  \font\sixgot=eufm10 at 6 truept 
   
  \newfam\gotfam
  \textfont\gotfam=\tengot\scriptfont\gotfam=\sevengot\scriptscriptfont\gotfam=\sixgot
  \def\got{\fam\gotfam\tengot}%

  
  \def\tit{%
  \textfont0=\seventeenrm\scriptfont0=\tenrm\def\rm{\fam0\seventeenrm}%
  \textfont1=\seventeenib\scriptfont1=\twelveib%
  \textfont2=\seventeensy\scriptfont2=\twelvesy\scriptscriptfont2=\ninesy
  \textfont3=\seventeenex
  \textfont\itfam=\seventeenti
  \def\it{\fam\itfam\seventeenti}%
  \textfont\bbfam=\seventeenbb \scriptfont\bbfam=\twelvebb
  \def\bb{\fam\bbfam\seventeenbb}%
  \textfont\msafam=\seventeenmsa\scriptfont\msafam=\twelvemsa
  \textfont\scalnfam=\seventeenscaln
  \def\pcap{\fam\pcapfam\seventeenpcap}
  \normalbaselineskip=25pt\normalbaselines\rm}

  \font\seventeenti=cmbxti10 scaled 1680
  
  \font\fourteenti=cmbxti10 at 14pt
  
  \font\twelveti=cmbxti10 scaled 1200
  \font\eleventi=cmbxti10 at 11pt

  %
  %
  \font\twelveit=cmti12	
  \font\elevenit=cmti10 scaled 1100
  \font\nineit=cmti9
  \font\eightit=cmti8
  \font\sevenit=cmti7

  %
  %
 
 \font\seventeenib=cmmib10 scaled 1680
  \font\fourteenib=cmmib10 scaled 1400
  \font\twelveib=cmmib10 scaled 1200
  \font\elevenib=cmmib10 scaled 1100
  \font\tenib=cmmib10
\font\eightib=cmmib10 scaled 800
  \font\nineib=cmmib10 scaled 900
\font\sevenib=cmmib10 scaled 700
\font\sixib=cmmib10 scaled 600
\font\fiveib=cmmib10 scaled 500

\ifx\ITAN\oui
\else
\innernewfam\cmmibfam
\textfont\cmmibfam=\tenib
\scriptfont\cmmibfam=\sevenib
\scriptscriptfont\cmmibfam=\fiveib
\def\ib{\fam\cmmibfam\tenib}
\fi

  %
  %
  \font\twelvei=cmmi10 scaled 1200
  \font\eleveni=cmmi10 scaled 1100
  \font\ninei=cmmi9
  \font\eighti=cmmi8 
  \font\seveni=cmmi7 			                
  \font\sixi=cmmi6
  
  \font\ninesl=cmsl9                    
  \font\eightsl=cmsl8 
  \font\sevensl=cmsl10 at 7pt

  \font\ninett=cmtt9                    
  \font\eighttt=cmtt8
  \font\seventt=cmtt10 scaled 700

  \font\seventeensy=cmsy10 scaled 1680    
  \font\fourteensy=cmsy10 scaled 1400
  \font\twelvesy=cmsy10 scaled 1176
  
  \font\ninesy=cmsy9                      
  \font\eightsy=cmsy8
  \font\sixsy=cmsy6
  \font\seventeenex=cmex10 at 17pt
  \font\fourteenex=cmex10 at 14pt
  \font\twelveex=cmex10 at 12pt
  \font\nineex=cmex10 at 9pt
  \font\eightex=cmex10 at 8pt
  \font\sevenex=cmex10 at 7pt
  \font\sixex=cmex10 at 6pt
  \font\fiveex=cmex10 at 5pt
  
   
  \font\fourteengp=cmmi10 at 14pt
  
  \font\twelvegp=cmmib10 at 12pt
  \font\elevengp=cmmib10 at 11pt
  \font\tengp=cmmib10                          
  \font\ninegp=cmmib10 at 9pt 
  \font\eightgp=cmmib8 
   
  \font\sixgp=cmmib6


  \def\gponze{\textfont0=\elevenbf\scriptfont0=\eightbf\scriptscriptfont0=\sixbf
           \textfont1=\elevengp\scriptfont1=\eightgp\scriptscriptfont1=\sixgp}
  \def\gpdouze{\textfont0=\twelvebf\scriptfont0=\tenbf\scriptscriptfont0=\ninebf
           \textfont1=\twelvegp\scriptfont1=\tengp\scriptscriptfont1=\ninegp}        
  
 \def\gpquatorze{\textfont0=\fourteenbf\scriptfont0=\twelvebf\scriptscriptfont0=\elevenbf
           \textfont1=\fourteengp\scriptfont1=\twelvegp\scriptscriptfont1=\elevengp}

  
  \expandafter\chardef\csname pre amssym.def at\endcsname=\the\catcode`\@
  \catcode`\@=11
  \def\undefine#1{\let#1\undefined}
  \def\newsymbol#1#2#3#4#5{\let\next@\relax
   \ifnum#2=\@ne\let\next@\msafam@\else
   \ifnum#2=\tw@\let\next@\bbfam@\fi\fi
   \mathchardef#1="#3\next@#4#5}
  \def\mathhexbox@#1#2#3{\relax
   \ifmmode\mathpalette{}{\m@th\mathchar"#1#2#3}%
   \else\leavevmode\hbox{$\m@th\mathchar"#1#2#3$}\fi}
  \def\hexnumber@#1{\ifcase#1 0\or 1\or 2\or 3\or 4\or 5\or 6\or 7\or 8\or
   9\or A\or B\or C\or D\or E\or F\fi}
  
  \def\setboxz@h{\setbox\z@\hbox}
  \def\wdz@{\wd\z@}
  \def\boxz@{\box\z@}
  
  \edef\msafam@{\hexnumber@\msafam}
  \mathchardef\dabar@"0\msafam@39
  
  \edef\bbfam@{\hexnumber@\bbfam}
  \def\widehat#1{\setboxz@h{$\m@th#1$}%
   \ifdim\wdz@>\tw@ em\mathaccent"0\bbfam@5B{#1}%
   \else\mathaccent"0362{#1}\fi}
  \def\widetilde#1{\setboxz@h{$\m@th#1$}%
   \ifdim\wdz@>\tw@ em\mathaccent"0\bbfam@5D{#1}%
   \else\mathaccent"0365{#1}\fi}
  \newsymbol\leqq 1335          
  \newsymbol\leqslant 1336
  \newsymbol\lessgtr 1337       
  \newsymbol\backprime 1038     
  \newsymbol\risingdotseq 133A  
  \newsymbol\fallingdotseq 133B 
  \newsymbol\succcurlyeq 133C   
  \newsymbol\geqq 133D          
  \newsymbol\geqslant 133E
  \newsymbol\nmid 232D
  \newsymbol\nexists 2040
  \newsymbol\smallsetminus 2272
  \newsymbol\varnothing 203F 
  \catcode`\@=\active

  \catcode`\@=11
  \newcount\typofr\typofr=1
  
  \catcode`\;=\active
  \def;{\ifnum\typofr=1\relax\ifhmode\ifdim\lastskip>\z@\unskip\fi
     \kern.2em\fi\string;\else\string;\fi}
  
  \catcode`\:=\active
  \def:{\ifnum\typofr=1\relax\ifhmode\ifdim\lastskip>\z@\unskip\fi
  \penalty\@M\ \fi\string:\else\string:\fi}
  
  \catcode`\!=\active
  \def!{\ifnum\typofr=1\relax\ifhmode\ifdim\lastskip>\z@\unskip\fi
     \kern.2em\fi\string!\else\string!\fi}
  
  \catcode`\?=\active
  \def?{\ifnum\typofr=1\relax\ifhmode\ifdim\lastskip>\z@\unskip\fi
     \kern.2em\fi\string?\else\string?\fi}

  \def\francais{\typofr=1\def\tpf{oui}}
  
  \def\oui{oui}
  \francais
  
  \catcode`\@=12
  

\ifx\textures\oui
\def\raye #1|{\leavevmode\setbox1=\hbox{#1}%
\raise .5pt\hbox to \wd1{\xleaders\hbox{\rge{$ \sct / $}%
\kern 1pt}\hfill\kern -1pt }\kern -\wd1 \unhbox1\relax }

\def\barre #1|{\leavevmode\setbox1=\hbox{#1}%
\rlap{\Red\vrule height 2.4pt depth -1.2pt width \wd1}\Black \unhbox1\relax}
\else
\def\raye #1|{\leavevmode\setbox1=\hbox{#1}%
\raise .5pt\hbox to \wd1{\xleaders\hbox{\rge{$ \sct / $}%
\kern 1pt}\hfill\kern -1pt }\kern -\wd1 \unhbox1\relax }

\def\barre #1|{\leavevmode\setbox1=\hbox{#1}%
\rlap{\color{red}\vrule height 2.4pt depth -1.2pt width \wd1}\color{black} \unhbox1\relax}

\fi
  

  
  \def\og{\leavevmode\raise.24ex\hbox{$\scriptscriptstyle\langle\!\langle\>$}}    
  \def\fg{\leavevmode\raise.24ex\hbox{$\scriptscriptstyle\>\rangle\!\rangle$}}    

  \def\d{\,{\rm d}}

  \def\z{{\bb Z}}
  \def\r{{\bb R}}
  
  \def\N{{\bb N}}
  
  \def\F{{\bb F}}

  \def\A{{\scal A}}
  
  \def\C{{\scaln C}}
  
  \def\E{{\scal E}}

  \let\LL=\L
  \def\L{{\scal L}}
  \def\M{{\scal M}}
  
  \def\O{{\scal O}}
  \def\P{{\scaln P}}

  \def\W{{\scal W}}

  \def\frac#1#2{{#1\over #2}}
  \def\di#1#2{\sct#1\atop{\sct#2}}

  \def\page#1{\rm p.\thinspace#1}
  \def\theoreme#1{\rm th.\thinspace#1}
  
  \def\qedbox{$\rlap{$\sqcap$}\sqcup$}           
  \def\qed{\nobreak\hfill\penalty250 \hbox{}\nobreak\hfill\qedbox\par }
  \def\Ssi{si, et seulement si, }

  \def\fa{fonc\-tion arith\-m\'etique}

  \def\¤{\S\thinspace}

  \def\¥{$\bullet$ }
  
  
  \def\e{{\rm e}}
  
  \def\md#1#2{\equiv#1\,({\rm mod\,}#2)}

  \def\epsilon{\varepsilon}

  \def\phi{\varphi}
  \def\theta{\vartheta}
  \def\rho{\varrho}
  \def\dm{{\textstyle{1\over 2}}}
  \def\txt{\textstyle}
  \def\dsp{\displaystyle}
  \def\sct{\scriptstyle}
  \def\pf{\noi{\it Proof. }}
  \def\nid{\ifnum\typofr=1\par\noindent{\it D\'emonstration. }\else\pf\fi}
  \def\noi{\noindent}
  \def\rem{\ifnum\typofr=1\noi{\it Remarque.}\ \else\noi{\it Remark.}\ \fi}
  \def\rems{\ifnum\typofr=1\noi{\it Remarques.}\ \else\noi{\it Remarks.}\ \fi}
  \def\re{{\Re e\,}}

  \def\emptyset{{\varnothing}}
  
  \def\pp{{\rm pp}}

  \def\1{{\bf 1}}
  \def\|{\Vert}

  \def\le{\leqslant}\def\leq{\leqslant}
  \def\ge{\geqslant}\def\geq{\geqslant}
  \def\wh{\widehat}
  \def\cf{{cf.}}
  \def\ie{{i.e.\ }}
  \def\eg{{e.g.}}
  \newsymbol\subsetneqq 2324
  \newsymbol\subsetneq 2328

  \def\log{\mathop{\rm log}\nolimits}
  \def\ft#1#2{{\txt{#1\over #2}}}




  \def\pmb#1{\setbox0=\hbox{#1}%
  \kern-.025em\copy0\kern-\wd0\kern.05em\copy0\kern-\wd0\kern-.025em\raise .0433em\box0 }

  
  \skewchar\eighti='177 \skewchar\sixi='177
  \skewchar\eightsy='60 \skewchar\sixsy='60
  
  \def\eightpoint{%
  \textfont0=\eightrm\scriptfont0=\sixrm\scriptscriptfont0=\fiverm
  \def\rm{\fam0\eightrm}%
  \textfont1=\eighti\scriptfont1=\sixi
  \scriptscriptfont1=\fivei\def\oldstyle{\fam1\seveni}%
  \textfont2=\eightsy\scriptfont2=\sixsy\scriptscriptfont2=\fivesy
  \textfont3=\eightex\scriptfont3=\sixex
  \textfont\itfam=\eightit
  \def\it{\fam\itfam\eightit}%
  \textfont\slfam=\eightsl
  \def\sl{\fam\slfam\eightsl}%
  \textfont\bbfam=\eightbb \scriptfont\bbfam=\sixbb\scriptscriptfont\bbfam=\fivebb
  \def\bb{\fam\bbfam\eightbb}%
  \textfont\msafam=\eightmsa\scriptfont\msafam=\sixmsa
  \textfont\scalnfam=\eightscaln
  \def\scaln{\fam\scalnfam\eightscaln}
  \textfont\ttfam=\eighttt
  \def\tt{\fam\ttfam\eighttt}%
\textfont\gotfam=\eightgot
  \textfont\bffam=\eightbf\scriptfont\bffam=\sixbf\scriptscriptfont\bffam=\fivebf
  \def\bf{\fam\bffam\eightbf}%
  \ifx\ITAN\oui\else\textfont\cmmibfam=\eightib
       \scriptfont\cmmibfam=\sixib
        \scriptscriptfont\cmmibfam=\fiveib
         \def\ib{\fam\cmmibfam\eightib}
   \fi
  \textfont\pcapfam=\eightpcap
  \def\pcap{\fam\pcapfam\eightpcap}
  \abovedisplayskip=2pt plus2pt minus 2pt
  \belowdisplayskip=2pt plus1pt minus 2pt
  \abovedisplayshortskip= 1pt plus 2pt minus 1pt
  \belowdisplayshortskip= 1pt plus 2pt minus 1pt
  \smallskipamount=2pt plus 1pt minus 2pt
  \medskipamount=3pt plus 2pt minus 2pt
  \bigskipamount=7pt plus 3pt minus 3pt
  \setbox\strutbox=\hbox{\vrule height 5pt depth 2pt width 0pt}%
  \normalbaselineskip=9pt\normalbaselines\rm}

  \def\({\left(}
  \def\){\right)}
  
  \def\footnoterule{\kern -2pt\hrule width 7truecm\kern 2.4pt}
  
  \def\xnotedef#1{\definexref{#1}{\noexpand\number\footnotenumber}{Note}}%

  
  
  \def\ninepoint{%
  \textfont0=\ninerm\scriptfont0=\sixrm\scriptscriptfont0=\fiverm
  \def\rm{\fam0\ninerm}%
  \textfont1=\ninei\scriptfont1=\sixi
  \scriptscriptfont1=\fivei\def\oldstyle{\fam1\ninei}%
  \textfont2=\ninesy\scriptfont2=\sixsy\scriptscriptfont2=\fivesy
  \textfont3=\nineex\scriptfont3=\sixex
  \textfont\itfam=\nineit
  \def\it{\fam\itfam\nineit}%
  \textfont\slfam=\ninesl
  \def\sl{\fam\slfam\ninesl}%
  \textfont\bbfam=\ninebb\scriptfont\bbfam=\sixbb\scriptscriptfont\bbfam=\fivebb
  \def\bb{\fam\bbfam\ninebb}%
  \textfont\msafam=\ninemsa\scriptfont\msafam=\sixmsa\scriptscriptfont\msafam=\fivemsa
  \textfont\scalnfam=\ninescaln
  \def\scaln{\fam\scalnfam\ninescaln}
  \textfont\ttfam=\ninett
  \def\tt{\fam\ttfam\ninett}%
  \textfont\bffam=\ninebf\scriptfont\bffam=\sixbf\scriptscriptfont\bffam=\fivebf
  \def\bf{\fam\bffam\ninebf}%
  \abovedisplayskip=3pt plus2pt minus 2pt
  \belowdisplayskip=3pt plus1pt minus 2pt
  \abovedisplayshortskip= 2pt plus 2pt minus 1pt
  \belowdisplayshortskip= 2pt plus 2pt minus 1pt
  \smallskipamount=2pt plus 1pt minus 2pt
  \medskipamount=3pt plus 2pt minus 2pt
  \bigskipamount=7pt plus 3pt minus 3pt
  \setbox\strutbox=\hbox{\vrule height 5pt depth 2pt width 0pt}%
  \normalbaselineskip=11pt plus.3pt minus.2pt\normalbaselines\rm}

  \def\sevenpoint{%
  \textfont0=\sevenrm\scriptfont0=\sixrm\scriptscriptfont0=\fiverm
  \def\rm{\fam0\sevenrm}%
  \textfont1=\seveni\scriptfont1=\sixi
  \scriptscriptfont1=\fivei\def\oldstyle{\fam1\seveni}%
  \textfont2=\sevensy\scriptfont2=\sixsy\scriptscriptfont2=\fivesy
  \textfont3=\sevenex\scriptfont3=\fiveex
  \textfont\itfam=\sevenit
  \def\it{\fam\itfam\sevenit}%
  \textfont\slfam=\sevensl
  \def\sl{\fam\slfam\sevensl}%
  \textfont\bbfam=\sevenbb \scriptfont\bbfam=\sixbb\scriptscriptfont\bbfam=\fivebb
  \def\bb{\fam\bbfam\sevenbb}%
  \textfont\msafam=\sevenmsa\scriptfont\msafam=\sixmsa
  \textfont\scalnfam=\sevenscaln
  \def\scaln{\fam\scalnfam\sevenscaln}
  \textfont\bffam=\sevenbf\scriptfont\bffam=\sixbf\scriptscriptfont\bffam=\fivebf
  \def\bf{\fam\bffam\sevenbf}%
  \textfont\ttfam=\seventt
  \abovedisplayskip=2pt plus2pt minus 2pt
  \belowdisplayskip=2pt plus1pt minus 2pt
  \abovedisplayshortskip= 1pt plus 2pt minus 1pt
  \belowdisplayshortskip= 1pt plus 2pt minus 1pt
  \smallskipamount=2pt plus 1pt minus 2pt
  \medskipamount=3pt plus 2pt minus 2pt
  \bigskipamount=7pt plus 3pt minus 3pt
  \setbox\strutbox=\hbox{\vrule height 5pt depth 2pt width 0pt}%
  \normalbaselineskip=9pt\normalbaselines\rm}

 \def\onzepts{%
 \textfont0=\elevenrm\scriptfont0=\ninerm
 \textfont1=\eleveni\scriptfont1=\ninei
}

\def\douzepts{%
  \textfont0=\twelverm\scriptfont0=\tenrm\def\rm{\fam0\twelverm}%
  \textfont1=\twelvei\scriptfont1=\teni%
  \textfont2=\twelvesy\scriptfont2=\tensy\scriptscriptfont2=\eightsy
  \textfont3=\twelveex
  \textfont\itfam=\twelveti
  \def\it{\fam\itfam\twelveti}%
  \textfont\bffam=\twelvebf\scriptfont\bffam=\tenbf\scriptscriptfont\bffam=\eightbf
  \def\bf{\fam\bffam\twelvebf}%
  \textfont\bbfam=\twelvebb \scriptfont\bbfam=\tenbb
  \def\bb{\fam\bbfam\twelvebb}%
  \textfont\msafam=\twelvemsa\scriptfont\msafam=\tenmsa
  \textfont\scalnfam=\twelvescaln
  \normalbaselineskip=15pt\normalbaselines\rm}

\def\quatorzepts{%
  \textfont0=\fourteenrm\scriptfont0=\twelverm\def\rm{\fam0\fourteenrm}%
  \textfont1=\fourteenib\scriptfont1=\twelveib%
  \textfont2=\fourteensy\scriptfont2=\twelvesy\scriptscriptfont2=\tensy
  \textfont3=\fourteenex
  \textfont\itfam=\fourteenti
  \def\it{\fam\itfam\fourteenti}%
  \textfont\bffam=\fourteenbf\scriptfont\bffam=\twelvebf\scriptscriptfont\bffam=\tenbf
  \def\bf{\fam\bffam\fourteenbf}%
  \textfont\bbfam=\fourteenbb \scriptfont\bbfam=\twelvebb
  \def\bb{\fam\bbfam\fourteenbb}%
  \textfont\msafam=\fourteenmsa\scriptfont\msafam=\twelvemsa
  \textfont\scalnfam=\twelvescaln
  \normalbaselineskip=18pt\normalbaselines\rm}


\def\AA{{\it Acta Arith.}}

\def\PLMS{{\it Proc. London Math. Soc.}}

\def\picture #1 by #2 (#3){\leavevmode\vbox to #2{
     \hrule width #1 height 0pt depth 0pt
      \vfill
       \special{picture #3}}}

\def\illustration #1 by #2 (#3) scaled #4{\dimen1=#2
  \divide\dimen1 by 1000
  \multiply\dimen1 by #4
  \vtop to \dimen1{\dimen1=#1
  \divide\dimen1 by 1000
  \multiply\dimen1 by #4
  \hsize=\dimen1\vss
  \noindent\special{illustration #3 scaled #4}}}

\ifx\couleurs\oui

\fi

\ifx\optionkeymacros\undefined\else \fi

\catcode`\Œ=\active\defŒ{{\aa}}       
\catcode`\º=\active\defº{\int}        
\catcode`\=\active\def{\c c}        
\catcode`\¶=\active\def¶{\partial}    
\catcode`\Ä=\active\defÄ{\oint}       
\catcode`\Æ=\active\defÆ{\triangle}   
\catcode`\Â=\active\defÂ{\neg}        
\catcode`\µ=\active\defµ{\mu}         
\catcode`\¿=\active\def¿{{\o}}        
\catcode`\¹=\active\def¹{\pi}         
\catcode`\Ï=\active\defÏ{{\oe}}       
\catcode`\§=\active\def§{{\ss}}       
\catcode`\ =\active\def {\dagger}     
\catcode`\Ã=\active\defÃ{\sqrt}       
\catcode`\·=\active\def·{\Sigma}      
\catcode`\Å=\active\defÅ{\approx}     
\catcode`\½=\active\def½{\Omega}      
\catcode`\£=\active\def£{{\it\$}}     
\catcode`\°=\active\def°{\infty}      
\catcode`\¤=\active\def¤{{\S}}        
\catcode`\¦=\active\def¦{{\P}}        
\catcode`\¥=\active\def¥{\bullet}     
\catcode`\»=\active\def»{\leavevmode\raise.585ex\hbox{\b a}}      
\catcode`\¼=\active\def¼{\leavevmode\raise.6ex\hbox{\b o}}        
\catcode`\­=\active\def­{\not=}       
\catcode`\²=\active\def²{\leq}        
\catcode`\³=\active\def³{\geq}        
\catcode`\Ö=\active\defÖ{\div}        
\catcode`\É=\active\defÉ{{\dots}}     
\catcode`\¾=\active\def¾{{\ae}}       
\catcode`\Ç=\active\defÇ{\og}         
\catcode`\Ò=\active\defÒ{``}          
\catcode`\Á=\active\defÁ{!`}          
\catcode`\¢=\active\def¢{\rlap/c}     
\catcode`\Ô=\active\defÔ{`}           
\catcode`\Õ=\active\defÕ{'}           


\catcode`\=\active\def{{\AA}}       
\catcode`\'=\active\def'{\c C}        
\catcode`\¯=\active\def¯{{\O}}        
\catcode`\¸=\active\def¸{\Pi}         
\catcode`\Î=\active\defÎ{{\OE}}       
\catcode`\®=\active\def®{{\AE}}       
\catcode`\×=\active\def×{\diamond}    
\catcode`\¡=\active\def¡{\accent'27}  
\catcode`\Ó=\active\defÓ{''}          
\catcode`\±=\active\def±{\pm}         
\catcode`\È=\active\defÈ{\fg}         
\catcode`\À=\active\defÀ{?`}          
\catcode`\Ð=\active\defÐ{--}          
\catcode`\Ñ=\active\defÑ{---}         


\catcode`\Š=\active\defŠ{\"a}        
\catcode`\'=\active\def'{\"e}        
\catcode`\•=\active\def•{\"{\i}}     
\catcode`\š=\active\defš{\"o}        
\catcode`\Ÿ=\active\defŸ{\"u}        
\catcode`\Ø=\active\defØ{\"y}        
\catcode`\å=\active\defå{\^A}        
\catcode`\€=\active\def€{\"A}        
\catcode`\…=\active\def…{\"O}        
\catcode`\†=\active\def†{\"U}        
\catcode`\‡=\active\def‡{\'a}        
\catcode`\Ž=\active\defŽ{\'e}        
\catcode`\'=\active\def'{\'{\i}}     
\catcode`\—=\active\def—{\'o}        
\catcode`\œ=\active\defœ{\'u}        
\catcode`\ƒ=\active\defƒ{\'E}        
\catcode`\æ=\active\defæ{\^E}        
\catcode`\é=\active\defé{\`E}        %
\catcode`\ˆ=\active\defˆ{\`a}        
\catcode`\=\active\def{\`e}        
\catcode`\"=\active\def"{\`{\i}}     
\catcode`\˜=\active\def˜{\`o}        
\catcode`\=\active\def{\`u}        
\catcode`\Ë=\active\defË{\`A}        
\catcode`\‹=\active\def‹{\~a}        
\catcode`\–=\active\def–{\~n}        
\catcode`\›=\active\def›{\~o}        
\catcode`\Ì=\active\defÌ{\~A}        
\catcode`\"=\active\def"{\~N}        
\catcode`\Í=\active\defÍ{\~O}        
\catcode`\‰=\active\def‰{\^a}        
\catcode`\=\active\def{\^e}        
\catcode`\"=\active\def"{\^{\i}}     
\catcode`\™=\active\def™{\^o}        
\catcode`\ž=\active\defž{\^u}        

\let\optionkeymacros\null

\ifx\corr\oui
\input CrayolaColors
\def\rge#1{\Red#1\Black}\else\def\rge#1{#1}\fi

\optionparag=1
\def\paradouze{oui}
\def\eqnumm#1#2{$({\rm I}{\cdot}#1{\cdot}#2)$}

\def\bull{\noi{$\bullet$}}

\dimart

\def\E{{\scal E}}
\def\F{{\scal F}}

\def\fa{fonction arith\-mŽ\-tique}

\def\\{\cr}
\def\mathcal{\scal}
\def\mathbb{\bb}

\def\LL{{\scal L}}

\def\RR{{\scal R}}

\def\paraunn#1{\paraun{#1}\writetocentry{section}{#1}}
\def\paradeuxn#1{\paradeuxb{#1}\writetocentry{subsection}{#1}}

\def\ak{\alpha_\kappa}
\def\jk{j_\kappa}
\def\lk{\lambda_\kappa}
\def\xk{\xi_\kappa(u)}
\def\psxy{\psi_f^*(x,y)}
\def\MM{{\got M}}
\def\L{{\got L}}
\def\nuf{\nu_{f,y}}
\def\deltaf{\delta_{f,y}}
\def\Deltaf{\Delta_{f,y}}

%
%
\newcount\paras\paras=0
\newcount\sparas\sparas=0
\def\tocsectionentry#1#2{\init{\sparas}\avance\paras
{\quad\bf\the\paras\quad }{\hskip-2mm #1\dotfill\hskip3mm\rm#2}\par }%
\def\tocsubsectionentry#1#2{\avance\sparas{\qquad\eightpoint\the\paras.\the\sparas}
{\hskip-2mm\eightpoint#1\dotfill\hskip3mm\rm#2}\par }%

\def\titrart{Moyennes de certaines~fonctions
multiplicatives~sur les entiers friables, 3}
\def\auteurs{G\'erald Tenenbaum \& Jie Wu}
\hautspages{\auteurs}{\titrart}

{\it\obeylines
Compositio Math. \bf144 \rm(2008) 339Ð376.}

\titrecentre{\titrart}
\bigskip
\centerline{\auteurs}
\bigskip
\bigskip
{\eightpoint\leftskip1cm\rightskip1cm
\noi{\bf Abstract.} We consider logarithmic 
averages, over friable integers, of non-negative 
multiplicative
functions. Under logarithmic, one-sided or 
two-sided hypotheses, we obtain sharp estimates 
that improve upon
known results in the literature regarding both 
the quality of the error term and the range of 
validity. The
one-sided hypotheses correspond to classical 
sieve assumptions. They are applied to provide an 
effective form of
the Johnsen--Selberg  prime power sieve.
\par\medskip
\noi{\bf AMS Subject classification.} 11N25, 11N35, 11N36, 11N56, 11C08.\par
\medskip
\noi{\bf Keywords.} Friable integers, Selberg's sieve, logarithmic averages of
multiplicative functions.\par}
\bigskip
\noi{\qquad\quad \bf Sommaire}

\smallskip
{\eightpoint\leftskip.6cm\rightskip.8cm
\readtocfile
}

\bigskip\medskip
\paraunn{Introduction}
Ce travail est le troisi\`eme volet d'une s\'erie 
consacr\'ee ˆ l'Žtude des valeurs moyennes de
certaines fonctions arithm\'etiques
sur les entiers friables.
\par
Dans la premire partie \citer{TW03},
nous avons
essentiellement examinŽ le cas de fonctions
multiplicatives $f(n)$ pour
   lesquelles les nombres $f(p)$ possdent,
lorsque $p$ dŽcrit la suite des nombres premiers,
une valeur moyenne
strictement positive. Dans la seconde 
\citer{HTW06},
Žcrite en collaboration avec
Guillaume Hanrot, nous avons explorŽ les
consŽquences de renseignements concernant
directement la sŽrie de Dirichlet
associŽe ˆ $f$, supposŽe analytiquement proche
d'un produit de fonctions ztas de Dedekind ---
ce qui a permis de
rel‰cher dans une certaine mesure l'hypothse de multiplicativitŽ.
\par
Nous nous proposons ˆ prŽsent de considŽrer des
moyennes de type logarithmique, {\it a priori}
plus rŽgulires, et
partant susceptibles d'tre contr™lŽes sous des
hypothses plus faibles concernant les nombres
$f(p)$.
Les moyennes de ce type sont importantes pour les
applications, notamment
dans les m\'ethodes de crible: voir par exemple~\citer{HR74} et \citer{G01}.

D\'esignons par $P(n)$
le plus grand facteur premier d'un entier naturel positif $n$,
avec la convention $P(1)=1$, et
par $\N_y:=\{n\in\N^*:P(n)\leqslant y\}$ l'ensemble des entiers $y$-friables.
Pour tous $x\ge 1$, $y\ge 1$, nous posons
$$S(x, y):=\N_y\cap [1, x],\qquad
S^*(x, y):=\N_y\,\cap\,]x, \infty[.$$
\par
Soit $f$ une fonction multiplicative. La sŽrie de
Dirichlet de $f\1_{\N_y}$ s'Žcrit, dans son domaine de convergence,
$$
F(s,y)
:= \sum_{n \in \N_y} {f(n)\over n^s}
= \prod_{p\le y} \sum_{\nu\ge 0} {f(p^\nu)\over p^{\nu s}}\cdot
\eqdef{Fsy}
$$

Nous supposons systŽmatiquement dans la suite que
$f$ est positive ou nulle et que $F(s,y)$
converge en $s=1$ --- et
donc sur toute la droite verticale $\re s=1.$ La
rarŽfaction ˆ l'infini des entiers
$y$-friables laisse alors augurer que  la fonction sommatoire pond\'er\'ee
$$\psi_f(x, y)
:= \sum_{n\in S(x, y)} {f(n)\over n}
\eqdef{psi}$$
converge rapidement vers $F(1, y)$.
Nous verrons que c'est effectivement le cas sous
des hypothses assez gŽnŽrales. Il est alors plus
pertinent (et bien entendu plus dŽlicat) d'Žvaluer le terme rŽsiduel
$$\psi_f^*(x, y):= F(1,y)-\psi_f(x,y)
= \sum_{n\in S^*(x, y)} {f(n)\over n}\cdot
\eqdef{psi*}$$
\par
Notons $u:=(\log x)/\log y$. Nos rŽsultats 
principaux consistent, d'une part, ˆ Žtablir, 
sous une hypothse de
majoration unilatŽrale en moyenne des nombres
$f(p)(\log p)/p$, une inŽgalitŽ du type
$$\psi_f^*(x,y)\leqslant 
F(1,y)\lambda_\kappa(u)u^{Cu/\log y} 
\eqdef{thunilatvague}$$
dans un trs vaste domaine en $(x,y)$, et, 
d'autre part, ˆ montrer, sous une hypothse 
d'estimation bilatŽrale
en moyenne des mmes quantitŽs, une formule asymptotique du type
$$\psi_f^*(x,y)= 
F(1,y)\lambda_\kappa(u)\big\{1+O(R)\big\} 
\eqdef{thbilatvague}$$
o $R$ est un terme d'erreur essentiellement 
comparable ˆ $u(\log u)/\log y$. Ici, 
$\lambda_\kappa$ dŽsigne une
fonction ˆ dŽcroissance trs rapide, dŽfinie 
comme la solution d'une Žquation diffŽrentielle 
aux diffŽrences. Nous
renvoyons le lecteur aux ŽnoncŽs respectifs des ThŽormes
\ref{sthmmajoration} et
\ref{sthmasymp}  pour une formulation prŽcise. 
Bornons-nous ici ˆ deux remarques: d'une part,
\eqref{thunilatvague} est valable dans les 
hypothses classiques du crible,  d'autre part, 
l'estimation
\eqref{thbilatvague} Žtablit l'optimalitŽ de 
\eqref{thunilatvague} dans l'intersection des 
domaines de validitŽ.
\par
Comme indiquŽ plus haut, le
comportement asymptotique de la somme pondŽrŽe $\psxy$ est par
nature plus rŽgulier que celui de la somme non pondŽrŽe
$$\Psi_f(x,y):=\sum_{n\in S(x,y)}f(n). $$
Un aspect remarquable de ce phŽnomne concerne,
ainsi qu'il a ŽtŽ soulignŽ dans~\citer{HTW06}, le
terme d'erreur en
$1/(\log  y)^\kappa$ dans la formule
$$\Psi_f(x,y)=C_\kappa(f)x(\log
y)^{\kappa-1}\bigg\{1+O\bigg({1\over (\log
y)^\kappa}+{\log (u+1)\over \log y}\bigg)\bigg\},$$
Žtablie dans \citer{TW03}, pour un domaine
adŽquat en $(x,y)$.\note{La constante
$C_\kappa(f)$ est prŽcisŽment dŽfinie ˆ la
formule \eqref{Ckf} {\it infra}.}
Lorsque
$f(p)$ est en moyenne  proche de $\kappa$, ce terme rŽsiduel n'est
en fait nŽcessaire que lorsque $x/y$ est
relativement petit et, partant, comme nous le
verrons plus loin, ne contribue pas,  lorsque
$\kappa<1$, ˆ un terme d'erreur spŽcifique dans la
formule asymptotique  correspondante pour $\psxy$.
\medskip
La suite de cet article est organisŽe comme suit. 
Le paragraphe \ref{prghist} contient les 
prŽrequis, relatifs
aux solutions d'Žquations diffŽrentielles aux 
diffŽrences, nŽcessaires ˆ l'ŽnoncŽ de nos 
rŽsultats, ainsi qu'une
brve prŽsentation des estimations de 
$\psi_f^*(x,y)$ antŽrieurement disponibles dans 
la littŽrature. Le paragraphe
\ref{prgresul} est dŽvolu ˆ la formulation 
prŽcise de nos ŽnoncŽs ainsi qu'ˆ la description 
d'une extension
possible. Nous y dŽcrivons Žgalement l'incidence 
sur les Žvaluations de $\psi_f^*(x,y)$ de la
connaissance d'hypothses plus fortes, mais 
standard, relatives aux moyennes non pondŽrŽes 
des nombres
$f(p)\log p$. Le paragraphe \ref{prgSelb} concerne l'application des
rŽsultats d'hypothses unilatŽrales au crible ˆ 
puissances de Gallagher--Selberg. Nous donnons un 
ŽnoncŽ gŽnŽral et
dŽtaillons un exemple pratique reprŽsentatif. 
Les paragraphes suivants sont consacrŽs aux 
dŽmonstrations. En
particulier, les Žquations fonctionnelles, qui 
forment le cÏur de la mŽthode, et leur traitement 
via l'Žquation
dite adjointe au problme, sont prŽsentŽs au paragraphe
\ref{prgeqfonc}.
\medskip
\paraunn{Historique}\drefun{hist}
\paradeuxn{Solutions d'Žquations diffŽrentielles aux diffŽrences}
L'ŽnoncŽ des rŽsultats de la bibliographie
nŽcessite quelques brefs rappels relatifs ˆ
certaines fonctions de
la thŽorie du crible dŽfinies comme solutions
d'Žquations diffŽrentielles aux diffŽrences.
\par
Nous d\'esignons par
$\varrho_\kappa$ la solution continue du systme
$$\cases{
\varrho_\kappa(u)=u^{\kappa-1}/\Gamma(\kappa)
& si $0<u\le 1$,
\cr\noalign{\vskip 1mm}
u\varrho_\kappa'(u)+(1-\kappa)\varrho_\kappa(u)+\kappa\varrho_\kappa(u-1)=0
& si $u>1$,
\cr}
\eqdef{rho}$$
o\`u $\Gamma$ dŽsigne la fonction d'Euler.
Ainsi, $\varrho_\kappa$ est, pour tout
$\kappa>0$, la puissance de convolution
fractionnaire de la fonction de Dickman
$\varrho:=\varrho_1$. Pour chaque $\kappa$,
la fonction $\varrho_\kappa$ hŽrite de la
propriŽtŽ de dŽcroissance rapide de $\varrho$. On
a par exemple
$$\varrho_\kappa(u)=(u\log u)^{-u}\e^{O_\kappa(u)}\qquad (u\to\infty).$$
Plus prŽcisŽment, d\'esignons par $\xi(u)$
l'unique solution r\'eelle non nulle
de $\e^\xi = 1 + u\xi$ si $u>0$, $u\not=1$ et posons $\xi(1)=\xi(0)=0$.
D\'efinissons alors
$$\eqalign{
\xi_\kappa(u)
& := \max\{1, \, \xi(u/\kappa)\},
\cr
I(s)
& := \int_0^s {\e^v-1\over v} \d v,
\cr
\sigma_j(u)
& := \kappa I^{(j)}(\xk),
\cr
}\eqdef{xis}$$
de sorte que, d'aprs le lemme III.5.8.1 de
\citer{Te95} et le lemme 4.5 de \citer{Sm91},
$$\eqalign{\xi_\kappa(u)
&= \log u + \log_2u + O\bigg({\log_2u\over \log u}\bigg)
\cr
\sigma_j(u)&=u\bigg\{1+O\bigg({1\over \log u}\bigg)\bigg\}\cr}
\qquad (u\to\infty),
\eqdef{6}$$
o, ici et dans la suite, $\log_k$ d\'esigne la
$k$-i\`eme it\'erŽe de la fonction logarithme.
Il rŽsulte alors du thŽorme 1 de \citer{Sm91} ou
du plus gŽnŽral thŽorme 2 de \citer{HT93} que
l'on~a
$$\varrho_\kappa(u)=\bigg\{1+O\bigg({1\over
u}\bigg)\bigg\}{\e^{\gamma\kappa-u\xi_\kappa(u)+\sigma_0(u)}\over
\sqrt{2\pi\sigma_2(u)}}\quad (u\to\infty),
\eqdef{farhok}$$
o $\gamma$ est la constante d'Euler.
\par
Nous posons encore
$$\lk(u)
:= \e^{-\gamma\kappa} \int_u^\infty\varrho_\kappa(v) \d v,\quad \jk(u)=1-\lk(u)
\qquad(u\ge 0),
\eqdef{lambda}$$
de sorte que $\lk(0)=\jk(\infty)=1$ --- voir par exemple \citer{He86}.
Ainsi qu'il a ŽtŽ notŽ dans~\citer{HTW06} ---  formule (4.12) ---, nous avons
$$\lk(u)=\bigg\{1+O\bigg({1\over  u}\bigg)\bigg\}
{\e^{-\gamma\kappa} \varrho_\kappa(u)\over \xk}
\qquad (u>0).
\eqdef{evallk}$$
Nous observons immŽdiatement, ˆ fins de rŽfŽrence
ultŽrieure, que cela implique, gr‰ce par exemple
aux estimations de
\citer{Sm91},
$$\lambda_\kappa(u-v)\ll\lk(u)\e^{v\xk}\qquad
(u\ge 1,\,0\le v\le u-\dm).\eqdef{comploclam}$$
\goodbreak
\paradeuxn{RŽsultats
antŽrieurs}

Pour $C>0$, $\delta>0$, $\kappa>0$, dŽsignons par
$\M(C,\delta,\kappa)$ la classe des fonctions
multiplicatives r\'eelles positives ou nulles $f$
satisfaisant aux conditions
$$\leqalignno{
& \bigg|\sum_{p\le z} f(p) {\log p\over p}
- \kappa \log z\bigg|\leqslant C(\log z)^{1-\delta}
\qquad
(z\ge 2),
& \eqdef{moyfpSong}
\cr
& \sum_{p}\sum_{\nu\ge 2}{f(p^\nu)\over p^\nu}\log p^\nu
\leqslant C.
& \eqdef{majfpnuSong}
\cr}$$
Sous l'hypothse supplŽmentaire $0<\delta<1$,
Song \citer{So01} a montr\'e que l'on a
$$\psi_f^*(x, y)
= F(1, y)
\bigg\{\lk(u)+O\bigg({\log(u+1)\over (\log y)^\delta}\bigg)\bigg\}
\eqdef{asympSong}$$
uniform\'ement pour
$f\in \M(C,\delta,\kappa),
\,
x\ge y\ge 2,
$
o, ici et dans la suite, nous notons systŽmatiquement
$$u:={\log x\over \log y}
\quad(x\ge 1, \, y>1),
\eqdef{defu}$$

\par
La formule \eqref{asympSong} peut \^etre
consid\'er\'ee comme une version quantitative
d'un r\'esultat dž ˆ de Bruijn \& van Lint \citer{deBL64} Žtablissant,
sous certaines hypothses de mme nature mais plus faibles,
une conclusion de la forme
$$\psi_f(x, y)
\sim (\log y)^\kappa \LL(\log y)
\qquad(y\to\infty)$$
uniform\'ement en $u$ sur tout compact de $]0,\infty[$, o
$\LL$ est une fonction ˆ croissance lente au sens de Karamata.
\par
En raison de la majoration triviale
$$\psi^*_f(x,y)\leqslant F(1,y),\eqdef{majtriv} $$
la formule \eqref{asympSong} est sans intŽrt hors du domaine
$1\le u\le \exp\big\{c(\log y)^\delta\big\}$ o
$c$ est une constante positive arbitraire.
Il n'est donc pas anecdotique de supprimer le
facteur $\log (u+1)$ dans le terme d'erreur.
C'est l'un des objets du rŽsultat de
Tenenbaum dans \citer{Te01}, qui fournit ainsi
une Žvaluation non triviale pour toute valeur de
$x$ et $y$, bien que
de moindre qualitŽ pour les grandes valeurs de 
$u$. Le mme travail contient Žgalement une 
extension au cas de
fonctions ˆ valeurs complexes, sous une hypothse 
plus forte concernant les moyennes des nombres 
$f(p)$.
\medskip
Une Žtude rŽcente due ˆ Greaves \citer{G05}
concerne essentiellement le cas $\delta=1$,
assujetti ˆ certaines
restrictions supplŽmentaires. Plus prŽcisŽment,
Žtant donnŽes deux constantes $A>0$, $\kappa>0$, et notant
    $$r_f(z):=\sum_{p\le z} f(p) {\log p\over p}- \kappa \log z\qquad
(z\ge 1),\eqdef{moyfp} $$
pour chaque fonction multiplicative positive ou nulle  $f$,
Greaves place son Žtude dans les hypothses
$$\leqalignno{
& r_f(z)-r_f(w)\leqslant A
\qquad(z\ge w\ge 1),
& \eqdef{eta1cote}
\cr\noalign{\vskip 2,8mm}
& f(p^\nu)=0
\qquad(p\geqslant 2,\,\nu\ge 2).
& \eqdef{fpnu0}
\cr}$$
   Pour toute constante $B\ge 1$ telle que
$$
{1\over \log z}\sum_{p\le z} {f(p)\log p\over p+f(p)}\leqslant B
\qquad(z\ge 2),\note{Il rŽsulte de
\eqref{eta1cote} que $B=A+\kappa$ est un choix
admissible.}
$$
il introduit alors
$$h_B(u) := \cases{
u\log(u/B)-u+B, & si $u>B$,
\cr\noalign{\smallskip}
0,              & dans le cas contraire.
\cr}$$

Avec ces notations, les rŽsultats principaux de  \citer{G05}
peuvent \^etre \'enonc\'es comme suit.

\proclaim Th\'eor\`eme A (\citer{G05}, {\theoreme 1}).
Soient $A$ et $\kappa$ deux nombres rŽels
positifs. Il existe une constante $c_\kappa$, ne
dŽpendant que de $\kappa$,
telle que, pour toute fonction multiplicative positive ou nulle  $f$
v\'erifiant
\eqref{eta1cote} et \eqref{fpnu0}, on ait
$$\psi_f^*(x, y)
\leqslant F(1, y)
\bigg(\lk(u)
+{\e^{c_\kappa A-h_B(u)}\over \log y}\bigg)\qquad (x\ge y\ge 2).
\eqdef{Greavesmaj}
$$

\proclaim Th\'eor\`eme B (\citer{G05}, {\theoreme 2}).
Soient $A$ et $\kappa$ deux nombres rŽels
positifs. Il existe une constante $c_\kappa$, ne
dŽpendant que de $\kappa$,
telle que, pour toute fonction multiplicative positive ou nulle  $f$
v\'erifiant \eqref{fpnu0} et
$$\leqalignno{
&|r_f(z)|
    \leqslant \dm A
\qquad(z\ge 1),
& \eqdef{eta2cotes}
\cr\noalign{\vskip 1mm}
&\sum_{p} {f(p)^2\log p\over p^2}
    \leqslant A,
& \eqdef{Hfp2log/p2}
\cr}$$
on ait
$$\psi_f^*(x, y)
= F(1, y)\bigg\{\lk(u)+O\bigg({u\,\e^{c_\kappa
A-h_B(u)}\over \log y}\bigg)\bigg\}\qquad (x\ge
y\ge 2).
\eqdef{Greavesasy}$$
La constante implicite ne d\'epend que de $\kappa$.
\par

\smallskip

\rems
(i) Il rŽsulte de \eqref{farhok} et
\eqref{evallk} que la minoration contenue dans la
formule
asymptotique \eqref{Greavesasy} n'est non triviale que si
$$u\ll (\log_2y)/\log_4y,\ \hbox{ \ie }\quad
y>x^{c(\log_4x)/\log_2x}.
\eqdef{domaineGreaves}$$
Dans ce mme domaine restreint,
\eqref{Greavesmaj} fournit donc une majoration
asymptotiquement optimale de
$\psi_f^*(x, y)$.
\smallskip\goodbreak
(ii) Comme l'hypo\-thse~\eqref{fpnu0} implique
que $\psxy=0$ pour $$x>N_y:=\prod_{p\le
y}p=\e^{y+o(y)},\eqdef{xgrandsfc}$$
l'inŽgalitŽ~\eqref{Greavesmaj} et la majoration
contenue dans \eqref{Greavesasy} sont triviales
ds que $x>N_y$.
\smallskip
\par(iii)
Les lemmes 1 et 2 de Rawsthorne \citer{R82},
permettent de dŽduire imm\'ediatement du ThŽorme
B une assertion de mme type que le ThŽorme A, mais lŽgrement plus
faible, ˆ savoir que l'on a, sous les mmes hypothses,
$$\psi_f^*(x, y)
\leqslant F(1, y)
\bigg(\lk(u)
+{u\,\e^{A_1-h_B(u)}\over \log y}\bigg)\qquad (x\ge y\ge 2),
\eqdef{Greavesmaj2}$$
o $A_1$ dŽpend au plus de $A$ et $\kappa$.
\smallskip

\paraunn{R\'esultats}\drefun{resul}
\paradeuxn{Objectifs}\drefdeux{obj}
Au vu des estimations dŽcrites plus haut, il est
naturel d'examiner, sous l'hypothse naturelle
$r_f(z)\ll1$, les possibilitŽs d'extension du 
domaine de validitŽ de la formule asymptotique
$$\psi_f^*(x,y)\sim F(1,y)\lambda_\kappa(u)$$
et, sous une hypothse unilatŽrale comme
\eqref{eta1cote}, de celles de la majoration
asymptotique correspondante
$$\psi_f^*(x,y)\leqslant \{1+o(1)\}F(1,y)\lambda_\kappa(u). \eqdef{inegasymp}$$
En effet, la dŽcroissance rapide de $\lk(u)$ ˆ
l'infini complique singulirement l'obtention
d'un terme
d'erreur relatif de bonne qualitŽ: on notera, par
exemple, que la formule \eqref{Greavesasy} peut
tre essentiellement
rŽŽcrite sous la forme
$$\psxy=F(1,y)\lambda_\kappa(u)\bigg\{1+O\bigg({\e^{u\log_2u+O(u)}\over
\log y}\bigg)\bigg\}.\eqdef{faGreaves}$$
\par
\par
Notre mŽthode Žtant relativement flexible, les 
hypothses concernant $r_f(z)$ peuvent tre 
notablement rel‰chŽes.
Nous prŽcisons sans dŽmonstration au paragraphe
\ref{prgextension} les rŽsultats qui peuvent tre 
obtenus sous l'hypothse plus faible 
\eqref{moyfpSong} ou sous
une hypothse unilatŽrale de mme type.
\par\smallskip
Une seconde motivation de ce travail rŽside dans
une clarification de la nature mme des hypothses nŽcessaires
concernant la fonction
multiplicative
$f$. \par
Alors qu'il appara"t indispensable, en vue
d'Žlargir le champ des applications, de
s'affranchir de la condition
\eqref{fpnu0}, qui restreint le support de $f$ aux
entiers sans facteur carrŽ, il est opportun de
distinguer
et d'Žtablir deux types de thŽormes, selon que
les hypothses portent sur une condition faible,
comme
\eqref{eta1cote} ou \eqref{eta2cotes}, relative aux nombres
$f(p)(\log p)/p$ ou sur une condition forte --- voir
\eqref{regfp} {\it infra} --- concernant le comportement
en moyenne de
$f(p)\log p$. Les ŽnoncŽs de la premire
catŽgorie fourniront des estimations valables
dans un domaine restreint en
$(x,y)$ mais applicables ˆ une vaste catŽgorie de
fonctions $f$, ceux du second, au contraire,
concerneront des
Žvaluations plus prŽcises, valables sous des
conditions plus restrictives pour le choix de $f$.
\medskip
\goodbreak
\paradeuxn{ƒnoncŽs}\drefdeux{enoncesTW}
\'Etant donn\'es des param\`etres
$A>0$,
$C>0$,
$\kappa>0$ et $\eta\in\,]0, \dm[$,
nous introduisons la classe $\MM_\kappa(A, C, \eta)$
--- resp. $\L_\kappa(A, C, \eta)$ ---
des fonctions multiplicatives r\'eelles positives ou nulles $f$
satisfaisant, avec la notation \eqref{moyfp},
aux conditions \eqref{eta2cotes}
--- resp. \eqref{eta1cote} --- et
$$
\sum_{p}\sum_{\nu\ge 2}{f(p^\nu)\over p^{(1-\eta)\nu}}
\leqslant C.
\eqdef{moyfpnu}
$$
Nous avons
$$
\MM_\kappa(A, C, \eta)\subset \L_\kappa(A, C, \eta).
\eqdef{MM+}$$

Nous Žtablissons le rŽsultat suivant.

\Propt{thmmajoration}
{Soient $A>0$, $C>0$, $\kappa>0$ et $\eta\in\,]0, \dm[$.
Pour une constante convenable $B>0$, nous avons
$$\psi_f^*(x, y)
\leqslant F(1, y)\lk(u)\e^{Bu\xk/\log y}
\eqdef{maj1}$$
sous la condition
$$f\in \L_\kappa(A, C, \eta),
\quad
x\ge 3,
\quad(\log x)^{3/\eta}\leqslant y\le x.
\eqdef{domainemajTW1}$$
En particulier, pour tout $c>0$ fixŽ, nous avons
$$\psi_f^*(x, y)
\leqslant F(1, y)\lk(u)\bigg\{1+O\bigg({u\log (u+1)\over \log y}\bigg)\bigg\}
\eqdef{maj1bis} $$
uniformŽment pour
$$f\in\L_\kappa(A,C,\eta),\quad x\geqslant
3,\qquad \exp\Big\{c\sqrt{\log
x\log_2x}\Big\}\leqslant y\leqslant x. $$
\PAR
Sous l'hypoth\`ese suppl\'ementaire \eqref{fpnu0},
l'inŽgalitŽ \eqref{maj1} est valable pour
$x\ge 2,$ $y\ge 2$.}
\goodbreak

Le \ref{thmmajoration} apporte en
particulier les prŽcisions suivantes sur le
ThŽo\-rme~A:
\par
(i)
Dans le sous-domaine \eqref{domainemajTW1},
l'hypoth\`ese contraignante \eqref{fpnu0} peut tre remplacŽe par
l'inoffensive condition  \eqref{moyfpnu} tout en
rŽduisant significativement le majorant.
\par
(ii)
Sous l'hypoth\`ese \eqref{fpnu0}, le  facteur
exponentiel du membre de droite de \eqref{maj1}
peut toujours tre remplacŽ par $ 1+\e^{O(u)}/\log y$: en effet,
comme $\psxy=0$ sous la condition \eqref{xgrandsfc}, nous
pouvons sans perte de gŽnŽralitŽ supposer que $\xk\ll\log y$.
Cela reprŽsente un renforcement consŽquent de \eqref{Greavesmaj}.
\par\goodbreak
(iii)
Le domaine de validit\'e \eqref{domaineGreaves} pour
une inŽgalitŽ asymptotiquement optimale de type
\eqref{inegasymp} a ŽtŽ Žtendu ˆ
$$u=o\bigg({\log y\over \log_2y}\bigg)\qquad (y\to\infty). $$

\goodbreak
\medskip
Introduisons  la notation
$$Z_j(y;f):=1+\sum_{p\le y}{f(p)^2(\log p)^j\over
p^2}\qquad (j=1,\,2). \eqdef{defWj}$$ Il rŽsulte
par exemple de
\eqref{majfp} et \eqref{scfp} {\it infra} que, sous la condition
\eqref{eta1cote}, on a
$$Z_1(y;f)\ll\log_2 y,\quad Z_2(y;f)\ll\log y,$$
alors que
$$Z_1(y;f)+Z_2(y;f)\ll 1, $$
si, par exemple, $f(p)\ll p/(\log p)^{2+\varepsilon}$ avec $\varepsilon>0$.
\smallskip
Sous l'hypo\-thse
bilatŽrale \eqref{eta2cotes}, nous obtenons une
formule asymptotique pour~$\psi^*_f(x,y)$.

\Propt{thmasymp}
{Soient
$A>0$, $C>0$, $\kappa>0$ et $\eta\in\,]0,\dm[$.
Il existe une constante positive $c_1=c_1(A, C, \kappa, \eta)$
telle que l'on ait
$$\psi_f^*(x, y) = F(1, y)\lk(u)
\big\{
1+O\big(E_{x, y}\big)
\big\}
\eqdef{asymTW}$$
uniform\'ement pour
$$
f\in \MM_\kappa(A, C, \eta),
\quad
x\ge 3,
\quad
\exp\Big\{c_1\sqrt{(\log x)\log_2x}\Big\}\leqslant y\le  x,
\eqdef{domaineTW}$$
o\`u l'on a posŽ
$$E_{x, y}
:=\min\Bigg(1,{Z_1(y;f)\over \log y}+
{u\log(u+1)\over \log y}\bigg\{1+u\log
\bigg(1+{Z_2(y;f)\log (u+1)\over \log
y}\bigg)\bigg\}\Bigg).
$$}

\medskip

Le \ref{thmasymp} rŽpond aux objectifs
prŽcŽdemment dŽcrits en renforant le ThŽorme B
dans quatre directions :

(i) La relation  \eqref{asymTW} fournit un Žquivalent asymptotique lorsque
$$x\ge 3,\quad \exp\Big\{w(x)(\log
x)^{2/3}(\log_2x\log_3x)^{1/3}\Big\}\leqslant y\le x $$
o $w(x)$ est une fonction arbitraire tendant
vers l'infini. De plus, ce domaine peut tre
Žtendu ˆ
$$x\ge 3,\quad \exp\Big\{w(x)\sqrt{(\log x)\log_2x}\Big\}\leqslant y\le x $$
ds que $Z_2(y;f)\ll1$, ce qui reprŽsente une
hypothse trs peu contraignante pour les
applications.
Dans tous les cas, le domaine
\eqref{domaineGreaves} a ŽtŽ significativement
agrandi: la borne supŽrieure autorisŽe pour $u$
est passŽe de
$\ll(\log_2y)/\log_4y$ ˆ  une puissance de $\log y$.
\par
(ii) Le terme d'erreur de \eqref{asymTW} est au
plus quadratique en $u\log u$, alors qu'il
cro"t exponentiellement en $u\log_2u$ dans celui de
\eqref{faGreaves}.
\par
(iii)
L'hypothse \eqref{moyfpnu}, modŽrant la
croissance des nombres $f(p^\nu)$ pour $\nu\ge
2$, remplace ˆ prŽsent la drastique condition
   \eqref{fpnu0}.
\par
(iv)
L'hypothse \eqref{Hfp2log/p2}, relative aux
moyennes pondŽrŽes des nombres $f(p)^2$, a ŽtŽ
supprimŽe. Cependant, par le biais d'une majoration de
$Z_2(y;f)$, toute information quantitative de ce
type permet de prŽciser le rŽsultat obtenu.

\medskip
D'un point de vue mŽthodologique, il est
intŽressant de noter qu'on ne peut faire appel
aux rŽsultats de Rawsthorne dans \citer{R82} pour dŽduire du
\ref{thmasymp} mme une version affaiblie du
\ref{thmmajoration}. En
effet, le principe mis en Ïuvre dans \citer{R82},
consistant essentiellement ˆ rŽduire le cas
$f\in\L_\kappa(A,C,\eta)$
au cas $f\in\MM_\kappa(A,C,\eta)$ gr‰ce au calcul des
variations, ne fonctionne, en l'Žtat actuel des
connaissances, que sous
l'hypothse \eqref{fpnu0}.

\smallskip

Nos dŽmonstrations des ThŽormes  \ref{sthmmajoration}  et
\ref{sthmasymp}
reposent sur la m\'ethode itŽrative
mise en Ïuvre dans \citer{Te01}.
Nous traitons directement $\psi_f^*(x, y)$ sans
prŽalablement approcher $f$ par une fonction de
Piltz $\tau_\kappa$,
mais cette modification technique n'est pas 
essentielle. En revanche, nous utilisons de 
manire cruciale une
majoration initiale prŽcise obtenue par la 
mŽthode de Rankin avec un choix essentiellement 
optimal des paramtres.

\smallskip
Ainsi qu'il a ŽtŽ mentionnŽ au paragraphe
\ref{prgobj}, il peut tre utile de disposer
d'Žvaluations de $\psi_f^*(x, y)$
sous des hypoth\`eses plus fortes concernant les
valeurs moyennes des nombres $f(p)$.
Dans cette perspective, la m\'ethode de
\citer{TW03} est directement applicable.
Au prix de quelques modifications mineures, nous
obtenons le rŽsultat suivant, pour l'ŽnoncŽ
duquel nous rappelons ˆ
prŽsent certaines notations.

Par souci de concision, nous renvoyons librement ˆ \citer{TW03}.
\par
Pour $\kappa>0$ et $b\in\,]0,\dm]$, nous
d\'esignons par $\RR(b,\kappa)$ la classe des fonctions croissantes
$R\in\C^1(]1,\infty[,\r^{+*})$  introduite dans
\citer{TW03} --- page 124. Pour la commoditŽ du
lecteur, rappelons
que la plupart des fonctions \og naturelles\fg\
$R(v)$ dont la vitesse de croissance est situŽe
entre celles de
$(\log_2v)^{1+\varepsilon}$ et
$
\exp\{(\log v)^{3/5-\varepsilon}\}$, o
$\varepsilon>0$ est arbitraire, appartiennent ˆ
$\RR(b,\kappa)$.
\par\goodbreak
\'Etant donn\'es des param\`etres
$A>0$, $C>0$, $\kappa>0$, $b\in\,]0,\dm]$,
$\eta\in\,]0,\dm[$, et une fonction
$R\in \RR(b,\kappa)$, nous dŽfinissons, comme dans \citer{TW03}, la classe
$$\M_\kappa=\M_\kappa(A, C,\eta;R)$$ des fonctions multiplicatives
r\'eelles positives ou nulles
$f$ satisfaisant \`a \eqref{moyfpnu} et
$$\eqalignno{
& \Big|\sum_{p\le z}f(p)\log p-\kappa z\Big|\leqslant Az/R(z)
\qquad (z>1).
& \eqdef{regfp}
\cr}$$
Notons toutefois, que, pour la cohŽrence des
dŽfinitions de $\M_\kappa(A,C,\eta;R)$ et
$\MM_\kappa(A,C,\eta)$, nous
avons permutŽ les r™les jouŽs par les constantes
$A$ et $C$ dans la dŽfinition de
$\M_\kappa(A,C,\eta;R)$.
\par
Pour tous
$b\in\,]0,\dm],$
$\varepsilon>0$, $\kappa>0$ et $R\in \RR(b;\kappa)$,
nous posons
$$U^*(y):=R\big(y^b\big)/\big\{1+\big|\log  R\big(y^b\big)\big|\big\}
\qquad{\rm et}\qquad
Y_\varepsilon^*:=y^{U^*(y)/\varepsilon},
\eqdef{UY}$$
et introduisons le domaine $J_\varepsilon^*(R)$ du plan en $x,\,y$,
d\'efini par la condition
$$2\le y\le x\le Y_\varepsilon^*.
\leqno(J_\varepsilon^*(R))$$

\Propt{thmHF}
{Soient
$A>0$, $C>0$, $\varepsilon>0$, $\kappa>0$,
$b\in\,]0,\dm]$, $\eta\in\,]0,\dm]$, $R\in \RR(b; \kappa)$.
On a
$$\psi_f^*(x, y)
= F(1,y)\lk(u) \{1 + O(\E_{x,y})\}
\eqdef{asymHF}$$
uniform\'ement pour $f\in\M_\kappa(A, C, \eta; R)$ et
$2\le y\le x\le Y_\varepsilon^*$,
o\`u l'on a pos\'e
$$\dsp\E_{x,y}
= {u\log(u+1)\over R(y^{b/2})}
+ {\log(u+1)\over \log y}\int_2^y{\d t\over tR(t^b)}
+ {\{u\log(u+1)\}^2\over R(y^b)\log y}.$$}

Pour illustrer les champs d'application
spŽcifiques des ThŽormes \ref{sthmasymp} et
\ref{sthmHF},
consid\'erons trois exemples.
\par
(i)
Lorsque $R(z)=\exp\{(\log z)^{(3-\varepsilon)/5}\}$ avec $0<\varepsilon<\dm$,
   la formule asymptotique~\eqref{asymHF} a lieu
uniform\'ement pour
$$x\ge 3,\quad \exp\big\{(\log_2x)^{5/3+\varepsilon}\big\}\leqslant y\le x
\leqno(H_\varepsilon)$$
et l'on peut choisir
$$\dsp\E_{x,y}:= {\log(u+1)\over \log y}.$$
Nous obtenons ainsi un renforcement considŽrable
de \eqref{asymTW} sous une hypothse beaucoup
plus restrictive.
\par\goodbreak
(ii)
Lorsque $R(z)=(\log z)^\delta$ avec $\delta>0$, il existe une constante $c>0$
telle que la formule asymptotique \eqref{asymHF} ait lieu
uniform\'ement pour
$$x\ge 3,
\qquad
\exp\big\{c(\log x\log_2x)^{1/D}\big\}\leqslant y\le x
\eqdef{domD}$$
avec
$D:=\min\{1+\delta,\dm(3+\delta)\}$ et
$$\dsp\E_{x,y}
= \cases{
\displaystyle
{u\log(u+1)\over (\log y)^\delta}
& si $\delta<1$,
\cr\noalign{\smallskip}
\displaystyle
{(u+\log_2y)\log(u+1)\over \log y}
& si $\delta=1$,
\cr\noalign{\smallskip}
\displaystyle
{\{u\log(u+1)\}^2\over (\log y)^{1+\delta}}
& si $\delta>1$.
\cr}$$
Notons que
l'hypoth\`ese \eqref{regfp} implique, avec la notation \eqref{moyfp},
$$r_f(z) \ll \cases{
(\log z)^{1-\delta}& si $\delta<1$,
\cr\noalign{\vskip 1mm}
\log_2z & si $\delta=1$,
\cr\noalign{\vskip 1mm}
1 & si $\delta>1$.
\cr}$$
Ainsi, lorsque $\delta>1$, le \ref{thmHF}
fournit-il en toute circonstance, moyennant une
hypothse lŽgrement plus
forte, un terme d'erreur relatif  sensiblement plus petit que celui du
\ref{thmasymp}.
\par\smallskip
(iii) Choisissons ˆ prŽsent $R(z):=(\log z)(\log_2z)^{1+c}$, avec $c>0$,
de manire ˆ approcher encore davantage
l'hypothse
$r_f(z)\ll1$ par une condition de
type~\eqref{regfp}. Nous obtenons alors la
validitŽ de \eqref{asymHF} dans le domaine
$$x\ge 3,\qquad \exp\big\{(\log 
x)^{1/2}(\log_2x)^{-c/2}\big\}\leqslant y\le x $$
avec
$$\E_{x,y}:={\log (u+1)\over \log
y}+{u\log (u+1)\over (\log
y)(\log_2y)^{1+c}}\bigg\{1+{u\log (u+1)\over \log
y}\bigg\}\cdot $$ Cela met en
Žvidence une cohŽrence opportune entre les
ThŽormes~\ref{sthmHF} et \ref{sthmasymp}. Il est
en effet ˆ noter que notre
hypothse sur les nombres $f(p)$ implique ˆ prŽsent
$Z_2(y;f)\ll1$.
\medskip
\paradeuxn{Description sommaire d'une extension}
\drefdeux{extension}
Si, dans les applications usuelles, l'on est 
effectivement amenŽ ˆ comparer la fonction 
$r_f(z)$ ˆ une
constante,  notre mŽthode s'adapte sans 
difficultŽ ˆ d'autres types d'hypothses. Dans 
cette
perspective, une Žtude exhaustive, parallle
ˆ celle effectuŽe dans \citer{TW03} pour le cas 
des moyennes standard, pourrait tre envisagŽe. 
Les majorants
des quantitŽs
$r_f(z)$ ou $r_f(z)-r_f(w)$ seraient alors 
susceptibles de varier dans une classe de 
fonctions relativement large,
dŽfinie par des conditions de croissance peu 
restrictives. Nous nous contentons ici de dŽcrire 
brivement les
rŽsultats obtenus dans le cas de majorants  de type \eqref{moyfpSong}.
\par\goodbreak
Introduisons donc, pour chaque valeur du 
paramtre $\delta$ dans $]0,1]$, la sur-classe
$\L_\kappa(A,C,\delta,\eta)$ de
$\L_\kappa(A,C,\eta)$ obtenue en affaiblissant la condition \eqref{eta1cote} en
$$r_f(z)-r_f(w)\leqslant A(\log z)^{1-\delta}
\qquad(z\ge w\ge 1)
\eqdef{eta1cotelog}
$$
et la sur-classe $\MM_\kappa(A,C,\delta,\eta)$ de 
$\MM_\kappa(A,C,\eta)$ dŽfinie en remplaant
\eqref{eta2cotes} par
$$ |r_f(z)|\leqslant \dm A(\log z)^{1-\delta}
\qquad(z\ge 1).
\eqdef{eta2coteslog}
$$
\par\goodbreak
Nous obtenons alors que, pour tous $A>0$, $C>0$, 
$\delta\in\,]0, 1]$, $\kappa>0$  et $\eta\in\,]0,
\dm[$, la majoration
$$
\psi_f^*(x, y)
\leqslant F(1, y)\lk(u)\e^{Bu\xk/(\log y)^\delta}
\eqdef{maj1log}$$
est valable, avec une constante convenable $B>0$,
sous la condition
$$f\in \L_\kappa(A, C, \delta, \eta),
\quad
x\ge 3,
\quad(\log x)^{3/\eta}\leqslant y\le x.
\eqdef{domainemajTW1log}$$
\par
Semblablement, nous pouvons Žtablir, pour tous 
$A,C,\delta,\kappa,\eta$ fixŽs comme indiquŽ plus 
haut,
l'existence d'une constante
$c_1$ telle que l'on ait
$$\psi_f^*(x, y) = F(1, y)\lk(u)
\big\{
1+O\big(E_{x, y}^*\big)
\big\}
\eqdef{asymTWlog}$$
uniform\'ement pour
$$
f\in \MM_\kappa(A, C, \delta, \eta),
\quad
x\ge 3,
\quad
\exp\big\{c_1(\log x\log_2x)^{1/(1+\delta)}\big\}\leqslant y\le  x,
\eqdef{domaineTWlog}$$
o\`u l'on a posŽ
$$
E_{x, y}^*
:=\min\bigg(1,
{Z_1(y;f)\over \log y}
+ {u\log(u+1)\over (\log y)^\delta}
\bigg\{1
+ u\log\bigg(1+{Z_2(y;f)\log (u+1)\over (\log 
y)^{2-\delta}}\bigg)\bigg\}\bigg).
$$

\medskip
\paraunn{Application au crible ˆ puissances de Johnsen--Selberg}\drefun{Selb}
Pour chaque puissance de nombre premier $p^\nu$ $(\nu\ge 1)$,
soit $\W(p^\nu)$
un ensemble de r\'esidus modulo $p^\nu$, 
identifiŽ ˆ la classe de ses reprŽsentants dans 
$\z$.
Supposons que $\W(p^\mu)\cap \W(p^\nu)=\emptyset$
si $\mu\not=\nu$.
Posons alors
$$\W(d):=\bigcap_{p^\nu\| d}\W(p^\nu),
\eqdef{h1}$$
de sorte que $n\in \W(d)$
\Ssi $n\in \W(p^\nu)$
ds que $p^\nu\| d$.
Convenons Žgalement de poser $\W(1)=\z$.
\par
Soient $\A$ une suite finie d'entiers (non nŽcessairement distincts),
$\P$ un ensemble de nombres premiers et
$z\ge 2$ un nombre r\'eel.
Posons
$$\P_z:=\P\cap[1,z].
$$
Nous cherchons une majoration de la quantit\'e
$$S(\A, \P; z)
:= \big|\big\{a\in \A : a\notin \W(p^\nu)
\;\;(p\in\P_z,\,\nu\ge 1)
\big\}\big|.$$
Le cas particulier o\`u ${\cal A}$ est un intervalle
a \'et\'e initialement considŽrŽ par Johnsen \citer{Jo71} sous des
hypothses sensiblement plus restrictives que
celles que nous avons effectuŽes
plus haut. Celles-ci apparaissent dans le travail
\citer{Sel77}, o Selberg Žtend sa mŽthode au
crible par des
puissances de nombres premiers.
Gallagher \citer{Gal72}, \citer{Gal73/74}, a fourni une preuve plus simple 
de l'estimation de Johnsen. Motohashi a
ensuite obtenu dans
\citer{Moto79} une nouvelle d\'emonstration, via
le grand crible, du r\'esultat de~\citer{Sel77},
rŽpondant ainsi ˆ une question posŽe dans le mme travail.
\par
La majoration finale pour $S(\A,\P;z)$ est
explicitŽe dans \citer{Sel77} lorsque $\A$ est
un intervalle. Nous nous
proposons, dans un premier temps, d'Žtendre la
formulation au cas gŽnŽral. Supposons que l'on
puisse
Žcrire
$$\sum_{\scriptstyle a\in \A
\atop\scriptstyle a\in \W(d)} 1
= {w(d)\over d}X+r_d\qquad (d\geqslant 1),
\eqdef{h2}$$
o\`u $X\ge 0$ est une quantit\'e ind\'ependante de $d$,
$w$ est une fonction multiplicative positive ou nulle,
et $r_d$ est, en un sens convenable, un terme
d'erreur. Nous pouvons supposer sans perte de
gŽnŽralitŽ que
$$w(p^\nu)=0 \qquad (p\notin\P_z,\,\nu\geqslant 1)\eqdef{h3}$$
et que
$$\sum_{\nu\geqslant 1}{w(p^\nu)\over
p^\nu}<1\quad (p\in\P_z), \eqdef{condcrible}$$
ce qui\rge{, dans les cas usuels d'application,} Žquivaut au fait que
$n\notin\cup_{\nu\geqslant 1}\W(p^\nu)$ a lieu
pour au moins un entier $n$.
\par
Soit $D>1$. Pour toute suite r\'eelle $\{\lambda_d\}$ v\'erifiant
$\lambda_1=1$ et $\lambda_d=0$ pour $d>D$,
on a
$$S(\A, \P; z)
\leqslant \sum_{a\in \A}
\Big(\sum_{a\in \W(d)}
\lambda_d\Big)^2.
\eqdef{f1}$$
En effet, le terme d'indice $a$ dans la somme de
\eqref{f1} est toujours positif ou nul (puisque
les $\lambda_d$ sont
rŽels) et il vaut $\lambda_1=1$ si
$a$ est comptŽ dans
$S(\A,\P;z)$.
\par\goodbreak
Introduisons la fonction multiplicative
(au sens de Selberg dans \citer{Sel77})
$$\varepsilon(p^{\mu}, p^{\nu})
:=\cases{
1 & si $\mu=\nu$ ou $\mu\nu=0$,
\cr\noalign{\vskip 1mm}
0 & dans le cas contraire.
\cr}$$
La condition initiale de disjonction des classes
$\W(p^\nu)$ implique que l'on ne peut avoir
$\W(p^\mu)\cap \W(p^\nu)\neq\emptyset$ que si
$\varepsilon(p^\mu,p^\nu)=1$. Par
multiplicativitŽ, il s'ensuit que
$\varepsilon(d,d')=1$ ds que $\W(d)\cap \W(d')\neq\emptyset$. En
d\'eveloppant le carr\'e de
\eqref{f1}, nous obtenons donc
$$\eqalign{S(\A, \P; z)
& \leqslant \sum_{d,\,d'\leqslant D}
\lambda_{d}\lambda_{d'}
\sum_{\di{a\in \A}
{a\in \W(d)\cap \W(d')}} 1
\cr
& = \sum_{d,\,d'\leqslant D}
\lambda_{d}\lambda_{d'} \varepsilon(d, d')
\sum_{\di{a\in \A}
{a\in \W([d, d'])}} 1\cr
&= X \sum_{d,\,d'\leqslant D}
\lambda_{d}\lambda_{d'}\varepsilon(d, d')
{w([d, d'])\over [d, d']}
+ R,\cr}
\eqdef{f2}$$
o\`u
$$R:=\sum_{d,\,d'\leqslant D}
\lambda_{d}\lambda_{d'}\varepsilon(d, d')r_{[d,d']}.
\eqdef{R}$$

Le calcul de Selberg (\citer{Sel77}, \pp.
239--240) pour minimiser le terme principal du
membre de droite de
\eqref{f2}  est valable sans changement. Posons
$$\theta(p^\nu) :=1-\sum_{1\le\mu\le
\nu}{w(p^\mu)\over p^\mu}>0\qquad (p\ge 2,\nu\ge
0). $$
Le
minimum est atteint lorsque l'on~a
$\lambda_d=\lambda_d^*$ pour tout
$d$, avec
$$\lambda_d^*:=
{
\sum_{m\le D} {t^*(m, d)t^*(m, 1)/g(m)}
\over
\sum_{m\le D} {t^*(m, 1)^2/ g(m)}}
\quad(d\le D),
\qquad
\lambda_d^*:=0
\quad(d>D),
\eqdef{lambdad*}
$$
o\`u $g$ et $t^*$ sont  les fonctions
multiplicatives d'une et deux variables
respectivement d\'efinies par
$$\eqalign{
g(p^\nu)
& := \{\theta(p^{\nu-1})-\theta(p^\nu)\}
{\theta(p^\nu)\over \theta(p^{\nu-1})},
\cr
t^*(p^{\mu}, p^{\nu})
& := \cases{
-\{\theta(p^{\mu-1})-\theta(p^\mu)\}/\theta(p^{\mu-1})
& si $\nu=0$,
\cr\noalign{\vskip 1mm}
\{\theta(p^{\mu-1})-\theta(p^\mu)\}/\theta(p^{\mu-1})
& si $0<\nu<\mu$,
\cr\noalign{\vskip 1mm}
1 & si $\mu=\nu$,
\cr\noalign{\vskip 1mm}
0 & si $\nu>\mu$.
\cr}
\cr}$$
Nous observons que $t^*(m,1)=0$ si $P(m)>z$. Un
calcul \'el\'ementaire montre que
$$\sum_{d\leqslant D,\,d'\leqslant D}
\lambda_{d}^*\lambda_{d'}^*\varepsilon(d, d')
{w([d, d'])\over [d, d']}
= {1\over
\sum_{\di{ d\le D}{ P(d)\leqslant z}}
\prod_{p^\nu\| d}\{1/\theta(p^\nu)-1/\theta(p^{\nu-1})\}}\cdot
$$
D'autre part, comme $\sup_d|\lambda_d^*|= 1$, nous pouvons Žcrire
$$\eqalign{
|R|
& \leqslant \sum_{\di{m\le D^2}{ P(m)\leqslant z}}
\sum_{[d, d']=m} \varepsilon(d, d')
|r_m| = \sum_{\di{m\le D^2}{ P(m)\leqslant z}}
3^{\omega(m)}|r_m|,
\cr}$$
o\`u $\omega(m)$ est
le nombre de diviseurs premiers distincts de $m$.

\smallskip

Nous avons ainsi obtenu le r\'esultat suivant.

\propt{CribleT}{ (Selberg)}
{Sous les hypoth\`eses \eqref{h2}, \eqref{h3} et \eqref{condcrible}, nous avons
$$S(\A, \P; z)
\leqslant {X\over \psi_{f}(D, z)}
+ \sum_{\di{m\le D^2}{ P(m)\leqslant z}}
3^{\omega(m)} |r_m|,
$$
o\`u $f$ est la fonction multiplicative d\'efinie par
$f(p^\nu):=p^\nu/\theta(p^\nu)-p^\nu/\theta(p^{\nu-1})$.
}
\goodbreak
Effectuons l'hypothse supplŽmentaire qu'il
existe une constante positive $\eta\in]0,\dm[$ 
\rge{et un entier $s\geqslant 1$ tels}
que
$$\eqalignno{
&\theta(p^\nu)
\geqslant \eta
\quad(p\in \P, \; \nu\ge 1),&\eqdef{h4}\cr
&\rge{\sum_{p}\sum_{1\leqslant \nu\leqslant s}{w(p^\nu)^2\log p\over
p^{2\nu}}+}\sum_{p}\sum_{\nu\rge{>s}}{w(p^\nu)\over 
p^{(1-\eta)\nu}}<\infty.&\eqdef{wpn2}\cr} $$
Nous pouvons alors utiliser les rŽsultats du
paragraphe \ref{prgenoncesTW} pour minorer
$\psi_f(D,z)$, en utilisant la dŽcroissance en 
$z$ de $S(\A,\P;z)$ pour traiter les trs petites 
valeurs de $z$.
\rge{Notant $f_s$  la fonction multiplicative dŽfinie par
$$f_s(p^\nu):=\cases{\dsp\sum_{1\leqslant 
j\leqslant s}{f(p^j)\over p^{j-1}}& si $\nu=1$,\cr
0& si $1<\nu\leqslant s$,\cr
f(p^\nu)& si $\nu>r$,\cr} $$
on vŽrifie aisŽment que l'on a
$$\psi_f^*(x,y)\leqslant \psi_{f_s}^*(x^{1/s},y) 
\qquad (x\geqslant 1, \,y\geqslant 1).$$}
L'hypothse
$\rge{f_s}\in\L_\kappa(A,C,\eta)$ est \rge{alors} satisfaite pour des
constantes convenables $\kappa>0$, $A$, $C$,
$\eta\in\,]0,\dm[$ ds que
$$\rge{\sum_{1\leqslant \nu\leqslant 
s}\sum_{y<p\leqslant t}{w(p^\nu)\log p\over
p^\nu}}\leqslant\kappa\log (t/y)+O(1)\qquad
(t\geqslant y\geqslant 1).
\eqdef{hypcrible}
$$
Posons
$$W(z):=\prod_{p\leqslant z}
\bigg(1-\sum_{\nu\geqslant 1}{w(p^\nu)\over p^\nu}\bigg),
$$
de sorte que, avec la notation \eqref{Fsy} pour
la sŽrie de Dirichlet friable associŽe ˆ la
fonction $f$ dŽfinie plus
haut, nous avons
$W(z)=1/ F(1,z)$.

\Propc{CribleC1}
{Soient $\kappa>0$, $\eta\in\,]0,\dm[$. Sous les
hypoth\`eses \eqref{h2}, \eqref{h3}, \eqref{h4}, \eqref{wpn2},
\eqref{hypcrible}, il existe une constante $B$ telle que l'on ait, uniformŽment
pour $2\le z\le\rge{D^{1/s}}$,
$$\eqalign{S(\A, \P; z)
& \leqslant {XW(z)\over j_\kappa(v)}
\bigg\{1+O\bigg(
{\lambda_\kappa^+(v)\over \log z}\e^{Bv(\log v)/\log z}\bigg)\bigg\}
+ \sum_{\di{m\le D^2}{P(m)\leqslant z}}
3^{\omega(m)} |r_m|
\cr}$$
o\`u nous avons posŽ $$v:=\min\big\{(\log D)/\rge{(s\log
z)},3(\log D)/(\rge s\eta\log_2D)\big\},\quad
\lambda_\kappa^+(v):=\lambda_\kappa(v)v\log (1+v).$$ }

Ce rŽsultat renforce significativement, ds que
$v\to\infty$, la majoration analogue dŽduite du ThŽorme A.

\smallskip

L'ŽnoncŽ  suivant fournit, ˆ titre
d'illustration, un cas concret d'application. 
\rge{Un autre exemple est fourni par le cas des 
entiers d'un
intervalle reprŽsentables comme somme de deux 
carrŽs; les dŽtails seront dŽveloppŽs dans un 
travail ultŽrieur.}
Lorsque
$G\in\z[X]$, nous d\'esignons par $\varrho(d; G)$ le nombre
des solutions dans $\z/d\z$  de la congruence $G(x)\md 0d$.

\Propc{CribleC2}
{Soient $N\in\N^*$, $I$ un intervalle contenant $N$ entiers,
\hbox{$q\in\N^*$}, et $G(X)\in\z[X]$ un polyn™me de degrŽ $g$, dont la
dŽcomposition
canonique en produit de polyn\^omes irr\'eductibles
s'Žcrit $G:=b\prod_{1\leqslant j\leqslant r}G_j$,
o $b\in\z^*$ et les $G_j$ sont distincts et unitaires.
Il existe une constante $B_G>0$, ne dŽpendant que de $G$, telle que l'on ait
$$
\sum_{\di{n\in I}{
p^\nu\| q\Rightarrow p^\nu\,\nmid G(n)}}1
\leqslant {NW(q;G)\over j_r(v)}
\bigg\{1+O\bigg(v\lambda_r(v)v_N\e^{B_Gv_N}+{(\log N)^{4g}\over
N^{\delta_q}}\bigg)\bigg\}\eqdef{estcriblepol}
$$
o\`u l'on a posŽ $\delta_q:={1/\log\{3+\omega(q)\}}$,
$$\eqalign{v
& :=\min\bigg({\log N\over 2\log P(q)},{\log N\over B_G\log_2N}\bigg),
\qquad
W(q;G)
    :=\prod_{p^\nu\| q}
\bigg(1-{\varrho(p^\nu; G)\over p^\nu}\bigg),
\cr}$$
et $\dsp v_N:={v\{\log (v+1)\}/\log N}$.
La constante impliqu\'ee ne d\'epend que de $G$.
}
\goodbreak
\nid Observons d'abord que, si $W(q;G)=0$, il
existe une puissance de nombre premier $p^\nu\|q$
telle que $G$ soit
identiquement nul modulo $p^\nu$. Le membre de
gauche de \eqref{estcriblepol} est donc nul.
\par
Dans toute la suite, nous supposons donc, sans
perte de gŽnŽralitŽ, que \hbox{$W(q;G)\neq0$}.
La minoration
$$W(q;G)\gg\Big({\varphi(q)\over q}\Big)^g \gg(\log \{3+\omega(q)\})^{-g},$$
o $\varphi$ dŽsigne l'indicatrice d'Euler,
rŽsulte alors immŽdiatement du fait que
$\varrho(p^\nu;G)$ est uniformŽment bornŽ et
n'excde pas $g$ sauf peut-tre
pour un nombre fini de valeurs de $p^\nu$.
\par\goodbreak
Appliquons le
\ref{CribleC1}  avec \rge{$s:=1$,}
$
\A:=\{G(n) : n\in I\}$,
$\P:=\{p : p\mid q\}$,  $z:=P(q)$,~et
$$
\W(p^\nu)
:=\cases{
\{0\}  & si $p^\nu\| q$,
\cr\noalign{\vskip 1mm}
\emptyset                  & dans le cas contraire.
\cr}
$$
Convenons d'Žcrire $d\|q$ pour signifier que $d|q$ et $(d,q/d)=1$.
Lorsque cette condition est rŽalisŽe, on a
$$\sum_{\di{a\in \A}{a\in \W(d)}} 1
= \sum_{\di{n\in I}{
G(n)\md0d}} 1
= {\varrho(d; G)\over d}N+O\big(\varrho(d; G)\big),
$$
et
le membre de gauche est nul lorsqu'elle ne l'est pas.
Les hypothses du \ref{CribleC1} sont donc v\'erifi\'ees avec
$$
w(p^\nu):=\cases{\varrho(p^\nu;G)
& si $p^\nu\| q$,
\cr
0 & dans les autres cas.
\cr}$$
Nous notons, en effet, que la validitŽ de la condition \eqref{h4} rŽsulte
   immŽdiatement des propriŽtŽs rappelŽes plus
haut de la fonction
$\varrho(p^\nu;G)$ --- voir, par exemple,  \citer{Te90}, \S2, pour les dŽtails.
Comme on a $$\varrho(p^\nu;G)=\sum_{1\leqslant
j\leqslant r}\varrho(p;G_j)$$ sauf  pour un
nombre fini de valeurs de
$p^\nu$ --- voir, par exemple \citer{Te90},
for\-mu\-le~(2.11)~\hbox{---,}  une application standard du
thŽorme des idŽaux premiers
--- \cf, \eg, \citer{Te90}, lemme~3.1 --- permet
de voir que  l'hypoth\`ese \eqref{hypcrible} est
satisfaite avec~$\kappa=r$.
\par
Choisissons alors $D:=\sqrt{N}$
dans le \ref{CribleC1}. Nous obtenons, dans les
conditions de l'ŽnoncŽ,
$$
\sum_{\di{n\in I}{p^\nu\| q\Rightarrow p^\nu\,\nmid G(n)}}1
\leqslant {NW(q;G)\over j_r(v)}
\bigg\{1+O\bigg({\lambda_r^*(v)\over \log N}+R\bigg)\bigg\}\eqdef{f3}
$$
avec
$$
\lambda_r^*(v):=\lambda_r(v)v^2\log (v+1),\quad R:=
\sum_{\di{d\le N}{d\|q}}
3^{\omega(d)} \varrho(d;G).
$$
Nous estimons $R$ par la mŽthode de Rankin. Pour
tout $\alpha\in\;]0,1[$, nous pouvons Žcrire,
notant $p_j$ le $j$-ime
nombre premier,
$$\eqalign{
R
& \ll N^\alpha \sum_{d\| q}
{(3g)^{\omega(d)} \over d^\alpha}
=N^\alpha\prod_{p^\nu\|q} \bigg(1+{3g\over p^{\alpha\nu}}\bigg)
\cr
&\ll N^\alpha\prod_{p\leqslant p_{\omega(q)}}\bigg(1+{3g\over p^\alpha}\bigg)
   \ll N^{1-\delta_q} (\log \{3+\omega(q)\})^{3g},
\cr}
\eqdef{f5}$$
pour le choix $\alpha:=1-\delta_q$.
\qed
\goodbreak\goodbreak
\paraunn{Estimations  relatives aux nombres $f(p)$}
Nous utilisons systŽmatiquement la notation \eqref{defu} et posons
$$
\ak=\ak(x,y) := 1 - {\xk\over \log y}\cdot
\eqdef{alphak}
$$
Pour les valeurs relatives de $x$ et $y$
considŽrŽes dans la suite, cette quantitŽ est une
bonne approximation du
point-selle associŽ au calcul de
$\Psi_f(x,y)$ ou
$\psxy$.
\Propl{fp}
{Soient $A$ et $\kappa$ deux constantes positives.
Pour toute \fa\  $f$
satisfaisant ˆ \eqref{eta1cote}, on a
$$f(p)\leqslant (A+\kappa)p/\log p
\eqdef{majfp}$$
pour tout nombre premier $p$.
De plus, il existe deux constantes positives
$c_4$ et $y_0$, ne dŽpendant que de $A$ et $\kappa$, telles que l'on ait
$$f(p)/p^{\ak}\leqslant \dm
\eqdef{majfp/pa}$$
   pour $y_0\le p\le y$ et $1\le u\le c_4(\log y)/\log_2y$.}
\goodbreak
\nid
La relation \eqref{eta1cote} appliquŽe avec $z=p$
et $w=p-1$ fournit immŽdiatement
$$f(p){\log p\over p}
=\kappa \log(1+1/p)+r_f(p)-r_f(p-1)
\leqslant \kappa + A.
\eqdef{majfplog/p}$$

Pour tout $t>0$, on a $t^{1-\ak}/\log
t=\e^{\vartheta(\log t)}$ o
$\vartheta(v):=(1-\ak)v-\log v$ est convexe.
Il s'ensuit que,
pour chaque $y_0$ fixŽ et $y_0\le p\le y$,
$$f(p)p^{-\ak}
\leqslant (A+\kappa) p^{1-\ak}/\log p
\leqslant (A+\kappa)
\big(y_0^{1-\ak}/\log y_0+y^{1-\ak}/\log y\big).$$
Soit $K>0$. Pour $c_4$ assez petite et $y_0$
assez grande, on a $\xk\le \log_2y-K$ dans le
domaine indiquŽ. La
majoration prŽcŽdente est donc $\ll (A+\kappa)\e^{-K}$, d'o l'on dŽduit bien
\eqref{majfp/pa} en choisissant convenablement $K$.
\qed

\Propl{sumfp2log/p2}
{Soient $A$ et $\kappa$ deux constantes positives.
Pour toute \fa\  $f$
satisfaisant ˆ \eqref{eta1cote}, on a
$$\eqalignno{
&\sum_{p} {f(p)^2\over p^2} \leqslant 4(A+\kappa)^2. & \eqdef{moyfp2/p2}
\cr
}
$$
}

\nid La relation \eqref{eta1cote} implique, pour tout entier $j\ge 1$,
$$ \sum_{2^j\le p<2^{j+1}}{f(p)\log p\over p}\leqslant A+\kappa.\eqdef{scfp}$$
Il suit, gr‰ce ˆ \eqref{majfp},
$$\sum_{p} {f(p)^2\over p^2}\leqslant
(A+\kappa)\sum_{p}{f(p)\over p\log p}\leqslant
{(A+\kappa)^2\over (\log 2)^2}\sum_{j\ge
1}{1\over j^2}\leqslant 4(A+\kappa)^2.
$$
\qed

\Propl{sumfp2alpha}
{Soient $A$ et $\kappa$ deux constantes positives.
Pour toute \fa\  $f$
satisfaisant ˆ \eqref{eta1cote} et tous nombres
rŽels $x$, $y$ satisfaisant ˆ $x\ge 2,$ $(\log
x)^3\le y\le x,$ on a
$$\eqalignno{\sum_{p\le y} {f(p)^2\over p^{2\ak}}
&\ll 1+{u^2\xk\over \log y},&\eqdef{fp2logp}\cr
\sum_{p\le y} {f(p)^2\log p\over p^{2\ak}}
&\ll Z_1(y;f)+u^2\xi_\kappa(u)\log
\bigg(1+{Z_2(y;f)\xi_\kappa(u)\over \log
y}\bigg). & \eqdef{majfp2alpha}\cr}
$$
}
\nid
Par  \eqref{majfp} et \eqref{scfp}, nous pouvons Žcrire
$$\eqalign{\sum_{p\le y} {f(p)^2\over
p^{2\ak}}&\leqslant (A+\kappa)\sum_{p\le
y}{f(p)p^{2(1-\ak)}\over
p\log p}
\cr
&\leqslant {4(A+\kappa)^2\over (\log 2)^{2}}\sum_{1\le
j\le (\log y)/\log 2}{2^{2(1-\ak)j}\over j^{2}
}\cdot\cr} $$
On obtient \eqref{fp2logp} en faisant appel ˆ l'estimation
$$\sum_{1\le j\le J}{2^{\beta j}\over j^{2}}\ll
{2^{\beta J}-1\over \beta J^{2}}+1\qquad (J\ge
1,\,\beta\in\,]0,1[),
$$
qu'on Žtablit aisŽment par sommation d'Abel.
\par
Montrons \eqref{majfp2alpha}. Comme la quantitŽ ˆ
estimer est une fonction croissante de $u$, nous
pouvons supposer
$u>u_0$, o $u_0$ est une constante arbitrairement grande. Nous avons
$$\eqalign{E:=\sum_{p\le y} {f(p)^2\log p\over
p^{2\ak}}&\leqslant  Z_1(y;f)+\sum_{p\le y}{f(p)^2(\log
p)^2\over
p^2}{p^{2-2\ak}-1\over \log p}\cdot\cr}
$$
Soient $t$ un paramtre vŽrifiant $1/\log y\ll
t\le \dm$ et $z:=y^{1-t}$. Scindons la dernire
somme au point $z$ et
observons que
$${p^{2-2\alpha_\kappa}-1\over \log p}\leqslant
{z^{2-2\alpha_\kappa}-1\over \log
z}\ll{\{u\xi_\kappa(u)\}^{2-2t}\over \log
y}\qquad (p\le z). $$ Nous
obtenons
$$E\ll Z_1(y; f)
+ \{u\xi_\kappa(u)\}^{2-2t}{Z_2(y; f)\over \log y} +E', $$
avec $$E':=\sum_{z<p\leqslant y}{f(p)^2(\log
p)^2\over p^2}{p^{2-2\alpha_\kappa}-1\over \log
p}.$$
Ensuite, nous avons
$$\eqalign{E'\ll {\{u\xi_\kappa(u)\}^2\over \log
y} \sum_{z<p\leqslant y}{f(p)\log p\over
p}\ll tu^2\xi_\kappa(u)^2,\cr}
$$
o la seconde estimation rŽsulte de l'hypothse \eqref{eta1cote}.
Comme nous avons supposŽ $u$ assez grand, nous
pouvons choisir $2t\log \{u\xi_\kappa(u)\}=\log
\big(1+Z_2(y;f)\xi_\kappa(u)/\log y\big)$, d'o \eqref{majfp2alpha}.
\qed
Pour $f\in\L_\kappa(A,C,\eta)$, nous posons
$$C_\kappa(f)
:= \prod_p (1-1/p)^\kappa \sum_{\nu\ge 0}f(p^\nu)/p^\nu,\eqdef{Ckf}$$
o le produit infini est considŽrŽ comme nul
lorsqu'il est divergent. Il rŽsulte
de~\eqref{eta2cotes} que cette dernire
ŽventualitŽ est exclue lorsque $f\in\MM_\kappa(A,C,\eta).$
\Propl{F1y}
{Soient $A>0$, $C>0$, $\kappa>0$, $\eta\in\,]0, \dm]$.
Il existe une constante $C_0=C_0(A, C, \kappa,
\eta)$ telle que,  pour
$f\in \L_\kappa(A, C, \eta)$,  on ait
$$\leqalignno{
& \e^{\gamma\kappa}C_\kappa(f)(\log y)^\kappa
\bigg(1-{C_0\over \log y}\bigg)
\leqslant F(1, y)
\leqslant C_0(\log y)^\kappa
\qquad(y\ge 2).
& \eqdef{cadreF1yL}
\cr
& {F(1, x)\over F(1, y)}
\geqslant  u^\kappa \bigg(1-{C_0\over \log x}\bigg)
\qquad(y\ge x\geqslant 2).
& \eqdef{rapportF1xyL}
\cr}$$
De plus, on a, uniformŽment pour $f\in \MM_\kappa(A, C, \eta)$,
$$F(1, y)
=\e^{\gamma\kappa}C_\kappa(f)(\log y)^\kappa
\bigg\{1+O\bigg({1\over \log y}\bigg)\bigg\}.
\eqdef{asympF1yM}$$
}

\nid Pour toute fonction $f$ de $\L_\kappa(A,C,\eta)$, nous pouvons Žcrire
$$\eqalign{F(1,y)&\leqslant \exp\Big\{\sum_{p\le
y}{f(p)\over p}+\sum_{p}\sum_{\nu\ge
2}{f(p^\nu)\over p^\nu}\Big\}\cr
&\leqslant
\exp\Big\{C+\kappa\log_2y-\kappa\log_22+G(y)\Big\}
\ll(\log y)^\kappa \e^{G(y)}\cr}$$
o l'on a posŽ, pour tout $z\ge 2$,
$$\eqalign{G(z)&:=\int_{2-}^z{\d\{ r_f(t)-r_f(2-)\}\over \log t}\cr
&={r_f(z)+\kappa\log 2\over \log
z}+\int_2^z{r_f(t)+\kappa\log 2\over t(\log
t)^2}\d t\le {A\over\log 2
}+\kappa,\cr}\eqdef{defW}
$$
en vertu de l'hypothse \eqref{eta1cote}. Cela
Žtablit l'inŽgalitŽ de droite de
\eqref{cadreF1yL}.
\par\goodbreak
Pour montrer l'inŽgalitŽ de gauche de
\eqref{cadreF1yL} et la formule asymptotique
\eqref{asympF1yM}, nous pouvons
clairement supposer la convergence du produit
infini $C_\kappa(f)$. Nous pouvons alors Žcrire
$$F(1, y)
= \prod_{p\le y}\sum_{\nu\ge 0}{f(p^\nu)\over p^\nu}
= C_\kappa(f) \prod_{p\le y}\bigg(1-{1\over p}\bigg)^{-\kappa}
\P_\kappa(f;y)^{-1}
$$
avec
$$\eqalign{\P_\kappa(f;y)
&:=\prod_{p>y}\bigg(1-{1\over p}\bigg)^{\kappa}
\sum_{\nu\ge 0}{f(p^\nu)\over p^\nu}
\cr
& =\exp\bigg\{\sum_{p>y}{f(p)-\kappa\over p}
+ O\bigg({1+f(p)^2\over p^2}
+ \sum_{\nu\ge 2}{f(p^\nu)\over p^\nu}\bigg)\bigg\},\cr}$$
o nous avons fait appel ˆ l'hypoth\`ese
\eqref{moyfpnu} et ˆ l'in\'egalit\'e
\eqref{majfp}.
\par
La formule de Mertens fournit d'abord, classiquement,
$$\prod_{p\le y}\bigg(1-{1\over p}\bigg)^{-\kappa}
= \e^{\gamma\kappa}(\log y)^\kappa
\bigg\{1+O_\kappa\bigg({1\over \log y}\bigg)\bigg\}.$$
\par
    Ensuite, gr‰ce ˆ la seconde formule de Mertens
$$\sum_{p\le t}{1\over p}=\log_2t+a+O\Big({1\over
\log t}\Big)\quad (t\ge 2), $$
o $a$ est une constante convenable, nous pouvons Žcrire
$$\sum_{p>y}{f(p)-\kappa\over
p}=\int_y^\infty{\d\{r_f(t)-r_f(2-)\}\over \log
t}+O\Big({1\over \log y}\Big), $$
et une intŽgration par parties analogue ˆ celle de \eqref{defW} fournit
$$\sum_{p>y}{f(p)-\kappa\over
p}\normalbaselineskip28pt\cases{\dsp\le{A+\kappa\over
\log y}+O\Big({1\over \log y}\Big)&
si
$f\in\L_\kappa(A,C,\eta)$,\cr
\ll\dsp{1\over \log y}& si $f\in\MM_\kappa(A,C,\eta).$\cr} $$
\par
D'autre part, il rŽsulte aisŽment de  \eqref{eta1cote} et \eqref{scfp} que
$$\sum_{p>y}{1+f(p)^2\over p^2}
\ll {1\over \log y}
$$
alors que la condition \eqref{moyfpnu} fournit
$$\sum_{p>y} \sum_{\nu\ge 2}{f(p^\nu)\over p^\nu}
\leqslant \sum_{p} \sum_{\nu\ge 2}{f(p^\nu)\over p^\nu}
\bigg({p\over y}\bigg)^{\nu\eta}
\ll {1\over y^{2\eta}}.$$
Nous dŽduisons donc de ce qui prŽcde que l'on a pour $y$ assez grand
$$\P_\kappa(f;y)\cases{
\leqslant \exp\{C_0/(2\log y)\}
& si $f\in \L_\kappa(A, C, \eta)$,
\cr\noalign{\vskip 2mm}
=\exp\{O({1/\log y})\}
& si $f\in \MM_\kappa(A, C, \eta)$.
\cr}$$
Cela implique bien l'inŽgalitŽ de gauche de \eqref{cadreF1yL}
et la formule \eqref{asympF1yM}.

Finalement, Žtablissons \eqref{rapportF1xyL}.
Avec les notations introduites plus haut, nous avons
$${F(1, y)\over F(1, x)}
= \prod_{x<p\le y}\sum_{\nu\ge 0}{f(p^\nu)\over p^\nu}
= \prod_{x<p\le y}\bigg(1-{1\over p}\bigg)^{-\kappa}
\P_\kappa(f;x,y)
$$
avec
$$\P_\kappa(f;x,y)
:=\prod_{x<p\le y}\bigg(1-{1\over p}\bigg)^{\kappa}
\sum_{\nu\ge 0}{f(p^\nu)\over p^\nu}
\leqslant \exp\bigg\{{C_0\over \log x}\bigg\}.$$
L'in\'egalit\'e annonc\'ee
dŽcoule donc de la formule de Mertens.
\qed

\medskip
\paraunn{Majoration prŽliminaire}

Les d\'emonstrations des ThŽormes 
\ref{sthmmajoration} et \ref{sthmasymp} reposent
de manire essentielle sur la
disponibilitŽ d'une \og bonne\fg\ majoration pour
$\psi_f^*(x,y)$. Dans le cas des fonctions de
$\L_\kappa(A,C,\eta)$,
une estimation relativement grossire est
suffisante. C'est ce que nous obtenons au
\ref{psi*premier}
ci-dessous. Dans le cas des fonctions de
$\MM_\kappa(A,C,\eta)$, en revanche, nous visons
une formule asymptotique et la
majoration optimale ˆ
   constante multiplicative prs
$$\psi_f^*(x, y)\ll F(1,y)\lk(u)\eqdef{majoptps}$$
s'avre indispensable pour prŽserver la qualitŽ
du rŽsultat final. Cette estimation nous sera en
fait fournie,
dans le domaine requis en
$(x,y)$, par le
\ref{thmmajoration}.
\par
Nous obtenons notre majoration prŽliminaire par une simple
application de la m\'ethode de Rankin. Gr‰ce au
choix prŽcis des paramtres issu de la mŽthode du
col, cela
   fournit incidemment une amŽlioration significative du lemme 1 de \citer{G05}.

\Propl{psi*premier}
{Soient $A>0$, $C>0$, $\kappa>0$ et $\eta\in\,]0,\dm[$.
On a
$$\psi_f^*(x, y)
\ll F(1, y) \lk(u)\,\xk\sqrt{u+1}\,\e^{Au\xk/\log y}
\eqdef{1eremajpsi*}$$
uniformŽment pour
$$f\in \L_\kappa(A, C, \eta),
\quad
y\ge 2,\quad
0\le u\le y^\eta/\log y.
\eqdef{3.2}$$
En particulier, on a
$$\psi_f^*(x, y)
\ll F(1, y) \lk(u) \, \xi_\kappa(u)\sqrt{u+1}
\eqdef{3.3}$$
uniformŽment pour
$f\in \L_\kappa(A, C, \eta),
\,
y\ge 3,\,
0\le u\le (\log y)/\log_2y.
$
\PAR
Sous l'hypoth\`ese suppl\'ementaire \eqref{fpnu0},
la majoration \eqref{1eremajpsi*} a lieu
pour tous $x>1$, $y\ge 2$.}
\goodbreak
\nid Nous pouvons supposer, sans perte de gŽnŽralitŽ, que $y$ est assez grand.
Nous appliquons la m\'ethode de Rankin. Soit $\sigma\in\,]0,1]$. Sous
l'hypothse
$$ K_\sigma:=\sum_{p\le y}\sum_{\nu\ge
2}f(p^\nu)/p^{\nu\sigma}<\infty,$$ nous pouvons Žcrire
$$\eqalign{\psi_f^*(x, y)
& \leqslant \sum_{\scriptstyle n>x\atop\scriptstyle P(n)\leqslant y}
{f(n)\over n} \bigg({n\over x}\bigg)^{1-\sigma}
\leqslant x^{\sigma-1}
\prod_{p\le y}\sum_{\nu\ge 0}{f(p^\nu)\over p^{\nu\sigma}}
\cr
& \leqslant x^{\sigma-1}
\exp\bigg\{
\sum_{p\le y}\sum_{\nu\ge 1}{f(p^\nu)\over p^{\nu\sigma}}
\bigg\}
   \leqslant x^{\sigma-1}
\exp\bigg\{K_\sigma+\sum_{p\le y} {f(p)\over p^{\sigma}}\bigg\}.
\cr}$$
D'aprs
l'hypoth\`ese \eqref{moyfpnu}, nous pouvons
choisir $\sigma=\ak\ge 1-\eta$ si \eqref{3.2} est
rŽalisŽe. Lorsque
\eqref{fpnu0} est satisfaite, le mme choix est
possible sans restriction puisque $\ak\in\,]0,1[$
si
$x\le N_y:=\prod_{p\le y}p$ et  $\psxy=0$ dans le cas contraire.
\par
Comme on a, gr\^ace \`a \eqref{moyfp2/p2},
$$F(1, y)
\geqslant \prod_{p\le y}\bigg(1+ {f(p)\over p}\bigg)
\gg
\exp\bigg\{\sum_{p\le y} {f(p)\over p}\bigg\},$$
il suit
$$\psi_f^*(x, y)
\ll F(1, y)
\exp\bigg\{-u\xk
+ \int_0^{1-\ak}\sum_{p\le y} {f(p)\log p\over p^{1-v}}\d v\bigg\}.
\eqdef{psi*xy/F1y}$$
La dernire somme en $p$ n'excde pas
$$\int_{1}^y{\kappa\over t^{1-v}}\d
t+\int_{2-}^yt^v\d\{r_f(t)-r_f(y)\}\leqslant\int_0^{v\log
y}\e^{w}\d w+Ay^v.  $$
D'o
$$\psxy\le F(1,y)\exp\bigg\{-u\xk
+ \int_0^{\xk}{\e^w-1\over w}\d w+{Au\xk\over \log y}\bigg\}.$$
Compte tenu de \eqref{farhok} et \eqref{evallk},
cela fournit bien le rŽsultat annoncŽ.
\qed

\paraunn{ƒquations fonctionnelles} \drefun{eqfonc}
\paradeuxn{Forme initiale}
Dans le lemme suivant, nous explicitons la forme
initiale de  l'Žquation fonctionnelle approchŽe
de type
Wirsing--Hildebrand pour la fonction
$\psi_f^*(x,y)$. ƒtant donnŽs une fonction
multiplicative $f$ et un entier $m$, nous
dŽfinissons une nouvelle fonction
multiplicative par la formule
$$f_m(n):=\cases{f(n) & si $(n,m)=1$,\cr
0 & si $(n,m)>1$.\cr} \eqdef{fm}$$
Nous convenons Žgalement de poser
$$r_f(z)=0 \quad (0\leqslant z\leqslant 1).\eqdef{conventionrf}$$
\Propl{Equpsi*}
{Soient $A>0$, $C>0$, $\kappa>0$,
$\eta\in\,]0,\dm[$ et $f\in\L_\kappa(A,C,\eta)$.
Pour tous $x\ge 2$, $y\ge 2$, on a
$$\psi_f^*(x, y) \log x + \int_x ^\infty {\psi_f^*(t, y)\over t} \d t
= \kappa \int_{x/y}^x {\psi_f^*(t, y)\over t} \d t
+ \sum_{1\le j\le 3}E_j(x,y;f),
\eqdef{equfonc}$$
o\`u l'on a posŽ
$$\left\{\eqalign{
E_1(x,y;f)&:=-\sum_{\di{n>x/y}{P(n)\leqslant y}}{f(n)\over
n}\Big\{r_f(y)-r_f\Big({x\over n}\Big)\Big\},\cr
E_2(x,y;f)&:=-\sum_{\nu\ge 1}\sum_{p\le
y}{f(p)f(p^\nu)\log p\over
p^{\nu+1}}\psi_{f_p}^*(x/p^{\nu+1},y),\cr
E_3(x,y;f)&:=\sum_{\nu\ge 2}\sum_{p\le y}{f(p^\nu)\log p^\nu\over
p^{\nu}}\psi_{f_p}^*(x/p^{\nu},y).\cr}\right.
\eqdef{RfM+}$$
}

\nid Nous avons d'une part
$$\eqalign{\sum_{\di{ n>x}{P(n)\leqslant y}}
{f(n)\log n\over n}&=\sum_{\di{ n>x}{P(n)\leqslant y}}
{f(n)\over n}\int_1^n{\d t\over t}\cr
&=\psxy\log x+\int_x^\infty{\psi_f^*(t,y)\over t}\d t.\cr}\eqdef{Hild1}$$
D'autre part, l'identitŽ $\log n=\sum_{p^\nu\| n}\log p^\nu$,
fournit immŽdiatement
$$\eqalign{\sum_{\di{ n>x}{P(n)\leqslant y}}
{f(n)\log n\over n}
& = \sum_{\di{mp^\nu>x
}{P(mp)\leqslant y,\,p\,\nmid\,m}}
{f(m)f(p^\nu)\log p^\nu\over mp^\nu}
\cr
& = \sum_{\nu\ge 1}\sum_{p\le y}{f(p^\nu)\log
p^\nu\over p^\nu}\psi_{f_p}^*(x/p^\nu,y).
\cr}
\eqdef{equation}$$
Nous considŽrons comme terme principal du membre de droite la quantitŽ
$$\eqalign{S&:=\sum_{p\le y}{f(p)\log p\over
p}\psi_f^*(x/p,y)=\sum_{\di{n>x/y}{P(n)\leqslant
y}}{f(n)\over n}\sum_{x/n<p\le
y}{f(p)\log p\over p}\cr
&=\sum_{\di{x/y<n\le x}{P(n)\leqslant 
y}}{f(n)\over n}\Big\{\kappa\log \Big({ny\over
x}\Big)+r_f(y)-r_f\Big({x\over
n}\Big)\Big\}+\psxy\big\{\kappa\log
y+r_f(y)\big\}.\cr}
$$
En observant que
$$\eqalign{\sum_{\di{x/y<n\le x}{P(n)\leqslant 
y}}{f(n)\over n}\log \Big({ny\over
x}\Big)&=\int_{x/y}^x{\psi_f^*(t,y)-\psi_f^*(x,y)\over t}\d t\cr
&=\int_{x/y}^x{\psi_f^*(t,y)\over t}\d t-\psi_f^*(x,y)\log
y,\cr}
$$
nous obtenons
$$S:=\kappa \int_{x/y}^x{\psi_f^*(t,y)\over t}\d t+E_1(x,y;f).$$
Par ailleurs, on vŽrifie aisŽment que
$$
\sum_{p\le y}{f(p)\log p\over p}\psi_{f_p}^*(x/p,y)-S=E_2(x,y;f).$$
Cela implique bien le rŽsultat annoncŽ.
\qed
\medskip
\paradeuxn{Exploitation de l'Žquation adjointe}

Il rŽsulte aisŽment de l'Žquation \eqref{rho} et
de la dŽfinition \eqref{lambda} que la fonction
$\lambda_\kappa$
est solution de l'\'equation diff\'erentielle aux diff\'erences
$$u\lambda'_\kappa(u)
= \kappa \big\{\lk(u) - \lambda_\kappa(u-1)\big\}.
\eqdef{eqfonclambda}$$
La relation \eqref{equfonc} est une version
discrte de cette Žquation Žcrite sous forme
intŽgrale.
C'est  la comparaison des deux formules  qui
conduira aux estimations des ThŽormes
\ref{sthmmajoration} et
\ref{sthmasymp}.

Comme dans \citer{Te01}, nous utilisons de
manire essentielle l'\'equation adjointe
de~\eqref{eqfonclambda}, soit
$$\big\{ug(u)\big\}'
= \kappa \big\{g(u+1) - g(u)\big\}
\eqdef{eqfoncmu}$$
dont une solution est la fonction dŽcroissante
$$\mu_\kappa(u)
:= \int_0^\infty
\exp\bigg(-uv + \kappa \int_0^v {1-\e^{-w}\over w}\d w\bigg) \d v
\qquad(u>0).
\eqdef{rk}
$$
Nous observons immŽdiatement qu'une manipulation
de routine fournit l'Žvaluation asymptotique
$$\mu_\kappa(u)\sim 1/u
\qquad(u\to\infty).
\eqdef{rkasymp}$$
\par
L'identitŽ
$$u\lambda_\kappa(u)\mu_\kappa(u)
=\kappa\int_{u-1}^u\lambda_\kappa(v)\mu_\kappa(v+1)\d
v\qquad (u\ge 1),\eqdef{pdtscal}$$
qui dŽcoule aisŽment de \eqref{eqfonclambda} et
\eqref{eqfoncmu} est un outil technique
particulirement bien adaptŽ
ˆ l'exploitation des diverses Žquations
fonctionnelles du problme considŽrŽ dans ce
travail.
\par\goodbreak

Pour $x>0$, $y\ge 2$,  et
$f\in \L_\kappa(A, C, \eta)$,
nous posons
$$\nuf(u):=\psi_f^*(x, y)/F(1,y),\eqdef{defnufyu}$$
o $F(s,y)$ est dŽfinie en \eqref{Fsy}.
\par
Nous posons, avec les notations \eqref{RfM+} et \eqref{defu},
$$R_f^*(u;y):=\sum_{1\le j\le 3}E_j(x,y;f).\eqdef{defRf*} $$

\Propl{Efonc}
{Soient $A>0$, $C>0$, $\kappa>0$, $\eta\in\,]0, \dm[$
et $f\in \L_\kappa(A, C, \eta)$.
Pour $y>1$, $u\ge 1$, nous avons
$$u \, \mu_\kappa(u)\nuf(u)
= \kappa \int_{u-1}^u \mu_\kappa(v+1) \nuf(v) \d v
-\int_u^\infty \mu_\kappa(v) \d\bigg\{{R_f^*(v;y)\over F(1, y)\log y}\bigg\}.
\eqdef{eqfoncnuf}$$
}

\nid
Reportons \eqref{defnufyu} dans \eqref{equfonc},
effectuons le changement de variables $v=(\log t)/\log y$
dans les deux int\'egrales et divisons par $F(1, y)\log y$.
Nous obtenons
$$\eqalign{
& u \,\nuf(u)
+ \int_u^\infty\nuf(v) \d v
= \kappa \int_{u-1}^u\nuf(v) \d v
+ {R_f^*(u;y)\over F(1, y)\log y}.
\cr}$$
Aprs diffŽrentiation relativement ˆ $u$,
multiplication par $\mu_\kappa(u)$
et intŽgration de $u$ ˆ l'infini, cette Žquation prend la forme
$$\eqalign{
\int_u^\infty &\mu_\kappa(v) \d\bigg\{{R_f^*(v;y)\over F(1, y)\log y}\bigg\}
\cr
& = \int_u^\infty v\mu_\kappa(v) \d\nuf(v)
- \kappa \int_u^\infty\mu_\kappa(v)\{\nuf(v)-\nuf(v-1)\} \d v.
\cr}\eqdef{eqfoncinterm}$$
Or, une int\'egration par parties  tenant compte des relations
$\lim_{u\to\infty}u\mu_\kappa(u)=1$,
$\lim_{u\to\infty}\nuf(u)=0$, fournit, gr\^ace \`a \eqref{eqfoncmu},
$$\int_u^\infty v\mu_\kappa(v) \d\nuf(v)
= - u\mu_\kappa(u)\nuf(u)
- \kappa \int_u^\infty \nuf(v)
\{\mu_\kappa(v+1)-\mu_\kappa(v)\} \d v.$$
Nous obtenons donc le r\'esultat annoncŽ en
reportant dans \eqref{eqfoncinterm}.
\qed

\paraunn{Estimations initiales}
Le processus itŽratif utilisŽ pour dŽmontrer les
ThŽormes \ref{sthmmajoration} et \ref{sthmasymp}
nŽcessite une
Žvaluation prŽcise des conditions initiales, qui
correspondent ici au cas $y=x$. Ces
renseignements sont obtenus au
\ref{psi*u=1} {\it infra}. Pour Žtablir ce rŽsultat, nous commen\c cons par
\'enoncer une \'equation fonctionnelle pour
$\psi_f(x):=\psi_f(x, x)$ analogue ˆ celle du
\ref{Equpsi*} pour $\psxy$.
Nous conservons la notation \eqref{fm} et posons pour $x\ge 1$
$$\left\{\eqalign{D_1(x;f)&:=\sum_{n\le
x}{f(n)\over n}r_f\Big({x\over n}\Big),\cr
D_2(x;f)&:=-\sum_{\nu\ge 1}\sum_{p^{\nu+1}\leqslant
x}{f(p)f(p^\nu)\log p\over
p^{\nu+1}}\psi_{f_p}(x/p^{\nu+1}),\cr
D_3(x;f)&:=\sum_{\nu\ge 2}\sum_{p^\nu\le x}{f(p^\nu)\log p^\nu\over
p^{\nu}}\psi_{f_p}(x/p^{\nu}),\cr
R_f(x)&:=\sum_{1\le j\le 3}D_j(x;f).\cr}\right.
\eqdef{termesD}$$
Nous omettons la dŽmonstration, qui est trs voisine de celle du \ref{Equpsi*}.

\Propl{Equapsixy}
{Pour toute fonction multiplicative $f$ et tout $x\ge 1$,
on a
$$\psi_f(x)\log x
= (\kappa+1)\int_1^x {\psi_f(t)\over t} \d t
+ R_f(x).
\eqdef{eqpsix}$$
}

Les dŽmonstrations des Th\'eor\`emes
\ref{sthmmajoration} et \ref{sthmasymp}  reposent
sur l'exploitation de
l'Žquation fonctionnelle \eqref{equfonc}. Le
lemme suivant permet fournit les valeurs
initiales. Rappelons que la
fonction
$j_\kappa$ dŽfinie en \eqref{lambda} vŽrifie
$$j_\kappa(1)=1-\lk(1)=\e^{-\gamma\kappa}/\Gamma(\kappa+1).$$

\Propl{psi*u=1}
{Soient $A>0$, $C>0$, $\kappa>0$ et $\eta\in\,]0,\dm[$.
\par
(i)
On a
$$\psi_f(x)
= F(1, x)j_\kappa(1)
\bigg\{1+O\bigg({Z_1(x;f)\over \log x}\bigg)\bigg\}
\eqdef{psifx}$$
uniform\'ement pour
$f\in \MM_\kappa(A, C, \eta)$ et $x\ge 2$.
\PAR
(ii)
Il existe une constante positive $C_0=C_0(A, C, \kappa, \eta)$
telle que l'on ait
$$\psi_f(x)
\geqslant F(1, x)j_\kappa(1)\bigg\{1-{C_0\over 1+\log x}\bigg\}\eqdef{psifxL}$$
pour
$f\in \L_\kappa(A, C, \eta)$ et $x\ge 1$.
}
\goodbreak
\nid
Il rŽsulte immŽdiatement de \eqref{termesD},
\eqref{moyfpnu}, et
\eqref{eta1cote} ou
\eqref{eta2cotes} que l'on~a
$$\eqalign{R_f(x)&\ll \psi_f(x)Z_1(x;f)\quad
\hbox{si\ } f\in\MM_\kappa(A,C,\eta),\cr
\noalign{\vskip-3mm}\cr
R_f(x)&\leqslant B\psi_f(x)\quad \hbox{si\ } f\in\L_\kappa(A,C,\eta),\cr} $$
o $B$ est une constante positive convenable. En
particulier, nous avons, pour $x_0$ convenable,
$$U_f(x):=R_f(x)/\psi_f(x)\leqslant \dm \log x\qquad (x\ge x_0)$$
et
nous pouvons dŽduire de \eqref{eqpsix} que la fonction
$$\vartheta_f(x):={\kappa+1\over (\log
x)^{\kappa+1}}\int_1^x {\psi_f(t)\over t} \d t$$
vŽrifie pour $x$ assez grand
$${\vartheta_f'(x)\over
\vartheta_f(x)}={(\kappa+1)U_f(x)\over x(\log
x)\{\log x-U_f(x)\}}\cdot \eqdef{thpr/th}$$
\par
Par intŽgration sur $[x,\infty[$, il suit,
lorsque $f\in\MM_\kappa(A,C,\eta)$, pour une
constante convenable~$K$,
$$\log \{\vartheta_f(x)/K\}\ll{Z_1(x;f)\over \log x} \qquad (x\ge x_0),$$
d'o, en reportant dans \eqref{eqpsix},
$$\psi_f(x)=K(\log x)^\kappa\Big\{1+O\Big({Z_1(x;f)\over \log x}\Big)\Big\}. $$
La valeur de $K$ peut tre obtenue de diverses
manires. Nous employons, par exemple, le
thŽorme taubŽrien de
Hardy--Littlewood--Karamata sous la forme donnŽe
au thŽorme II.7.8 de \citer{Te95}, qui fournit,
en posant
$\F(s):=\sum_{n\ge 1}f(n)/n^s$,
$$K:={1\over
\Gamma(\kappa+1)}\lim_{\sigma\to1+}{\F(\sigma)(\sigma-1)^{\kappa}}
={C_\kappa(f)\over\Gamma(\kappa+1)}\cdot
$$ L'assertion (i) de notre ŽnoncŽ rŽsulte alors de \eqref{asympF1yM}.
\par\goodbreak
Lorsque $f\in\L_\kappa(A,C,\eta)$, nous
commenons par fixer $x\ge 2$ et remplacer
$f(p^\nu)$ par
$\tau_\kappa(p^\nu):={\kappa+\nu-1\choose
\nu}$ pour $p>x$, ce qui ne modifie pas la valeur
de $\psi_f(x)$. On vŽrifie alors aisŽment, par
exemple gr‰ce au
thŽorme taubŽrien de Hardy--Littlewood--Karamata, que
$$\eqalign{\lim_{t\to\infty}{\psi_f(t)\over (\log
t)^\kappa}&=\lim_{t\to\infty}\vartheta_f(t)
\cr&={C_\kappa(f)\over
\Gamma(\kappa+1)}={1\over
\Gamma(\kappa+1)}\prod_{p\le
x}\big(1-1/p\big)^\kappa\sum_{\nu\ge
0}f(p^\nu)/p^\nu\cr&={j_\kappa(1)F(1,x)\over (\log
x)^\kappa}\Big\{1+O\Big({1\over \log x}\Big)\Big\}.\cr}\eqdef{Karam}$$
\par
Ensuite, nous rŽŽcrivons \eqref{eqpsix} sous la forme
$${\psi_f(t)\over (\log
t)^\kappa}=\vartheta_f(t)+{R_f(t)\over (\log
t)^{\kappa+1}} \eqdef{psithR}$$
d'o nous dŽduisons
$$ \d\Big\{{\psi_f(t)\over (\log
t)^\kappa}\Big\}={\d R_f(t)\over (\log
t)^{\kappa+1}}={\d \{R_f(t)-R_f(x)\}\over (\log
t)^{\kappa+1}},$$
et donc, pour $z\ge x$,
$${\psi_f(z)\over (\log
z)^\kappa}-{\psi_f(x)\over (\log
x)^\kappa}={R_f(z)-R_f(x)\over (\log
z)^{\kappa+1}}+(\kappa+1)\int_x^z{R_f(t)-R_f(x)\over t(\log
t)^{\kappa+2}}\d t. \eqdef{complocps/log}$$
Par \eqref{termesD},
\eqref{moyfpnu}, et
\eqref{eta1cote}, nous avons, pour une constante
convenable $D$ et tout $t\ge x$,
$$R_f(t)-R_f(x)\leqslant D\psi_f(t)\leqslant
DF(1,t)\leqslant DF(1,x)\Big({\log t\over
\log x}\Big)^\kappa\Big\{1+O\Big({1\over \log x}\Big)\Big\}.$$
En notant que \eqref{Karam} et \eqref{psithR} impliquent
$$\lim_{z\to\infty}{R_f(z)-R_f(x)\over (\log z)^{\kappa+1}}=0, $$
nous obtenons, en faisant tendre $z$ vers l'infini dans \eqref{complocps/log},
$${C_\kappa(f)\over
\Gamma(\kappa+1)}-{\psi_f(x)\over (\log
x)^\kappa}\leqslant {\e^{\gamma \kappa}DC_\kappa(f)\over
\log x}\Big\{1+O\Big({1\over \log x}\Big)\Big\}.$$
Compte tenu de la valeur de $C_\kappa(f)$ donnŽe
en \eqref{Karam}, cela Žtablit bien la validitŽ de
l'as\-ser\-tion~(ii).
   \qed

\medskip

\paraunn{DŽmonstration du \ref{thmmajoration}}
Posons
$\delta_{f,y}(u):=\nuf(u)/\lambda_\kappa(u)=\psi_f^*(x,y)/\{F(1,y)\lambda_\kappa(u)\}$ 

$(u\geqslant 0)$ et
$$\deltaf^*(u):=1+\sup_{1/2\leqslant v\leqslant
u}\{\deltaf(v)-1\}^+\qquad (u\geqslant \dm). $$
Nous devons Žtablir l'existence d'une constante
$B=B(A,C,\eta,\kappa)$ telle que, pour toute
fonction $f$ de
$\L_\kappa(A,C,\eta)$, on ait
$$\deltaf^*(u)\leqslant \e^{Bu\xk/\log y}\qquad
(1\leqslant u\leqslant y^{\eta/3}),
\eqdef{majdeltaf*L}$$
l'inŽgalitŽ Žtant valable sans restriction sous
l'hypothse supplŽmentaire \eqref{fpnu0}.

\Propl{rkL}
{Soient $A>0$, $C>0$, $\kappa>0$ et $\eta\in\,]0, \dm[$.
Il existe une constante positive $M_0$
telle que l'on ait
$$-\int_u^\infty \mu_\kappa(v)
\d\bigg\{{R_f^*(v;y)\over F(1, y)}\bigg\}
\leqslant M_0\lk(u)\Big\{\xk\deltaf(u-1)+ \e^{Au\xk/\log y}\Big\}
\eqdef{majrkL}$$
   pour
$$f\in \L_\kappa(A, C, \eta),
\quad
y\ge 3
\quad\hbox{et}\quad
1\le u\le y^{\eta/3}.
\eqdef{domainerkL}$$
De plus, sous l'hypoth\`ese suppl\'ementaire \eqref{fpnu0},
l'in\'egalit\'e \eqref{majrkL} est valable pour
$f\in \L_\kappa(A, C, \eta)$, $x>1$,
$y\ge 2$, en omettant le terme exponentiel du membre de droite.}
\nid
Posons, avec les notations
\eqref{RfM+},
$$I_j:={-1\over
F(1,y)}\int_u^\infty\mu_\kappa(v)\d
E_j(y^v,y;f)\qquad (1\le j\le 3),\eqdef{defIj}$$
et notons d'emblŽe que, puisque la fonction $E_2(x,y;f)$ dŽfinie en
\eqref{RfM+} est dŽcroissante en $x$, nous avons
$I_2\leqslant  0$. Nous pouvons donc dŽduire de
\eqref{eqfoncnuf} que
$$u\,\mu_\kappa(u)\nuf(u)\leqslant \kappa
\int_{u-1}^u \mu_\kappa(v+1) \nuf(v) \d
v+{I_{1}+I_3\over
\log y}.\eqdef{inegbasenu}$$
\par
Nous avons d'abord, par sommation d'Abel,
$$I_{1}=-\int_u^\infty
\{E_{1}(x,y;f)-E_1(y^v,y;f)\}\mu'_\kappa(v)\d v.
\eqdef{represintI1}$$
Le terme entre accolades vaut
$$\eqalign{
& \sum_{\di{n>x/y}{P(n)\leqslant y}}
{f(n)\over n}\big\{r_f(y)-r_f(x/n)\big\}-
\sum_{\di{n>y^{v-1}}{P(n)\leqslant y}}{f(n)\over
n}\big\{r_f(y)-r_f(y^v/n)\big\}
\cr
& \quad =\sum_{\di{x/y<n<y^{v-1}}{P(n)\leqslant y}}
{f(n)\over n}\big\{r_f(y)-r_f(x/n)\big\}
+\sum_{\di{n>y^{v-1}}{P(n)\leqslant y}}{f(n)\over n}
\big\{r_f(y^v/n)-r_f(x/n)\big\}
\cr
&\quad \leqslant A\psi_f^*(x/y,y).
\cr}$$
Il s'ensuit que
$$\eqalign{I_1
&\leqslant {A\mu_\kappa(u)\psi_f^*(x/y,y)\over F(1,y)}
= {A\mu_\kappa(u)\lambda_\kappa(u-1)\deltaf(u-1)}
\cr
&\leqslant A_1\lambda_\kappa(u)\xk\deltaf(u-1)
\cr}
\qquad (x>1,\,y>1),
\eqdef{majI1L}
$$
pour une constante convenable $A_1$, en vertu de
\eqref{comploclam} et \eqref{rkasymp}.
\par\goodbreak
Une int\'egration par parties et
le fait que
$\mu_\kappa'(v)\leqslant 0$ permettent d'Žcrire
$$\eqalign{
I_3
& \leqslant {\mu_\kappa(u)\over F(1,y)}
\sum_{\nu\ge 2}
\sum_{p\le y}{f(p^\nu)\over p^{\nu}}
\psi_{f}^*\bigg({x\over p^\nu},y\bigg).
\cr}\eqdef{majI3Lbase}$$
Majorons $\psi_{f}^*(x/p^\nu, y)$ en appliquant
\eqref{1eremajpsi*}. Posant $u_p:=(\log p)/\log
y$, nous obtenons
$$\eqalign{\psi_{f}^*(x/p^\nu, y)
& \ll F(1,y)\lk(u-\nu u_p) \xk u^{1/2}\e^{Au\xk/\log y}
\cr
& \ll F(1,y)\lambda_\kappa(u)p^{\eta\nu/2}
\xk u^{1/2}\e^{Au\xk/\log y}
\cr}$$
lorsque $p^\nu\le x/\sqrt{y}$. Dans le cas
contraire, nous utilisons l'inŽgalitŽ triviale
$\psi_{f}^*(x/p^\nu, y)\leqslant  F(1,y)$.
Il vient, pour des constantes positives
$M_j=M_j(A, C, \kappa, \eta)$ convenables,
$$\eqalign{I_3
& \leqslant {M_1\lk(u)\xk\e^{Au\xk/\log y}\over u^{1/2}}
\sum_{\nu\ge 2}\sum_{p\le y}{f(p^\nu)\over p^{\nu(1-\eta)}}
+\sum_{\nu\ge 2}\!\sum_{\di{p\le y}{p^\nu>x/\sqrt{y}}}
{f(p^\nu)\over p^{\nu(1-\eta)}x^{\eta/2}}
\cr
& \leqslant M_2\lk(u)\e^{Au\xk/\log y}+x^{-\eta/2}
\leqslant  M_3\lk(u)\e^{Au\xk/\log y},
\cr}
$$
puisque
$x^{-\eta/2}
=\exp\{-(\eta/2)u\log y\}
\leqslant \exp\{-(3/2)u\log u\}
\ll \lk(u)$ dans le domaine $u\le y^{\eta/3}$.
Cela achve la dŽmonstration de l'inŽgalitŽ \eqref{majrkL}.
\par
Sous l'hypoth\`ese \eqref{fpnu0}, l'intŽgrale
$I_3$ est nulle. Le rŽsultat annoncŽ dŽcoule
alors, par les mmes calculs,
de la validitŽ inconditionnelle de~\eqref{1eremajpsi*}.
\qed
\goodbreak
Nous sommes ˆ prŽsent en mesure de complŽter la
dŽmonstration du \ref{thmmajoration}.
\par
Nous commenons par observer que,
d'aprs
\eqref{psifxL} et
\eqref{rapportF1xyL}, nous avons, pour une
constante positive convenable $D$ et tous $x,\,y$
tels que $2\le x\le y$,
$$\eqalign{\psi_f^*(x, y)
& = F(1, y)-\psi_f(x)
\leqslant F(1, y)- F(1, x) j_\kappa(1)
\bigg(1-{D\over 1+\log x}\bigg)
\cr
& = F(1, y)\bigg\{1- j_\kappa(1){F(1, x)\over F(1, y)}
\bigg(1-{D\over1+ \log x}\bigg)\bigg\}
\cr
& \leqslant F(1, y)\bigg\{1-j_\kappa(1)u^\kappa
\bigg(1-{2D\over 1+\log x}\bigg)\bigg\}.
\cr}$$
Comme $\lambda_\kappa(u)=1-j_\kappa(1)u^\kappa$
pour $0\leqslant u\leqslant 1$, nous avons donc
obtenu, l'existence
d'une constante $M=M(A,C,\eta,\kappa)>0$
telle que l'on ait
$$\delta_{f, y}(u)
\le1+{Mu^{\kappa}\over1+u \log y}\qquad (0\leqslant u\leqslant 1)
\eqdef{psi*u<1L}$$
pour $f\in \L_\kappa(A, C, \eta)$,
$y\ge 2$  et en particulier
$$\deltaf^*(1)=1+\sup_{1/2\leqslant u
\leqslant 1}\{\deltaf(u)-1\}^+
\leqslant 1+{2M\over \log y}\cdot
\eqdef{majdeltaf1/2-1}$$

Lorsque $1\leqslant w\leqslant u\leqslant \ft32$,
nous avons, en vertu de \eqref{inegbasenu}, de la
premire inŽgalitŽ
\eqref{majI1L} et de \eqref{majI3Lbase},
$$w\,\mu_\kappa(w)\lambda_\kappa(w)\deltaf(w)\leqslant
\kappa \int_{w-1}^w
\mu_\kappa(v+1)\lambda_\kappa(v) \deltaf(v) \d
v+{M_4\over \log y}\cdot
$$ Ici et dans toute la suite de cette
dŽmonstration nous dŽsignons par $M_j$
$(j\geqslant 4)$ des constantes positives
convenables. Scindons l'intŽgrale ˆ $w-\dm$ et
majorons la contribution de l'intervalle
$w-1\leqslant v\leqslant w-\dm$
   en estimant
$\deltaf(v)$ par~\eqref{psi*u<1L}. Il suit, pour une constante convenable
$M_5$,
$$\eqalign{\deltaf(w)
&\leqslant 1+a(w)\{\deltaf^*(u)-1\}+{M_5\over \log y}
\cr
&\leqslant \dm+\dm \deltaf^*(u)+{M_5\over \log y},\cr}
$$
o nous avons posŽ
$$a(w):= {\kappa\over w\mu_\kappa(w)\lambda_\kappa(w)}
\int_{w-1/2}^w\mu_\kappa(v+1)\lambda_\kappa(v)\d v
\leqslant \dm,$$
  la dernire inŽgalitŽ rŽsultant de \eqref{pdtscal} et
de la dŽcroissance de la fonction $v\mapsto \mu_\kappa(v+1) \lambda_\kappa(v)$.

En prenant le supremum pour $1\leqslant
w\leqslant u$, et en tenant compte de
\eqref{majdeltaf1/2-1}, nous obtenons donc
que
$$\deltaf(\ft32)=1+\sup_{1/2\leqslant u\leqslant
3/2}\{\deltaf(u)-1\}^+\leqslant 1+{M_6\over \log
y}\cdot\eqdef{majdeltaf1/2-3/2}  $$
\par

Pour tout $u$ de $[\ft32,y^{\eta/3}]$ et tout $w$ de $[u-1,u]$,
nous pouvons Žcrire,
gr‰ce ˆ \eqref{eqfoncnuf} et \eqref{majrkL},
$$\eqalign{
w\mu_\kappa(w)\lambda_\kappa(w)\delta_{f, y}(w)
& \leqslant \kappa \int_{w-1}^{w} \mu_\kappa(v+1) 
\lambda_\kappa(v)\deltaf(v) \d v
\cr
& \quad
+ {M_0\lk(w)\over \log
y}\Big\{\xi_\kappa(w)\deltaf(w-1)+
\e^{Aw\xi_\kappa(w)/\log y}\Big\}.
\cr}$$
Majorons $\deltaf(v)$ par $\deltaf^*(u-\dm)$ si
$v\leqslant u-\dm$ et par $\deltaf^*(u)$ si
$u-\dm <v\leqslant u$ puis $\delta_{f,y}(w-1)$ par $\deltaf^*(u-\dm)$.
Divisons ensuite l'inŽgalitŽ obtenue par $w\mu_\kappa(w)\lambda_\kappa(w)$.
Il suit
$$\deltaf(w)
\leqslant \bigg(1-a(w)+{M_7\xk\over \log y}\bigg)
\deltaf^*(u-\dm)
+a(w)\deltaf^*(u)
+{M_7\over \log y}\e^{Aw\xi_\kappa(w)/\log y}.$$

\par
Ë ce stade, observons que si
$\deltaf^*(u)>\deltaf^*(u-\dm)$, alors
$$\deltaf^*(u)=\sup_{u-1/2\leqslant w\leqslant
u}\deltaf^*(w).$$   D'o\`u, dans cette circonstance,
$$\deltaf^*(u)
\leqslant \bigg(1+{2M_7\xk\over \log y}\bigg)
\delta_{f, y}^*(u-\dm)
+ {2M_7\e^{Au\xk/\log y}\over \log y}\cdot
\eqdef{24nu3}$$
Par une rŽcurrence immŽdiate, cela implique, pour
un choix convenable de $B$, la validitŽ de
\eqref{majdeltaf*L} dans le domaine indiquŽ.
\par
Sous l'hypothse supplŽmentaire \eqref{fpnu0}, le second terme du membre de
droite de~\eqref{24nu3} peut tre omis et l'inŽgalitŽ est
valable pour $x\geqslant y\geqslant 2$.
L'itŽration de~\eqref{24nu3}
conduit donc au rŽsultat indiquŽ dans ce cas.
\qed

\paraunn{D\'emonstration du \ref{thmasymp}}
La fonction $\nuf(u)$ Žtant dŽfinie en
\eqref{defnufyu}, il est commode ici  de poser
$$\nuf(u)=\lk(u) + {\Delta_{f, y}(u)\over \log y}\cdot\eqdef{defDeltavy} $$
Notre dŽmonstration utilise de manire cruciale
la majoration issue du \ref{thmmajoration}
$$\psi_f^*(x,y)\ll F(1,y)\lambda_\kappa(u)\qquad
(1\leqslant u\ll(\log y)/\log_2y)
\eqdef{majpsi*}$$
valable uniformŽment pour $f\in\L_\kappa(A,C,\eta)$.
Nous exploitons ˆ nouveau l'Žquation
fonctionnelle \eqref{eqfoncnuf}, rŽŽcrite sous la
forme
$$u \, \mu_\kappa(u)\Deltaf(u)
= \kappa \int_{u-1}^u \mu_\kappa(v+1) \Deltaf(v) \d v
-\int_u^\infty \mu_\kappa(v) \d\bigg\{{R_f^*(v;y)\over F(1, y)}\bigg\},
\eqdef{eqfoncDeltaf}$$
mais nous avons ˆ prŽsent besoin d'une estimation bilatŽrale du terme d'erreur.

Dans toute la suite de la dŽmonstration, nous utilisons la notation
$$V_f(x,y):=\xk+{Z_1(y;f)\over u}+u\xk\log
\bigg(1+ {Z_2(y;f)\xk\over \log
y}\bigg).\eqdef{defVfxy} $$

\Propl{rkM}{Soient $A>0$, $C>0$, $\kappa>0$ et $\eta\in\,]0, \dm[$.
Il existe une constante positive $c_9$
telle que l'on ait
$$\bigg|\int_u^\infty \mu_\kappa(v)
\d\bigg\{{R_f^*(v;y)\over F(1, y)}\bigg\}\bigg|
\ll\lk(u)V_f(x,y)
\eqdef{majrkM}$$
uniform\'ement pour
$$f\in \MM_\kappa(A, C, \eta),
\quad
y\ge 3
\quad\hbox{et}\quad
1\le u\le c_9(\log y)/\log_2y.
\eqdef{domainerkM}$$
}
\goodbreak
\nid Conservons les notations \eqref{defIj} pour
les intŽgrales $I_j$ $(1\leqslant j\leqslant 3)$
dont la somme est
Žgale ˆ l'intŽgrale de \eqref{majrkM}.
Les assertions ˆ  Žtablir rŽsultent donc de la
validitŽ, dans le domaine \eqref{domainerkM}, des
estimations
$$I_j\ll \lk(u)V_f(x,y)\qquad (1\le j\le 3).\eqdef{restej} $$
\par
Le cas $j=1$ est presque immŽdiat. Il rŽsulte en
effet de la reprŽsentation \eqref{represintI1} et
du calcul qui la suit
que l'on a
$$I_1 \ll{\psi_f^*(x/y,y)\over uF(1,y)}\ll \lk(u)\xk,\eqdef{majI1M}$$
o l'on a utilisŽ \eqref{rkasymp}, \eqref{majpsi*} et \eqref{comploclam}.
\par
Le traitement de $I_2$ est analogue mais plus
dŽlicat. Nous observons d'abord qu'il rŽsulte
immŽdiatement de la
dŽfinition de $E_2(x,y;f)$ que
$$I_2\ll{1\over uF(1,y)}\sum_{\nu\ge 1}\sum_{p\le
y}{f(p)f(p^\nu)\log p\over
p^{\nu+1}}\psi_{f}^*(x/p^{\nu+1},y). $$
\par
DŽsignons par $I_{21}$ la sous-somme correspondant ˆ $\nu=1$
et par $I_{22}$ la sous-somme complŽmentaire.
Si $u\le 3$, nous avons trivialement
$$I_{21}\ll\sum_{p\le y}{f(p)^2\log p\over p^2}\ll \lk(u)V_f(x,y). $$
   Si $u>3$, nous dŽduisons de \eqref{majpsi*} et \eqref{comploclam} que, notant
$u_p:=(\log p)/\log y$,
$$\eqalign{I_{21}&\ll {1\over u}\sum_{p\le
y}{f(p)^2\log p\over
p^2}\lk(u-2u_p)\ll{\lk(u)\over u}\sum_{p\le
y}{f(p)^2\log
p\over p^{2\ak}}\cr
&\ll\lk(u)V_f(x,y),\cr}$$
d'aprs \eqref{majfp2alpha}.
\par
Nous avons ensuite, par \eqref{majfp},
$$I_{22}\ll {1\over uF(1,y)}\sum_{\nu\ge
2}\sum_{p\le y}{f(p^\nu)\over
p^{\nu}}\psi_{f}^*(x/p^{\nu+1},y).$$
Majorons $\psi_{f}^*(x/p^{\nu+1},y)$ par $$\ll
F(1,y)\lk(u-(\nu+1)u_p)\ll\lambda_\kappa(u)p^{\eta\nu}$$
lorsque $p^\nu\le y^{u-1/2}$ et, trivialement,
par $\ll F(1,y)$ dans le cas contraire. Il vient
$$\eqalign{I_{22}&\ll {\lk(u)\over u}\sum_{\nu\ge
2}\sum_{p\le y}{f(p^\nu)\over
p^{\nu(1-\eta)}}+\sum_{\nu\ge
2}\sum_{\di{p\le
y}{p^\nu>x/\sqrt{y}}}{f(p^\nu)\over
p^{\nu(1-\eta)}x^{\eta/2}}\cr
&\ll{\lk(u)\over u}+x^{-\eta/2}\ll\lk(u),\cr}
$$
puisque $x^{-\eta/2}\ll \exp\{-\eta
u^2/(2c_9)\}\ll\lk(u)$ dans le domaine
\eqref{domainerkM}.
\par
Cela Žtablit bien \eqref{restej} pour $j=2$.
\par\goodbreak
Le cas $j=3$ est relativement plus simple. Nous avons
$$\eqalign{
|I_3|
& \leqslant \mu_\kappa(u)\sum_{\di{p\le y}{\nu\ge 2}}
{f(p^\nu)\log p^\nu\over p^\nu}
\psi_f^*(x/p^\nu,y)
\cr}\eqdef{majI4M}
$$
et nous obtenons \eqref{restej} avec $j=3$ comme
prŽcŽdemment en employant \eqref{majpsi*} ou la
majoration triviale
selon que l'on a ou non $p^\nu\le
\sqrt{x}$.
\par
Cela ach\`eve la d\'emonstration.
\qed
Nous dŽduisons immŽdiatement de \eqref{majrkM} et
de \eqref{eqfoncDeltaf} une inŽquation
fonctionnelle pour $\Deltaf(u)$.
\Propl{inegDeltaM}
{Soient $A>0$, $C>0$, $\kappa>0$ et $\eta\in\,]0, \dm[$.
Il existe deux constantes positives
$C_0$ et $c_1$
telles que l'on ait
$$\eqalign{u \mu_\kappa(u)|\Delta_{f, y}(u)|
& \leqslant \kappa \int_{u-1}^u \mu_\kappa(v+1) |\Delta_{f, y}(v)| \d v
+ C_0\lk(u)V_f(x,y)
\cr}
\eqdef{iterationM}$$
pour
$f\in \MM_\kappa(A, C, \eta)$,
$y\ge 2$ et $1\le u\le c_1(\log y)/\log_2y$.
}

Le rŽsultat suivant fournit une estimation
prŽliminaire permettant d'initialiser le
traitement rŽcursif final.

\Propl{psi*u<u0}
{Soient $A>0$, $C>0$, $\kappa>0$,
$\eta\in\,]0,\dm[$, $u_0>1$. Posons
$\kappa_1:=\min(1,\kappa)$ et dŽsignons par~$c_1$
la constante apparaissant au \ref{inegDeltaM}. Sous les conditions
$$y\ge 2, \quad 0<u\le u_0,\eqdef{condmajinit}$$
la formule asymptotique
$$
\psi_f^*(x, y)
= F(1, y)
\bigg\{\lk(u)
+O\bigg({1+u^{\kappa_1-1}Z_1(y;f)\over \log y}\bigg)\bigg\}
\eqdef{valini}$$
est valable uniform\'ement pour
$f\in \MM_\kappa(A, C, \eta)$.
}

\nid
   Nous dŽduisons de \eqref{psifx} et \eqref{asympF1yM}
que, pour $2\le x\le y$,
$$\eqalign{\psi_f(x, y)
& = \psi_f(x)
= F(1, x) j_\kappa(1)
\bigg\{1+O\bigg({Z_1(x;f)\over \log x}\bigg)\bigg\}
\cr
& = {C_\kappa(f)\over \Gamma(\kappa+1)}
(u\log y)^\kappa
\bigg\{1+O\bigg({Z_1(y;f)\over u\log y}\bigg)\bigg\}.
\cr}$$
Cela implique bien \eqref{valini} pour $0<u\le 1$,
gr\^ace \`a une nouvelle application de~\eqref{asympF1yM}.
Nous avons donc obtenu, avec la notation
\eqref{defDeltavy}, l'existence d'une constante
$C_0=C_0(A,C,\eta)$ vŽrifiant
$$|\Delta_{f, y}(u)|
\leqslant C_0\{1+u^{\kappa-1}Z_1(y;f)\}
\eqdef{psi*u<1}$$
   pour tous $f\in \MM_\kappa(A, C, \eta)$,
$y\ge 2$ et $0<u\le 1$.
\par\goodbreak
Posons alors $\Delta_{f,y}^*(u):=\sup_{1/2\le
v\le u}|\Delta_{f,y}(v)|$. Nous allons montrer
par rŽcurrence sur
l'entier $j$, $1\leqslant j\leqslant  2u_0$, que
$$\Delta_{f,y}^*(u)\ll Z_1(y;f)\qquad (1\le
u\le  \dm
j+1).\eqdef{Deltafj}$$ D'aprs \eqref{psi*u<1}, nous avons d'abord
$$ \sup_{1/2\leqslant v\leqslant
1}|\Delta_{f,y}(v)|\ll
Z_1(y;f).\eqdef{Deltaf1/2-1}$$
   Ensuite, par
\eqref{psi*u<1} et le
\ref{inegDeltaM}, nous voyons qu'il existe une
constante $C_1$ telle~que, pour tous $u$, $w$
tels que $1\leqslant
w\leqslant u\leqslant
\ft32$,
$$\eqalign{w\mu_\kappa(w)|\Delta_{f,y}(w)|&\leqslant
\kappa
C_0\int_{w-1}^{w-1/2}\mu_\kappa(v+1)\{1+v^{\kappa-1}Z_1(y;f)\}\d
v\cr
&\quad +\Delta_{f,y}^*(u)\kappa\int_{w-1/2}^w\mu_\kappa(v+1)\d
v+C_1Z_1(y;f).\cr}
$$
Comme l'Žquation fonctionnelle \eqref{eqfoncmu} implique
$$0\le \kappa\int_{w-1/2}^w\mu_\kappa(v+1)\d v\le
\kappa\int_{w-1}^w\mu_\kappa(v+1)\d
v=w\mu_\kappa(w)-1\quad (w\ge 1),$$
il suit, pour une constante convenable $C_2$,
$$w\mu_\kappa(w)|\Deltaf(w)|\leqslant
\{w\mu_\kappa(w)-1\}\Delta_{f,y}^*(u)+
C_2Z_1(y;f). $$
En divisant par $w\mu_\kappa(w)$ et prenant le
supremum du membre de gauche lorsque $w\in[1,u]$,
nous obtenons
immŽdiatement que, si
$\Deltaf^*(u)>C_2Z_1(y;f)$, alors
$\Deltaf^*(u)$ est Žgal au membre de gauche de
\eqref{Deltaf1/2-1}.\par
Il s'ensuit que la majoration \eqref{Deltafj} est vŽrifiŽe pour $j=1$.
\par
Supposons-la satisfaite
pour
$j\ge 1$.  Sous l'hypothse $1+\dm j\le u\le
\ft12+\dm (j+1)$, nous avons alors, d'aprs le
\ref{inegDeltaM}, pour tout $w$ de $[1+\dm
j,\ft32+\dm j]$
$$\eqalign{w\mu_\kappa(w)|\Delta_{f,y}(w)|&\leqslant\kappa
\int_{w-1}^w\mu_\kappa(v+1)\Delta_{f,y}^*(u)\d
v+C_1Z_1(y;f)\cr
&=\{w\mu_\kappa(w)-1\}\Delta_{f,y}^*(u)+C_1Z_1(y;f),\cr}
$$
et donc, en faisant tendre $|\Deltaf(w)|$ vers $\Deltaf^*(u)$,
$$\Delta_{f,y}^*(u)\leqslant C_1Z_1(y;f). $$
Cela Žtablit que la propriŽtŽ annoncŽe est bien
inductive et achve la dŽmonstration de la
formule \eqref{valini}.
\qed

\par\goodbreak
Nous sommes ˆ prŽsent en mesure de terminer la preuve du \ref{thmasymp}.
\par

Soit $u_0>1$ tel que $\inf_{v\ge u_0}v\mu_\kappa(v)\geqslant \dm$.
D'apr\`es  le \ref{psi*u<u0},
il existe une constante
$C_3=C_3(A, C, \kappa, \eta, u_0)>8C_0$ telle que
$$|\Delta_{f, y}(u)|
\leqslant C_3\lk(u)uV_f(x,y)
\qquad(1\le u\le u_0, \, y\ge 2).$$
\par
Pour $y\ge 2$ fix\'e, posons
$$u_2=u_2(y)
:= \inf\big\{v\ge 1 :
|\Delta_{f, y}(v)|
>C_3v\lambda_\kappa(v)V_f(y^v,y)\big\}$$
et notons $x_2:=y^{u_2}$.
Nous avons trivialement  $u_2\ge u_0$.
Si $u_2\le c_1(\log y)/\log_2y$, nous dŽduisons
de \eqref{eqfoncmu}, \eqref{iterationM} et
\eqref{pdtscal}
que
$$\eqalign{
C_3u_2^2&\mu_\kappa(u_2)\lambda_\kappa(u_2)V_f(x_2,y)
   \leqslant u_2\mu_\kappa(u_2)|\Delta_{f, y}(u_2)|
\cr
& \leqslant \kappa
\int_{u_2-1}^{u_2}  \mu_\kappa(v+1) |\Delta_{f, y}(v)| \d v
+ C_0\lambda_\kappa(u_2)V_f(x_2,y)
\cr
&\leqslant  \Big\{\kappa C_3
\int_{u_2-1}^{u_2}  \mu_\kappa(v+1) v\lambda_\kappa(v) \d v
+ C_0\lambda_\kappa(u_2)\Big\}V_f(x_2,y).
\cr}\eqdef{contrad}$$
Comme la fonction $v\mapsto
\mu_\kappa(v+1)\lambda_\kappa(v)$ est
dŽcroissante, nous avons
$$\kappa\int_{u_2-1}^{u_2-1/2}\mu_\kappa(v+1)\lambda_\kappa(v)\d
v\ge \dm \kappa
\int_{u_2-1}^{u_2}\mu_\kappa(v+1)\lambda_\kappa(v)\d
v=\dm u_2\mu_\kappa(u_2)\lambda_\kappa(u_2), $$
d'o
$$\eqalign{\kappa\int_{u_2-1}^{u_2}&
\mu_\kappa(v+1) v\lambda_\kappa(v) \d v\cr&\leqslant
\kappa(u_2-\dm)\int_{u_2-1}^{u_2-1/2}\mu_\kappa(v+1)\lambda_\kappa(v)\d
v+\kappa
u_2\int_{u_2-1/2}^{u_2}\mu_\kappa(v+1)\lambda_\kappa(v)\d v\cr
&\leqslant (u_2-\ft14)u_2\mu_\kappa(u_2)\lambda_\kappa(u_2).\cr} $$
En reportant dans \eqref{contrad}, il suit, aprs
division par $\lambda_\kappa(u_2)V_f(x_2,y)$,
$$\ft18 C_3\le \ft14C_3u_2\mu_\kappa(u_2)\leqslant C_0,$$
ce qui contredit la dŽfinition de $C_3$.
Cela ach\`eve la d\'emonstration.
\bigskip

\paraunn{D\'emonstration du \ref{thmHF}}
\paradeuxn{ƒtude d'un cas particulier}

Nous \'etablissons trois propositions pr\'eliminaires,
qui sont respectivement des modifications
mineures des propositions 4.2, 4.3 et 4.4 de
\citer{TW03}.
Dans toute la suite les num\'eros de
formules $(a{\cdot}b)$ de \citer{TW03}
sont identifi\'{e}s sous la forme \eqnumm{a}{b}.

Comme dans \citer{TW03}, nous dŽsignons par
$\M_\kappa^*(A;R)$ la classe des fonctions
   exponentiellement multiplicatives,
\ie\ telles que
$$f(p^\nu)=f(p)^\nu/\nu !
\quad(p\ge2, \nu\ge 1),$$
satisfaisant \eqref{regfp}.
Nous conservons Žgalement les notations
$$\Psi_f(x, y)
:= \sum_{n\in S(x, y)} f(n),
\qquad
R_b(y)
:=\int_{3/2}^y {\d t\over tR(t^b)}.$$
\par
Par souci de concision, dans toute la suite de ce
paragraphe, nous renvoyons librement ˆ
\citer{TW03} et nous nous
limitons ˆ signaler les modifications n\'ecessaires
des ŽnoncŽs ou des dŽmonstrations de
\citer{TW03}.

\Propl{L1}
{Sous la condition suppl\'ementaire $x\ge y$,
la proposition 4.2 de \citer{TW03} reste valable
en rempla\c cant $E_0(x,y)$ par
$$E_0'(x,y)
:= {\xi_\kappa(u)\over \log y}\qquad (x\ge y\ge 3).$$}

\nid
La preuve de la proposition 4.2 de \citer{TW03}
reste valable, moyennant les changements
suivants.\par
\noi\bull{\it Page 150, ligne 3 de \citer{TW03}.}
D'apr\`es les estimations \eqnumm{3}{35} et \eqnumm{4}{12}, on a
$$\eqalign{W_{j2}
& \ll {\Psi_f(3x/2,y)\over x}
\ll (3/2)^{\alpha} \varrho_\kappa(u)(\log y)^{\kappa-1}
\cr
& \ll \lk(u)(\log y)^{\kappa} {\xi_\kappa(u)\over \log y}
\quad (j=1,2),
\cr}$$
d'o\`u
$${\xi_\kappa(u)^2\over (\log y)^2}W_{12}+W_{22}
\ll\lk(u)(\log y)^{\kappa} E'_0(x,y).$$

\bull{\it
Page 150, ligne $-2$ de \citer{TW03}.}
Les relations \eqnumm{4}{18} et \eqnumm{4}{12}  impliquent
$$\eqalign{{\xi_\kappa(u)^2\over (\log y)^2}W_{11}+W_{21}
& \ll \varrho_\kappa(u)(\log
y)^{\kappa-1}\bigg(1+{\xi_\kappa(u)^2\log_2(3x)\over (\log y)^2}\bigg)
\cr\noalign{\smallskip}
& \ll\lk(u)(\log y)^{\kappa} E'_0(x,y).
\cr}$$

\bull{\it Page 151, lignes 1-5 de \citer{TW03}.} On a
$$W_{j3}\ll \cases{
\varrho_\kappa(u) (\log y)^{\kappa-1} \log_2(3x) & si $j=1$,
\cr\noalign{\smallskip}
\varrho_\kappa(u) (\log y)^{\kappa-1}            & si $j=2$.
\cr}$$
Il s'ensuit, comme pr\'ec\'edemment, que
$$\eqalign{{\xi_\kappa(u)^2\over (\log y)^2}W_{13}+W_{23}
& \ll \varrho_\kappa(u)(\log
y)^{\kappa-1}\bigg({\xi_\kappa(u)^2\log_2(3x)\over (\log y)^2}+1\bigg)
\cr\noalign{\smallskip}
& \ll\lk(u)(\log y)^{\kappa} E'_0(x,y).
\cr}$$

\bull{\it Page 151, lignes $-7$ \`a $-1$ de \citer{TW03}.}
Estimons pr\'ecis\'ement les trois termes apparaissant
\`a la ligne $-4$ de la page 151 de \citer{TW03}.

En utilisant \eqnumm{3}{35} avec $(t, x) = (x, Y_\varepsilon^*)$,
nous avons d'abord
$$\normalbaselineskip15pt\eqalign{{x\Psi_f(xY_\varepsilon^*,y)\over
(xY_\varepsilon^*)^2}
& \ll {x\over (xY_\varepsilon^*)^2}
x^{\alpha} Y_\varepsilon^*
\varrho_\kappa\bigg({\log Y_\varepsilon^*\over
\log y}\bigg) (\log y)^{\kappa-1}
\cr\noalign{\vskip 1mm}
& \ll (Y_\varepsilon^*)^{-1}
\varrho_\kappa(u) (\log y)^{\kappa-1}
\cr\noalign{\vskip 2mm}
& \ll \lk(u) (\log y)^{\kappa} E'_0(x,y).
\cr}$$
Pour la seconde estimation,
nous avons utilis\'e l'in\'egalit\'e $\alpha\le 1$ et le fait que
$v\mapsto \varrho_\kappa(v)$ est positive et d\'ecroissante ($\kappa<1$)
ou unimodale ($\kappa\ge 1$).
La troisi\`eme estimation r\'esulte trivialement
de $Y_\varepsilon^*\geqslant 1$ et \eqnumm{4}{12}.
\par
Semblablement, nous pouvons \'ecrire
$$\eqalign{
x\int_{Y_\varepsilon^*}^{xY_\varepsilon^*} {\Psi_f(t,y)\over t^3} \d t
& = {x\over (Y_\varepsilon^*)^2}
\int_{1}^{x}{\Psi_f(vY_\varepsilon^*,y)\over v^3} \d v
\cr
& \ll {x\over (Y_\varepsilon^*)^2} \int_{1}^{x}
{v^{\alpha} Y_\varepsilon^*
\varrho_\kappa(\log Y_\varepsilon^*/\log y)
(\log y)^{\kappa-1}\over v^3} \d v
\cr\noalign{\vskip 1mm}
& \ll \varrho_\kappa(u) (\log y)^{\kappa-1}
\ll \lk(u) (\log y)^{\kappa} E'_0(x,y),
\cr}$$
o\`u la premi\`ere majoration r\'esulte
d'une application de \eqnumm{3}{35}
avec $(t, x) = (v, Y_\varepsilon^*)$.
En utilisant \eqnumm{3}{26} avec $x=t$, on a de plus
$$\eqalign{
x\int_x^{Y_\varepsilon^*} {\Psi_f(t,y)\over t^3}\d t
& \ll x(\log y)^{\kappa-1} \int_x^{Y_\varepsilon^*}
{\varrho_\kappa(\log t/\log y)
\over t^2} \d t
\cr
& \ll \varrho_\kappa(u) (\log y)^{\kappa-1}
\ll \lk(u) (\log y)^{\kappa}
E'_0(x,y).
\cr}$$
\par
Nous obtenons ainsi
$$x\int_x^{xY_\varepsilon^*}{\Psi_f(t,y)\over t^3}\d t
\ll \lk(u) (\log y)^{\kappa} E'_0(x,y).$$
\par
Il reste \`a examiner le troisi\`eme terme.
D'apr\`es  \eqnumm{3}{25} et \eqnumm{4}{12}, nous avons
$$\eqalign{{\rm e}^{-u\xi_\kappa(u)} (\log
y)^{\kappa-1}{\sqrt{u+1}\log y\over Y_\varepsilon^*}
& \asymp \varrho_\kappa(u)
(\log y)^{\kappa-1}
\wh\varrho(-\xi_\kappa(u))^{-\kappa}{u\log y\over Y_\varepsilon^*}
\cr
& \ll \varrho_\kappa(u) (\log y)^{\kappa-1}
    \ll \lk(u) (\log y)^{\kappa} E'_0(x,y).
\cr}$$
Cela ach\`eve la preuve du \ref{L1}.
\qed

\Propl{L2}{
Sous la condition suppl\'ementaire $x\ge y$,
la proposition 4.3 de~\citer{TW03} reste valable en rempla\c cant
$E_\kappa^*(x,y)$ par
$$E_\kappa^\dagger(x,y)
:= {u \xi_\kappa(u)\over R(y^b)} + {\xi_\kappa(u) R_b(y)\over\log y}.
\eqdef{1}$$}

\nid
Nous allons \'etablir, avec les notations \eqnumm4{33} et \eqnumm4{34},
les estimations
asymptotiques
$$I_1
= \lk(u)
\bigg\{1 + O\bigg({1\over (\log y)^{1+\kappa}}\bigg)\bigg\},
\eqdef{2}$$
$$I_2\ll \lk(u) {\varepsilon(u,y)\over \sqrt u}.
\eqdef{3}$$
Le r\'esultat annonc\'e d\'ecoule imm\'ediatement
de ces estimations au vu de
\eqnumm{4}{32} et~\eqnumm{4}{20}.

\bull {\it Preuve de \eqref{2}.}
Il suffit de d\'emontrer que l'on a
$$\eqalign{I_1'
& := \int_{|t|\geqslant {1\over 2}\log y}
\widehat j_\kappa(w){\rm e}^{uw}\d t
\ll {\lk(u)\over (\log y)^{1+\kappa}},
\cr}$$
o\`u, par convention, $w=-\xi_\kappa(u)+it$.

Soit $T := \max\{{1\over 2}\log y, \, 1+u\xi_k(u)\}$.
D\'esignons par $I_{11}'$ et $I_{12}'$
les contributions respectives \`a $I_1'$ des domaines
${1\over 2}\log y\le |t|\leqslant T$ et $|t|> T$.

Pour estimer $I'_{11}$, nous pouvons supposer
$T = 1+u\xi_k(u)\geqslant {1\over 2}\log y$
car le domaine d'int\'egration est vide dans le cas contraire.
Nous avons alors, d'apr\`es
\eqnumm{4}{24}
et le lemme 4.10 de
\citer{Sm91},
$$\eqalign{I_{11}'
& \ll \e^{-u\xi_\kappa(u)}\wh\varrho(-\xi_\kappa(u))^{\kappa}
\int_{{1\over 2}\log y\le |t|\leqslant T}
\e^{-u/(\xi_\kappa(u)^2+\pi^2)}  {\d t\over t}
\cr
& \ll \varrho_\kappa(u)\,\sqrt{u}\,\e^{-u/(\xi_\kappa(u)^2+\pi^2)} \log T
\ll {\lk(u)\over (\log y)^{1+\kappa}},
\cr}$$
en vertu de \eqnumm{3}{25} et \eqnumm{4}{12}.\par\goodbreak

D'autre part, on peut \'ecrire, gr\^ace \`a \eqnumm{4}{35},
$$I_{12}'
= \e^{-\gamma\kappa-u\xi_\kappa(u)}
\int_{|t|\geqslant T}
\bigg\{{1\over w^{1+\kappa}}
+ O\bigg({1+u\xi_\kappa(u)\over w^{2+\kappa}}\bigg)\bigg\}
\e^{iut} \d t.$$
La contribution du terme d'erreur \`a la derni\`ere
int\'egrale est
$$\ll u\xi_\kappa(u)/T^{1+\kappa}
\ll u\xi_\kappa(u)/(\log y)^{1+\kappa}.$$
Celle du terme principal $\e^{iut}/w^{1+\kappa}$
peut \^etre estim\'ee par
la seconde formule de la moyenne: nous avons
$$\eqalign{\int_{|t|\geqslant T} {\e^{iut}\over w^{1+\kappa}} \d t
& = \int_{|t|\geqslant T}
{\e^{iut}\over \{-\xi_\kappa(u) + it\}^{1+\kappa}} \d t
\cr
& = {1\over i^{1+\kappa}} \int_{|t|\geqslant T} {\e^{iut}\over |t|^{1+\kappa}}
\bigg\{1+O\bigg({\xi_\kappa(u)\over |t|}\bigg)\bigg\} \d t
\cr
& \ll {\xi_\kappa(u)\over uT^{1+\kappa}}
\ll {\xi_\kappa(u)\over u(\log y)^{1+\kappa}}.
\cr}$$
Nous avons donc \'etabli que
$$I_{12}'
\ll \e^{-u\xi_\kappa(u)} {u\xi_\kappa(u)\over (\log y)^{1+\kappa}}
\ll {\lk(u)\over (\log y)^{1+\kappa}}\cdot$$

\bull {\it Preuve de \eqref{3}.} Compte tenu de
\eqnumm{4}{38}, il suffit de d\'emontrer que l'on a
$$I_{22}
:= {1\over 2\pi}
\int_{T_2\le |t|\leqslant {1\over 2}\log y}
D_y\Big(1+{w\over \log y}\Big) \widehat j_\kappa(w) {\rm e}^{uw} \d t
\ll \lk(u) {\varepsilon(u,y)\over \sqrt u}.$$
Nous allons proc\'eder comme dans \citer{TW03} dans le cas  $\kappa=1$,
$u>1$.\note{Voir
\citer{TW03} \page156, lignes 6-19.}

Similairement \`a \eqnumm{4}{40},
nous pouvons \'ecrire, en vertu de la formule asymptotique~\eqnumm{4}{35},
$$I_{22}=I_{221}+I_{222}$$
avec
$$I_{221}:={-{\rm e}^{-u\xi_1(u)-\gamma}\over \pi}
\int_{T_2}^{{1\over 2}\log y}
\re \bigg(D_y\bigg(1+{w\over \log y}\bigg){\rm e}^{uit}\bigg) {\d
t\over t^{1+\kappa}},$$
puisque l'on a  $D_y(\overline s) = \overline {D_y(s)}$ pour $|s-1|<1$, et
$$\eqalign{I_{222}
& \ll { T_2 \varepsilon(u,y)\over \e^{u\xi_\kappa(u)}}
\int_{T_2}^{\log y}{t^{-1-\kappa}\d t\over 1+\varepsilon(u,y)t}
\cr
& \ll {\rm e}^{-u\xi_\kappa(u)} T_2 \varepsilon(u,y)
   \ll \lk(u) {\varepsilon(u,y)\over \sqrt u}\cdot
\cr}$$
Une int\'egration par parties fournit ensuite,
compte tenu de \eqnumm{4}{28},
\eqnumm{4}{29}
et~\eqnumm{4}{41},
$$\eqalign{I_{221}
& \ll {\rm e}^{-u\xi_\kappa(u)}
\bigg\{{\varepsilon(u,y)\over T_2^\kappa}
+
\!\!
\int_{T_2}^{{1\over 2}\log y}
\!\!
\bigg({|D_y'(1+w/\log y)|\over
t^{1+\kappa}\log y}
+ {|D_y(1+{w/\log y})|\over t^{2+\kappa}}\bigg) \d t\bigg\}
\cr
& \ll {\rm e}^{-u\xi_\kappa(u)}
\int_{T_2}^{{1\over 2}\log y}{\varepsilon(u,y)\over t^{1+\kappa}} \d t
   \ll {\rm e}^{-u\xi_\kappa(u)} \varepsilon(u,y)
   \ll \lk(u) {\varepsilon(u,y)\over \sqrt u}\cdot
\cr}$$
Cela ach\`eve la d\'emonstration.
\qed

Les Lemmes \ref{sL1} et \ref{sL2}
impliquent imm\'ediatement le r\'esultat suivant.

\Propl{L3}
{Sous la condition suppl\'ementaire $x\ge y$, la proposition 4.4
de~\citer{TW03} reste
valable avec $E_\kappa^\dagger(x,y)$ \`a la place de
$E_\kappa^*(x, y)$.}

Le rŽsultat suivant fournit la conclusion du
\ref{thmHF} dans le cas particulier
$f\in\M_\kappa^*(A;R)$. Nous verrons au
paragraphe suivant que le cas gŽnŽral s'en dŽduit
aisŽment gr‰ce ˆ un argument de convolution.

\Propl{L4}
{Soient $A>0$, $\varepsilon>0$,
$\kappa>0$, $b\in\,]0,\dm]$,
$R\in \RR(b;\kappa)$.
On a
$$\psi_f^*(x, y)
=F(1,y)\lambda_\kappa(u)\big\{1+O\big(E_\kappa^\ddagger(x, y)\big)\big\},
\eqdef{asymf*}$$
uniform\'ement pour $f\in\M_\kappa^*(A;R)$,
$y\ge 2$ et $\sqrt{y}\leqslant x\le Y_\varepsilon^*$,
o\`u
$$E_\kappa^\ddagger(x, y)
:= {u\xi_\kappa(u)\over R(y^{b/2})}
+ {\xi_\kappa(u) R_b(y)\over \log y}.$$}

\nid
Puisque
$\e^{-\gamma \kappa}/\Gamma(\kappa+1)\leqslant j_\kappa(u)\leqslant 1$,
le \ref{L3} implique que l'on a, uniform\'ement
pour $f\in \M_\kappa^*(A;R)$ et $2\le y\le x\le Y_\varepsilon$,
$$\eqalign{\psi_f(x, y)
& = \sum_{P(n)\leqslant y} {f(n)\over n} - \sum_{\di{
n>x}{P(n)\leqslant y}} {f(n)\over n}
\cr
& = F(1,y) j_\kappa(u) \big\{1 + O\big(\lambda_\kappa(u)
E_\kappa^\dagger(x,y)\big)\big\}.
\cr}
\eqdef{5}$$
Nous allons maintenant \'etablir que,
sous l'hypoth\`ese $\sqrt{y}\leqslant x\le y$,
cette formule reste valable en rempla\c cant $E_\kappa^\dagger(x, y)$ par
$E_\kappa^\ddagger(x, y)$.
\`A cette fin, nous observons que le cas particulier $y=x$ de \eqref{5},
\eqnumm{2}{7} et la formule
$$j_\kappa(u)
= {\e^{-\gamma \kappa}\over \Gamma(\kappa+1)} u^\kappa
\qquad(0\le u\le 1)$$
impliquent, pour $\sqrt{y}\leqslant x\le y$,
$$\eqalign{F(1, x) j_\kappa(1)
& = \e^{\gamma \kappa} C_\kappa(f) (\log x)^\kappa
\bigg\{1+O\bigg({1\over R(x^b)}\bigg)\bigg\}
{\e^{-\gamma \kappa}\over \Gamma(\kappa+1)}
\cr
& = F(1, y) j_\kappa(u) \bigg\{1+O\bigg({1\over R(y^{b/2})}\bigg)\bigg\}
\cr}$$
et
$$\eqalign{\lambda_\kappa(1) E_\kappa^\dagger(x, x)
& \ll {1\over R(x^{b})} + {R_b(x)\over \log x}
\ll \lambda_\kappa(u) \bigg({u\xi_\kappa(u)\over R(y^{b/2})}
+ {\xi_\kappa(u) R_b(y)\over \log y}\bigg).
\cr}$$
Cela fournit bien l'extension annonc\'ee.
Il est \`a noter que
$E_\kappa^\ddagger(x,y)$ est, \`a
$y$ fix\'e, une fonction croissante de $x$.
Ainsi on a
$$\psi_f(x, y)
= F(1,y) j_\kappa(u)
\big\{1 + O\big(\lk(u) E_\kappa^\ddagger(x, y)\big)\big\}$$
uniform\'ement pour $f\in \M_\kappa^*(A;R)$,
$y\ge 2$ et $\sqrt{y}\leqslant x\le Y_\varepsilon^*$.
Cela \'equivaut au r\'esultat annoncŽ.
\qed

\goodbreak\smallskip

\paradeuxn{Extension au cas gŽnŽral}

Maintenant consid\'erons une fonction quelconque
$f\in \M_\kappa(A, C, \eta; R)$.
On peut \'ecrire $f=f^**g$ avec
$f^*\in \M_\kappa^*(A;R)$, et
$$g(p^\nu)
:=\sum_{j+k=\nu}(-1)^jf(p^k)f(p)^j/j!
\qquad
(p\ge 2,\nu\ge 1).$$

Posons $z := \max\{x/y, \, \sqrt{y}\}$,
$u_m:=(\log m)/\log y$ $(m\ge 1)$ et
$$F^*(s, y):=\sum_{P(n)\leqslant y}f^*(n)/n^s.$$
En appliquant le \ref{L4} \`a $f^*$,\note{La
formule \eqref{asymf*} est bien applicable pour
Žvaluer
$\psi_{f^*}(x/m,y)$
   puisque $m\le z$ implique
$\log(x/m)/\log y\ge \dm.$} nous obtenons
$$\eqalign{\quad
\psi_f^*(x, y)
& = \sum_{\di{m\le x}{P(m)\leqslant y}}
{g(m)\over m} \psi_{f^*}^*\bigg({x\over m}, y\bigg)
+ F^*(1, y) \sum_{\di{m>x}{P(m)\leqslant y}} {g(m)\over m}
\cr
& = F^*(1, y)
\big\{
\Upsilon_1
+ O\big(E_\kappa^\ddagger(x, y)\Upsilon_1
+ \Upsilon_2\big)\big\},
\cr}
\eqdef{Decomppsi*}$$
avec
$$\eqalign{
\Upsilon_1
& := \sum_{\di{m\le z}{P(m)\leqslant y}}
{g(m)\over m} \lambda_\kappa(u-u_m),
\cr
\Upsilon_2
& := \sum_{\di{m>z}{P(m)\leqslant y}} {|g(m)|\over m},
\cr}\eqdef{Upsj}$$
o nous avons employŽ l'estimation triviale
$\psi_{f^*}^*(x/m, y)\!\leqslant\! F^*(1, y)$ pour \hbox{$z<m\le x$}.
\par\goodbreak
Posons
$$\delta_\kappa(u, u_m)
:= \lk(u) - \lambda_\kappa(u-u_m)
= \int_0^{u_m} \varrho_\kappa(u-v) \d v.$$
Il rŽsulte aisŽment de \eqref{comploclam} que
$$\delta_\kappa(u, u_m)
\ll \cases{
\displaystyle \lk(u) \xi_\kappa(u)
\dsp{\log m\over \log y}      & si $u_m\le 1/(2\xi_\kappa(u))$,
\cr\noalign{\medskip}
\lk(u) m^{1-\alpha}           & si $u_m\le u-\dm$.
\cr}
\eqdef{7}$$
\par
Nous majorons $\Upsilon_2$ par la m\'ethode de Rankin.
Lorsque $x/y\ge \sqrt{y}$, nous avons gr\^ace \`a
\eqnumm{4}{51}
$$\eqalign{\Upsilon_2
& \ll \sum_{\di{m>x/y}{ P(m)\leqslant y}}
{|g(m)|\over m}
\bigg({x\over my}\bigg)^{\alpha-1} {\log m\over \log (x/y)}
\cr
& \ll {\e^{-(u-1)\xi_\kappa(u)}\over (u-1)\log y}
\sum_{P(m)\leqslant y} {|g(m)|\log m\over m^{\alpha}}
\cr
& \ll {\e^{-u\xi_\kappa(u)} \xi_\kappa(u)\over \log y}
\bigg(1+{u^2\xk\over R(y^b)}\bigg)
\cr\noalign{\vskip 1mm}
& \ll \lk(u) {\xi_\kappa(u)\over \log y}
\bigg(1+{u^2\xk\over R(y^b)}\bigg).
\cr}
\eqdef{MajU21}$$
Similairement si $x/y<\sqrt{y}$, on a
$$\eqalign{\Upsilon_2
& \ll \sum_{\di{m>\sqrt{y}}{ P(m)\leqslant y}}
{|g(m)|\over m}
\bigg({\sqrt{y}\over m}\bigg)^{\alpha-1} {\log m\over \log y}
\cr
& \ll {y^{-(1-\alpha)/2}\over \log y}
\sum_{P(m)\leqslant y} {|g(m)|\log m\over m^{\alpha}}
\cr\noalign{\vskip 0,5mm}
& \ll \lk(u) {\xi_\kappa(u)\over \log y}
\bigg(1+{u^2\xk\over R(y^b)}\bigg),
\cr}
\eqdef{MajU22}$$
puisque $x/y<\sqrt{y}\Rightarrow u\le {3\over 2}$.

Il reste \`a \'evaluer le terme principal de \eqref{Decomppsi*}.
\`A cette fin, nous \'ecrivons
$$\Upsilon_1
= \lk(u) \sum_{m\ge 1} {g(m)\over m}
+ O(\lk(u)\Upsilon_2 + \Upsilon_3),
\eqdef{10}$$
o\`u $\Upsilon_2$ est d\'efinie par \eqref{Upsj} et
$$\Upsilon_3
:= \sum_{\di{m\le z}{P(m)\leqslant y}}
{|g(m)\delta_\kappa(u, u_m)|\over m}.$$

Posons $w := y^{1/(2\xk)}$.
D'apr\`es la premi\`ere majoration de \eqref{7} et \eqnumm{4}{51},
nous~avons
$$\eqalign{\sum_{\di{m\le w}{P(m)\leqslant y}}
{|g(m) |\delta_\kappa(u, u_m)|\over m}
& \ll \lk(u) {\xi_\kappa(u)\over \log y}
\sum_{m\le z}{|g(m)|\log m\over m}
\cr
& \ll \lk(u) {\xi_\kappa(u)\over \log y}
\bigg(1+{u^2\xk\over R(y^b)}\bigg).
\cr}$$
En utilisant l'in\'egalit\'e $1\le \log m/\log w$,
la seconde estimation de \eqref{7} et \eqnumm{4}{51},
nous obtenons de plus que
$$\eqalign{
\sum_{\scriptstyle w<m\le z\atop\scriptstyle P(m)\leqslant y}
{|g(m)\delta_\kappa(u, u_m)|\over m}
& \ll {\lk(u)\over \log w}
\sum_{P(m)\leqslant y} {|g(m)| \log m\over m^{\alpha}}
\cr
& \ll \lk(u) {\xi_\kappa(u)\over \log y}
\bigg(1+{u^2\xk\over R(y^b)}\bigg).
\cr}$$
En rassemblant ces estimations, nous d\'eduisons que
$$\Upsilon_3
\ll \lk(u) {\xi_\kappa(u)\over \log y}
\bigg(1+{u^2\xk\over R(y^b)}\bigg)
\eqdef{MajU3}$$
et finalement
$$\Upsilon_1
= \lk(u) \sum_{P(m)\leqslant y} {g(m)\over m}
+ O\bigg(\lk(u){\xk\over \log y}\bigg(1+{u^2\xk\over R(y^b)}\bigg)\bigg).
\eqdef{EvalU1}$$
Maintenant le \ref{thmHF} d\'ecoule donc de
\eqref{Decomppsi*},
\eqref{MajU21},
\eqref{MajU22},
\eqref{MajU3},
\eqref{EvalU1} et les faits que
$$F^*(1, y)\sum_{P(m)\leqslant y} {g(m)\over m}
=F(1, y),
\qquad
F^*(1, y)\asymp F(1, y).
$$
Cela ach\`eve la d\'emonstration.
\qed

\goodbreak

\vskip 10mm

\centerline{\twelvebf Bibliographie}
\bigskip
\eightpoint{\leftskip9mm\rightskip6mm

\bibtem{deBL64}
N.G. de Bruijn \& J.H. van Lint,
Incomplete sums of multiplicative fonction, I, II,
{\it Nederl. Akad. Wetensch. Proc. Ser. A} {\bf 67} (1966), 339--347;
348--359.

\bibtem{Gal72}
P. X. Gallagher,
Sieving by prime powers,
in: {\it Proceedings of the 1972 Number Theory Conference}
(Univ. Colorado, Boulder, Colo., 1972),  pp. 95--99.
Univ. Colorado, Boulder, Colo., 1972.

\bibtem{Gal73/74}
P. X. Gallagher,
Sieving by prime powers,
Collection of articles dedicated to Carl Ludwig Siegel
on the occasion of his seventy-fifth birthday, V.
{\it Acta Arith.} {\bf 24} (1973/74), 491--497.

\bibtem{G01}
G. Greaves,
{\it Sieves in number theory},
Ergebnisse der Mathematik und ihrer Grenzgebiete (3)
[Results in Mathematics and Related Areas (3)], {\bf 43}.
Springer-Verlag, Berlin, 2001. xii+304 pp.

\bibtem{G05}
G. Greaves,
Incomplete sums of non-negative multiplicative functions,
\AA, {\bf 118} (2005), n$^\circ$ 4, 337--357.

\bibtem{He86} D. Hensley,
    The convolution powers of the Dickman 
function, {\it J. London Math. Soc. \rm
(2) \bf 33} (1986), 395--406.

\bibtem{HR74}
H. Halberstam \& H.-E. Richert,
{\it Sieve methods},
London Mathematical Society Monographs, n$^\circ$ 4,
Academic Press [A subsidiary of Harcourt Brace Jovanovich, Publishers],
London-New York, 1974. xiv+364 pp.

\bibtem{HTW06}
G. Hanrot, G. Tenenbaum \& J. Wu,
Moyennes de certaines fonctions multiplicatives
sur les entiers friables, 2, \PLMS, ˆ para"tre.



\bibtem{HT93}
A. Hildebrand \& G. Tenenbaum,
On a class of difference differential equations arising in number theory,
{\it J. d'Analyse \bf 61} (1993),  145--179.

\bibtem{Jo71} J. Johnsen, On the large sieve 
method in $GF(q, x)$, {\it Mathematika \bf18} 
(1971), 172--184.

\bibtem{Moto79}
Y. Motohashi,
A note on the large sieve. III,
{\it Proc. Japan Acad. Ser. A Math. Sci.} {\bf 55} (1979),
n$\circ$ 3, 92--94.

\bibtem{R82}
D. A. Rawsthorne,
Selberg's sieve estimate with a one sided hypothesis,
\AA\ {\bf 61} (1982), 282--289.

\bibtem{Sel77}
A. Selberg,
Remarks on multiplicative functions,
in: {\it Number theory day}
(Proc. Conf., Rockefeller Univ., New York, 1976),  pp. 232--241.
Lecture Notes in Math., Vol. 626, Springer, Berlin, 1977.

\bibtem{Sm91}
H. Smida,
Sur les puissances de convolution de la fonction de Dickman,
\AA\ {\bf 59}, n$^{\circ}$\thinspace2  (1991), 124--143.


\bibtem{So01}
J. M. Song,
Sums of nonnegative multiplicative functions over integers free of large
prime factors I,
\AA, \bf97\rm, 4 (2001), 329--351.

\bibtem{Te90}
G. Tenenbaum,
Sur une question d'Erd\H os et Schinzel,
in: {\it A tribute to Paul Erd\H os},
405--443, Cambridge Univ. Press, Cambridge, 1990.

\bibtem{Te95}
G. Tenenbaum,
{\it Introduction \`a la th\'eorie analytique et
probabiliste des nombres},
Cours sp\'ecialis\'es, n$^{\circ}$\thinspace1,
Soci\'et\'e Math\'ematique de France (1995), xv + 457 pp.

\bibtem{Te01}
G. Tenenbaum,
Note on a paper by Joung Min Song,
\AA\ \bf97\rm, 4 (2001), 353--360.

\bibtem{TW03}
G. Tenenbaum \& J. Wu,
Moyennes de certaines fonctions multiplicatives sur les entiers friables,
{\it J. reine angew. Math.} {\bf 564} (2003), 119--166.
\par
}

\bigskip
\nobreak
\medskip

\vbox{\sevenrm\baselineskip9pt\obeylines
Institut \'Elie Cartan, UMR 7502
Universit\'e Henri Poincar\'e--Nancy 1
BP 239\par
54506 Vand\oe uvre-l\`es-Nancy
France
\medskip
\seventt gerald.tenenbaum@iecn.u-nancy.fr
wujie@iecn.u-nancy.fr}

\bye